\documentclass[11pt,openany,leqno]{smfart}
\usepackage[english, francais]{babel}
\usepackage{smfthm}
\usepackage{amsmath,amsthm,amsfonts,amssymb,amscd,url}

\usepackage{graphics}
\usepackage{float} 
\input{xy} \xyoption{all}
\usepackage[latin1]{inputenc}
\usepackage{enumerate}
\usepackage{color}
\usepackage{epsfig} 
\input{xy} \xyoption{all}

\begin{document}

\begin{otherlanguage}{english}
\def\qand{\quad \text{and}\quad }
\def\st{{\mathfrak t}}
\def\sT{{\mathfrak T}}
\def\sI{{\mathfrak I}}
\def\sR{{\mathfrak R}}
\def\sM{{\mathfrak M}}
\def\sA{{\mathfrak A}}
\def\sB{{\mathfrak B}}
\def\sY{{\mathfrak Y}}
\def\sE{{\mathfrak E}}
\def\sP{{\mathfrak P}}
\def\sG{{\mathfrak G}}
\def\sa{{\mathfrak a}}
\def\sw{{\mathfrak w}}
\def\se{{\mathfrak e}}
\def\sb{{\mathfrak b}}
\def\sc{{\mathfrak c}}
\def\sg{{\mathfrak g}}
\def\sd{{\mathfrak d}}
\def\sq{{\mathfrak q}}
\def\sp{{\mathfrak p}}
\def\ss{{\mathfrak s}}

\def\bm{\, \boxminus }
\def\bp{\, \boxplus }
\def\bd{\, \boxdot }

\def\Diff{\text{Diff}}
\def\Max{\text{max}}
\def\cD{\mathcal D}
\def\C{\mathbb C}
\def\R{\mathbb R}
\def\K{\mathbb K}
\def\D{\mathbb D}
\def\N{\mathbb N}
\def\M{\mathbb M}
\def\Z{\mathbb Z}
\def\a{{\underline a}}
\def\b{{\underline b}}
\def\c{{\underline c}}
\def\Log{\text{log}}
\def\loc{\text{loc}}
\def\inta{\text{int }}
\def\det{\text{det}}
\def\exp{\text{exp}}
\def\Re{\text{Re}}
\def\lip{\text{Lip}}
\def\leb{\text{Leb}}
\def\dom{\text{Dom}}
\def\diam{\text{diam}\:}
\def\supp{\text{supp}\:}
\newcommand{\ovfork}{{\overline{\pitchfork}}}
\newcommand{\ovforki}{{\overline{\pitchfork}_{I}}}
\newcommand{\Tfork}{{\cap\!\!\!\!^\mathrm{T}}}
\newcommand{\whforki}{{\widehat{\pitchfork}_{I}}}
\newcommand{\marginal}[1]{\marginpar{{\scriptsize {#1}}}}

\theoremstyle{plain}
\newenvironment{sublemm}{\begin{enonce}{Sublemma}}{\end{Sublemma}}
\newenvironment{Importantremark}{\begin{enonce}{Important remark}}{\end{enonce}}
\newenvironment{prob}{\begin{enonce}{Problem}}{\end{enonce}}
\newenvironment{fact}{\begin{enonce}{Fact}}{\end{enonce}}
\newenvironment{ques}{\begin{enonce}{Question}}{\end{enonce}}
\theoremstyle{remark}
\newenvironment{Claim}{\begin{enonce}{Claim}}{\end{enonce}}
\newenvironment{PartAns}{\begin{enonce}{Partial Answer}}{\end{enonce}}
\newenvironment{remas}{\begin{enonce}{Remarks}}{\end{enonce}}
\newenvironment{Examples}{\begin{enonce}{Examples}}{\end{enonce}}
\newenvironment{exam}{\begin{enonce}{Example}}{\end{enonce}}
\newtheorem{theorem}{Theorem}
\renewcommand*{\thetheorem}{\Alph{theorem}}
%\newtheorem{theo}{\bf Theorem}[section]
%\newtheorem{lemm}[theo]{\bf Lemma}
%\newtheorem{ques}[theo]{\bf Question}
%\newtheorem{sublemm}[theo]{\bf Sublemma}
%\newtheorem{IH}[theo]{\bf Extra induction hypothesis}
%\newtheorem{prop}[theo]{\bf Proposition}%[section]
%\newtheorem{coro}[theo]{\bf Corollary}%[section]
%\newtheorem{Property}[theo]{\bf Property}%[section]
%\newtheorem{Claim}[theo]{\bf Claim}
%\theoremstyle{remark}
%\newtheorem{rema}[theo]{\bf Remark}
%\newtheorem{important rema}[theo]{\bf Important remark}
%\newtheorem{remas}[theo]{\bf Remarks}
%\newtheorem{fact}[theo]{\bf Fact}
%\newtheorem{partial answer}[theo]{\bf Partial answer}
%\newtheorem{exem}[theo]{\bf Example}
%\newtheorem{Examples}[theo]{\bf Examples}
%\newtheorem{defi}[theo]{\bf Definition}

%\newcommand\relatif{{\rm \rlap Z\kern 3pt Z}}
%\makeatletter
%\renewcommand\theequation{\thesection.\arabic{equation}}
%\@addtoreset{equation}{section}
%\makeatother

\title{Zoology in the Hénon family: twin babies and Milnor's swallows}

\author{Pierre Berger
%\footnote{Université Paris 13, Sorbonne Paris Cité, LAGA, CNRS, UMR 7539, F-93430, Villetaneuse, France.} 
}
\date{}

\maketitle
\begin{flushright}
\itshape En mémoire de Jean-Christophe.\\
\end{flushright}

\begin{abstract} 
We study $C^{d,r}$-Hénon-like families $(f_{a\, b})_{a\, b}$ with two parameters $(a,b)\in \R^2$. We show the existence of an open set of parameters $(a,b)\in \mathcal D$, so that a renormalization chart conjugates an iterate of $f_{a\, b}$ to a perturbation of $(x,y)\mapsto ((x^2+c_1)^2+c_2,0)$. We prove that the map $(a,b)\in \mathcal D\mapsto (c_1,c_2)$ is a $C^d$-diffeomorphism; as first numerically conjectured by Milnor in 1992.

Furthermore, we show the existence of an open set of parameters $(a,b)$ so that $f_{a\, b}$ displays exactly two different renormalized Hénon-like maps whose basins union attracts Lebesgue a.e. point with bounded forward orbit. A great freedom in the choice of the renormalized parameters enables us to deduce in particular the existence of a (unperturbed) Hénon map with exactly $2$ attracting cycles (an answer to a Question by Lyubich).

The proof is based on a generalization of puzzle pieces for Hénon-like maps, and on a generalization of both the affine-like formalism of Palis-Yoccoz and the cross map of Shilnikov. The distortion bounds enable 	us to define (for the first time) $C^{r}$ and $C^{d,r}$-renormalizations and multi-renormalizations with bounds on all the derivatives.
\end{abstract}
\begin{altabstract} 
Nous considérons des $C^{d,r}$-familles $(f_{a\, b})_{a\, b}$ d'applications d'allure Hénon $f_{a\, b}$ à deux paramètres $(a,b)\in \R^2$. Nous montrons l'existence d'un ouvert de paramètres $(a,b)\in \mathcal D$ tel qu'une carte de renormalisation conjugue un itéré de $f_{a\, b}$ avec une perturbation de $(x,y)\mapsto ((x^2+c_1)^2+c_2, 0)$. Nous montrons que l'application $(a,b)\mapsto (c_1,c_2)$ est un $C^d$-difféomorphisme ; cela prouve une conjecture numérique de Milnor. 

Aussi, nous montrons l'existence d'un ouvert de paramètres $(a,b)$ tel que $f_{a\, b}$ possède exactement deux applications de Hénon renormalisées dont l'union des deux bassins attire presque tout point ayant une orbite bornée. Une grande liberté dans le choix des paramètres des deux renormalisations nous permet en particulier de montrer l'existence d'une application (non-perturbée) ayant exactement deux cycles attractifs (une réponse à une question de Lyubich). 

La preuve est basée en particulier sur une généralisation des pièces de puzzle pour les applications d'allure Hénon ainsi qu'une généralisation des représentations d'allure affine de Palis-Yoccoz et des applications croisées de Shilnikov. Les bornes de distorsion obtenues nous permettent de définir (pour la première fois) des $C^r$ et $C^{d,r}$-renormalisations et multi-renormalisations avec des bornes sur toutes les dérivées.
\end{altabstract}
 
\tableofcontents
Jean-Christophe Yoccoz used to create and develop powerful analytical and topologically-combinatorial tools to explore and then describe concrete dynamical systems in a very sharp way. In this work in his memory, we will develop two paradigms of his approach.

The topological and combinatorial notion of puzzle pieces \cite{BY,Y97} and the affine-like iterate \cite{PY01,PY09} will be here generalized to a concept of piece (involving strips which are not necessarily wide). 
 In the Hénon-like context, we will obtain a similar objects to those in \cite{berhen}.

The distortion bounds of \cite{PY01,PY09} on affine-like representations are extended to the $C^{d,r}$-topology. We show that they are preserved for arbitrarily long ($\star$)-product of hyperbolic pieces. The affine-like representation was priory (and independently) studied by the Shilnikov school in the $C^r$-topology \cite{Sh67, GST93, GST08}. By merging these two viewpoints, we will obtain both more precise and more general bounds. They will enable us to define Hénon-like and a (new) multi-renormalization (see Theorem \ref{Charthenonlike} \textsection \ref{sectionthmC} and Theorem \ref{Charthenonlike2} \textsection \ref{sectionthmD}) without loss of regularity. As a matter of fact, for the first time in dimension 2, infinitely renormalizable maps can be defined in the $C^{r}$-topology. 

Such tools, applied to the Hénon-like families, will enable us to solve a few problems (Theorems \ref{thm B} \textsection\ref{texsectionthmB} and \ref{thm A} \textsection\ref{texsectionthmA}). 
\section{Statement of Theorems \ref{thm B} and \ref{thm A} on Hénon-like families}
\subsection{Topologies}
For $r<\infty$, the space of $C^r$-maps $g$ from $\R^n $ to $\R^m$ is endowed with the norm $\|g\|_{C^r} = \max_{0\le j\le r} \sup_x \|D^j g(x)\|$. For $d\le r$, a family $(g_a)_a$ is of class $C^{d, r}$ if all the derivatives 
$\partial_a^iD^j g_a(x)$ are well defined and continuous for $i+j\le r$ and $i\le d$, with $D$ is differential w.r.t. the variable $x$. The $C^{d,r}$-norm of $(g_a)_a$ is $\max\{\sup_{a,x} \| \partial_a^iD^j g_a(x)\|:{i+j\le r,\, i\le d}\}$. 

If $d=r$, we will say that the family is of class $C^r$. Note that a family is of class $C^r$, if and only if the map $(a,x)\mapsto g_a(x)$ is of class $C^r$. 

For $r=\infty$, we endow the space of $C^\infty$-maps with the distance $d(g,h)= \sum_{s\ge 0} \min (\|g-h\|_{C^{s}},2^{-s})$ and similarly for the space of $C^{\infty}$-families of maps.

\subsection{Hénon-like family}
The main results are proved for 2-parameter families of Hénon-like diffeomorphisms. These families are generalizations of the Hénon family $(h_{a\, b})_{a\, b}$:
\[h_{a\, b}: \; (x,y)\in \R^2\mapsto (x^2+a-b y, x)\in \R^2\; .\]

\begin{defi}[Hénon-like map]A map $f $ of a domain $D\subset \R^2$ into $\R^2$ is \emph{$C^r$-Hénon-like} with parameter $(a,b)\in \R^2$ if it is a perturbation of $h_{a\; b}$ of the form.
\[(x,y)\mapsto h_{a\, b}(x,y) + (\zeta , \xi) (x,b y)\; ,\]
where $\zeta$ and $\xi $ are $C^r$-functions from $D$ into $\R$. If the $C^r$ norms of $\zeta$ and $\xi $ are smaller than $\delta$, then $f$ is \emph{$\delta$-$C^r$-Hénon-like}. 
\end{defi}

\begin{rema} The determinant of the Hénon-like map $f $ is in $[(1-3\delta)b,(1+3\delta)b]$. Hence the determinant is small whenever $b$ is small. 
On the other hand, given a Hénon-like map, the parameter $(a,b)$ is not uniquely defined.
\end{rema}

\begin{defi}[Hénon-like family]
 For any $d\le r\le \infty$, a \emph{$\delta$-$C^{d,r}$-Hénon-like family} $(f_{a\, b})_{I\times J}$ (of multiplicity $m=1$) consists of mappings of the form:
 $$f_{a\, b}(x,y)=(x^2+a-b y+\zeta_{a\, b}(x,b y), x+\xi_{a\, b}(x,b y)) \, ,$$
where $(\zeta_{a\, b})_{I\times J}$ and $(\xi_{a\, b})_{I\times J}$ are $C^{d,r}$-families of mappings of norm at most $\delta$. The domain $I\times J$ of the parameters $(a,b)$ is the product of two segments $I$ and $J$. 
 \end{defi}
 For instance the original Hénon family $(h_{a\, b})_{a\, b}$ is $\delta$-$C^{d,r}$-Hénon-like, for every $\delta>0$ and $r\ge 1$.

Let us define the Hénon-like families of higher multiplicities which appear at every 2 parameters, non-degenerated unfolding of a homoclinic tangency (see Example \ref{PTrevisited}). These families also appear as the renormalizations of Hénon-like families of lower multiplicities (see in Theorem \ref{Charthenonlike} \textsection \ref{sectionthmC}, Theorem \ref{Charthenonlike2} \textsection \ref{sectionthmD},
 and Examples \ref{examren1} and \ref{examren2}).

 \begin{defi}[Wide Hénon-like family of multiplicity $m\ge 1$] For any $d\le r\le \infty$, a \emph{$\delta$-$C^{d,r}$-Hénon-like family $(f_{a\, b\, m})_{I\times J}$ of multiplicity $m$ } consists of mappings of the form:
 $$f_{a\, b\, m}(x,y)=(x^2+a-b^m y+\zeta_{a\, b}(x,b^m y), x+\xi_{a\, b}(x,b^m y)) \, ,$$
where $(\zeta_{a\, b})_{I\times J}$ and $(\xi_{a\, b})_{I\times J}$ are $C^{d,r}$-families of mappings of norms at most $\delta$. 

The family is \emph{stretched} if the interval $I$ contains $[-2,1/1/4]$. 
The family is \emph{$L$-wide} if it is stretched and $\max_J |b|^m \ge L \min_J |b|^m$. 
\end{defi}
 Sometime, we will denote the parameter $(a,b)$ by $p$. 
\begin{center}
\emph{Thanks to the following proposition we will always assume $\xi_{a\, b}=0$.}\end{center}
\begin{prop}\label{zetanull}
For every $C^{d,r}$-Hénon-like family $(f_{a\, b\, m})_{a\, b}$ of multiplicity $m\ge 1$, up to a coordinate change, we may assume that $\xi_{a\, b}=0$: $f_{a\, b\, m}(x,y)=(x^2+a-b^m y+\zeta_{a\, b}(x,b^m y),x)$. 
\end{prop}
\begin{proof}
The implicit function theorem applied to $(x,y,t)\mapsto \xi_{a\, b}(x-t,y)-t$ gives the existence of a $C^{d,r}$-small family of functions $(\rho_{a\, b})_{a\, b}$ such that $\xi_{a\, b}(x-\rho_{a\, b}(x,y),y)=\rho_{a\, b}(x,y)$. Then the diffeomorphism $(x,y)\mapsto (x-\rho_{a\, b}(x,b y),y)$ conjugates $f_{a\, b\, m}$ to a Hénon-like map with the requested form. 
\end{proof}
\subsection{Swallow-like family}
In \cite{Mil92}, Milnor studied the composition of two quadratic maps $Q_{a\, b}= Q_b\circ Q_a$ with $Q_a(x)= x^2+a$ and $Q_b(x)= x^2+b$. He observed that the following parameter space looks like a swallow (see fig. \ref{swallow_bassin}):
 $$\mathcal S:=\{(a,b)\in \R^2: Q_{a\, b}^n (0)\not\to \infty\text{ or }Q_{b\, a}^n (0)\not\to \infty\}$$ 
 We can split $\mathcal S$ into two regions. The \emph{wings} of the swallow is the set of parameters $(a,b)$ so that exactly one of the orbits $(Q_{a\, b}^n(0))_{n\ge 0}$ and $(Q_{b\, a}^n(0))_{n\ge 0}$ is unbounded. The \emph{body} of the swallow is the set of parameters $(a,b)$ so that both orbits $(Q_{a\, b}^n(0))_{n\ge 0}$ and $(Q_{b\, a}^n(0))_{n\ge 0}$ are bounded.
\begin{prop}
The body of the swallow is bounded by the three analytic curves:
\[\mathcal C_1= \{(a,b): Q_{a\, b}(b)=-b \text{ and }DQ_{a\, b}(b)<-1\},\; \mathcal C_2= \{(a,b): Q_{b\, a}(a)=-a \text{ and }DQ_{b\, a}(a)<-1\},\]
\[\mathcal C_3= \{(a,b): \exists X Q_{a\, b}(X)=X\; ,\; DQ_{a\, b}(X)=1\text{ and } 
 DQ_{a}(X)>0 \}\; .\]
\end{prop}
\begin{figure}[h!]
 \centering
 \includegraphics[width=10cm]{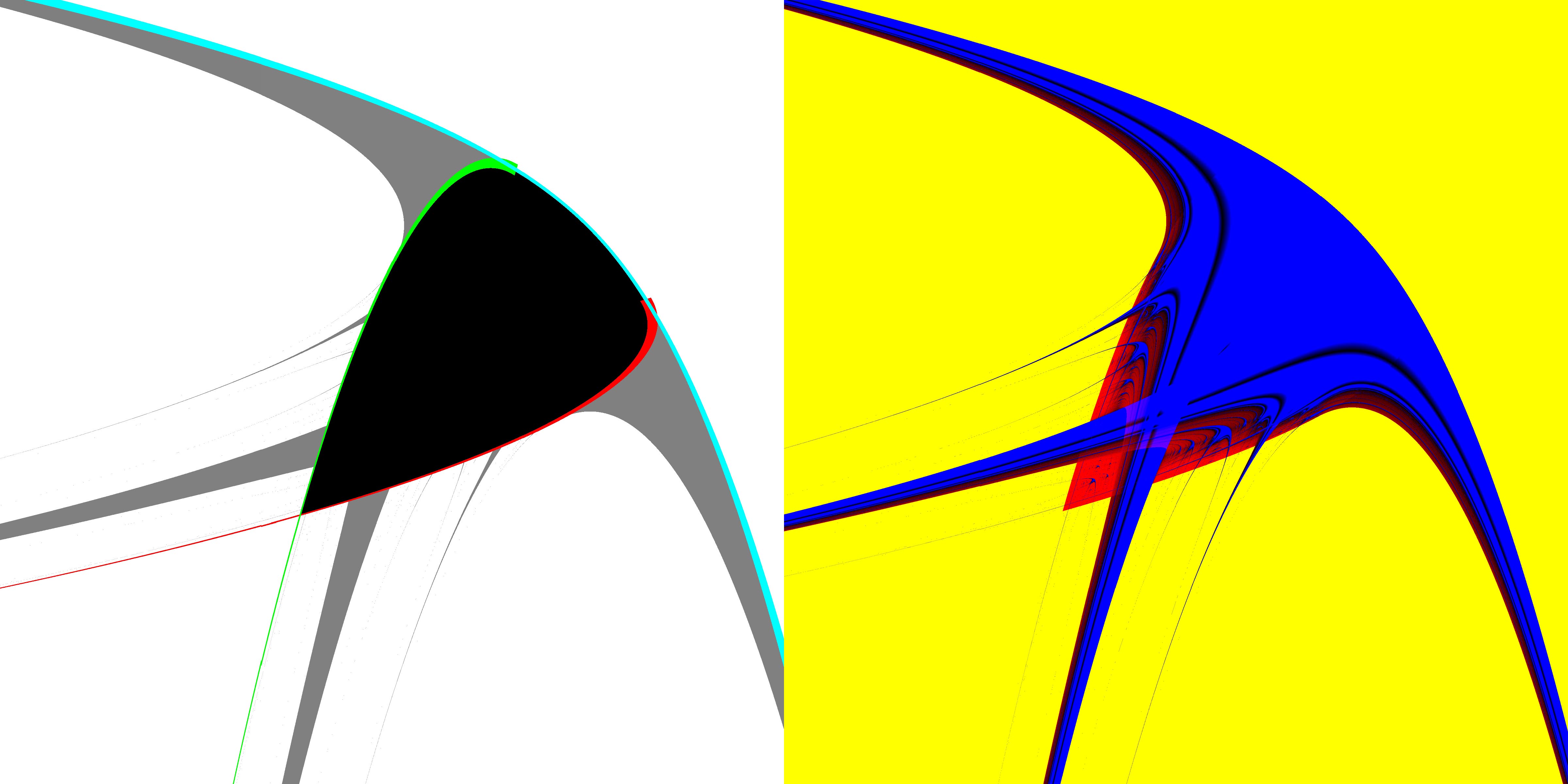}
 \caption{
 Two pictures of the Swallow parameter space $\mathcal S$, with $a$ the $x$-coordinate and $b$ the $y$-coordinate. 
 At the Left is drawn the body of the swallow in black, the wings in gray, the curve $\mathcal C_1$ in red, the curve $\mathcal C_2$ in green and the curve $\mathcal C_3$ in cyan.
At the right, the swallow is depicted in function of the values of $\frac1N \log |DQ_{b\, a}^{N}|(b)$ and $\frac1N \log |DQ_{a\, b}^{N}|(b)$, with $N=10\, 000$. It is black if both are close to 0, red if one is negative, blue if one is positive. 
 }\label{swallow_bassin}
 \end{figure} 
\begin{proof}For every $(a,b)$ in the swallow, 
the immediate basin of the infinity for $Q_{a\; b}$ is the complement of a segment $[-\gamma,\gamma]$. The boundary $\{-\gamma,\gamma\}$ of the immediate basin complement is sent into itself by $Q_{a\, b}$, and $Q_{a\, b}$ sends $[b,+\infty]$ into itself since it is $Q_b(\R)$. Thus $\gamma$ is a fixed point of $Q_{a\, b}$. Furthermore 
$D_\gamma Q_{a,\, b}\ge 1$. Let $\bar \gamma:= Q_a(\gamma)$. Note that $DQ_a(\gamma)>0$. 

Let $(a,b)$ be in the body of the swallow.
The image $Q_a([-\gamma, \gamma])$ is $[a,\bar \gamma]$. 
If $a<-\bar \gamma$, then its image by $Q_b(a)=Q_{a\; b}(0)$ would be greater than $\gamma$ and so $0$ would be in the basin of the infinity. Thus $Q_a([-\gamma,\gamma])$ is included in $ [-\bar \gamma,\bar \gamma]$. Similarly, 
$Q_b([-\bar \gamma,\bar \gamma])$ is included in $ [- \gamma,\gamma]$. These two inclusions properties are also sufficient conditions for $(a,b)$ to be in the body. They are robust properties except when:
\begin{itemize}
\item $b = -\gamma$ and so $Q_{a\; b}(b)=-b$. Also $D Q_{a\; b}(b)=-D Q_{a\; b}({\gamma})< -1$. This condition defines $\mathcal C_1$. 
\item $a= -\bar \gamma$ and so $Q_{b\; a}(a)=-a$. Also $D Q_{b\; a}(a)=-D Q_{b\; a}({\bar \gamma})< -1$. This condition defines $\mathcal C_2$. 
\item $\gamma$ is a parabolic fixed point, and so $(a,b)\in \mathcal C_3$. 
\end{itemize}
\end{proof}

\begin{defi}[Swallow-like family]
A 2-parameter family $(S_{a\, b})_{(a,b)}$ of planar mappings is \emph{swallow-like} if it is of the form: 
\[S_{a\, b}:\; (x,y)\mapsto ((x^2+b)^2+a,0)+\varsigma_{a\; b}(x,y)\; .\]
It is $\delta$-$C^{d,r}$-swallow-like if the family $(\varsigma_{a\; b})_{(a,b)}$ is of class $C^{d,r}$ and $\delta$-$C^{d,r}$-small. It is \emph{$R$-wide}, if 
$S_{a\, b}(x,y)$ is well defined for every $(x,y,a,b)\in [-R,R]^4$. 
\end{defi}

\begin{prop}\label{compo2Swallow}
Let $I$ be a neighborhood of $0\in \R$, let $K>0$, let $c:=(a,a')\in I^2\mapsto (b_c, b'_c)\in \R^2$ be a $C^d$-small map. Let $(\zeta_c, \xi_c)_c$ and $(\zeta'_c,\xi'_c)_c$ be $C^{d,r}$-small families. 
Put $f_c(x,y)= (x^2+a -b_c +\zeta_c(x,b_c y), x+\xi_c(x,b_c y))$ and
 $f'_c(x,y)= (x^2+a' -b'_c +\zeta'_c(x,b'_cy), x+\xi'_c(x,b'_c y))$. 

Then the family $(f_{c} \circ f'_{c})_c$ is $C^{d,r}$-conjugated to a $C^{d,r}$-swallow-like family. Moreover the latter family is $R$-wide if $I$ contains $[-R,R]$. 
\end{prop}
 \begin{proof}
 Let $\epsilon>0$ be small but large compared to $\max (\|b_c\|_{C^d}, \|b'_c\|_{C^d})$. Let $\psi :(x,y)\mapsto (x, \epsilon\cdot y)$. Then we conclude by observing that $(\psi\circ f_{c}\circ \psi^{-1})_{I^2}$ and $(\psi\circ f'_{c}\circ \psi^{-1})_{I^2}$ are $C^{d,r}$-close to respectively $((x,y)\mapsto (x^2+a,0))_{I^2}$ and $((x,y)\mapsto (x^2+a',0))_{I^2}$. \end{proof}
 
\subsection{Embedding swallows into the Hénon parameter space} \label{texsectionthmB}
Based on numerical simulations and bifurcation diagrams of El Hamouly-Mira \cite{EHM82}, Milnor \cite{Mil92} observed that the Swallow parameter space seems embedded into the parameter space of the Hénon family. In figure \ref{detailparahenonetshallow}, we can see that indeed many swallows seem embedded into the Hénon parameter space. 
Through an intuitive explanation, Milnor conjectured that given a parameter in the embedded swallow, an iterate of the dynamics of the Hénon map would have locally and approximately the same dynamics as the corresponding composed quadratic mapping. This is now established by the following:
\begin{figure}[h!]
 \centering
 \includegraphics[width=16cm]{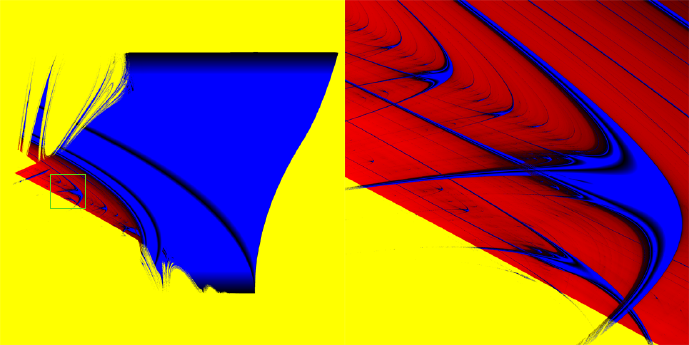}
 \caption{Parameter space of the Hénon family $(h_{a,b})_{a\; b}$ at the left, and a zoom on the green box subset at the right. 
 Again, we put $a$ at the $x$-coordinate, $b$ at the $y$-coordinate, and the color is yellow when the points $0$ escapes to infinity. Otherwise, the color varies in function of the value of $\frac1N \log \|D_0f^N(0,1)\|$ for $N=10\, 000$.}
\label{detailparahenonetshallow}
\end{figure}
\begin{theorem}[First main result]\label{thm B}
There exist $L>0$, $\delta>0$, $\hat b>0$ such that for every $0\le d\le r-3\le \infty$ and $\delta^*>0$, for every $L$-wide, $\delta$-$C^{d,r}$-Hénon-like family $(f_{a\, b})_{I\times J}$ with $J\subset [-\hat b,\hat b]$, there exists $n\ge 1$ and:
\begin{itemize}
\item a $C^d$-diffeomorphism $P$ from $(-R,R)^2$ onto an open subset of $I\times J$,
\item a $C^{d,r}$-family $(\phi_c)_{c\in [-R,R]^2}$ of diffeomorphisms $\phi_c$ from $[-R,R]^2$ onto a domain $D_{c}\subset \R^2$,
 \end{itemize}
 satisfying that $(S_c:= \phi_{c}^{-1}\circ f_{P(c)}^n\circ \phi_{c})_{c \in [-R,R]^2} $ is $\delta^*$-$C^{d,r}$-swallow-like and $R$-wide. 
 \end{theorem}
Actually the family $(S_c)_c$ is given by a composition of two Hénon-like families satisfying the hypotheses of Proposition \ref{compo2Swallow}. 
This theorem provides a new pattern (which is the Milnor Swallow) to describe the parameter space of a Hénon-like family with 2-parameters. A previous way was given by renormalization strips \cite{CLM, Ha11} equal to the domains of a Hénon-like renormalizations that we shall develop in the next section.

\subsection{Renormalization}
A unimodal map is \emph{renormalizable} if there exists an open interval $I$ sent into itself by an iterate $P^N$ of the dynamics such that $P^N|I$ is (conjugated to) a unimodal map. 
This definition can be generalized to Hénon-like maps.

\begin{defi} A \emph{pre-renormalization domain $D$ of period} $N\ge 2$ of a Hénon-like map $f$ is an open set $D$ such that there exist
a $C^r$-Hénon-like map $\mathcal R f$ and an embedding $\phi\in C^r(D\cup f^{N}(D),\mathbb R^2)$ for which the following diagram commutes:
\begin{displaymath}
 \xymatrix{
 D\ar[d]_{\phi}\ar[r]^{ f^N} & f^{N}(D)\ar[d]_{\phi}\\
 \R^2 \ar[r]_{\mathcal Rf} & \R^2 \\
 }
\end{displaymath}
\end{defi}

There are many domains which are not canonical, for instance any $D$ such that $f^N(D)\cap D=\emptyset$ is a (not so interesting) pre-renormalization domain. Let us define those which are canonical. 

We observe that for any $\delta$-$C^r$-Hénon-like map $f_{a\, b}$, if $b$ and $\delta$ are small, 
 with $\psi(x,y)=(x,-by)$, the map 
$F_{a\, b}:= \psi \circ f_{a\, b} \circ \psi^{-1}$ is a $C^r$-perturbation of:
\[\hat Q_a : (x,y)\mapsto (x^2+a+y,0).\]
Suppose $a<1/4$ so that $Q_a(x):=x^2+a$ has two different hyperbolic fixed points $\alpha_0<\beta_0$. Note that $(\beta_0,\beta_0)$ is close to a fixed point $\beta$ of $f_{a\, b}$. 
 Also $(\beta_0,0)$ is a hyperbolic fixed point $\hat Q_{a}$, with the following local stable manifold: 
 \[W^s_{loc}(\beta_0,\hat Q_a):=\left\{(x,y)\in \mathbb R\times [-1,\infty): x^2+a+y=\beta_0\right\}.\]
The curve $W^s_{loc}(\beta_0,\hat Q_a)$ intersects $\mathbb R\times\{0\}$ at $\beta_0$ and $-\beta_0$. If $\delta$ is small compared to $1/4-a$, then the fixed point persists $(\beta_0,0)$ to the fixed point $\psi(\beta)\approx (\beta_0(a),0)$ of $F_{a\, b}$, and the curve $W^s_{loc}(\beta_0,\hat Q_a)$ persists as a local stable manifold $W^s_{loc}(\psi(\beta),F_{a\, b})$ of $\psi(\beta)$. We put $W^s_{loc}(\beta,f_{a\, b}):=\Psi^{-1}(W^s_{loc}(\psi(\beta),F_{a\, b}))$. 

\begin{defi}[Canonical domain $Y_D$]\label{YD} 
For every $C^1$-$\delta$-Hénon-like map $f_{a\, b}$ with $\delta$ and $b$ small compared to $1/4-a>0$, the curve $W^s_{loc}(\beta,f_{a\, b})$ and 
the two lines $\{y=\pm 1/(8b)\}$ bound a subset $Y_D(f_{a\, b})\subset \R^2$ which is diffeomorphic to the filed square $[-1,1]^2$. 
\end{defi}
 
If moreover $a+\beta> 0$ and moreover $\delta $ and $b$ are small compared to 
$\beta+a$, the box $Y_D$ is sent into itself by $f_a$. This condition is equivalent to $a\in (-2,1/4)$ and $(b,\delta)$ small compared to $\max(a+2,1/4-a)$.

\begin{defi}\label{canrendo} A \emph{canonical pre-renormalization domain} $D_{can}$ of a Hénon-like map $f$ is a pre-renormalization domain such that $\phi(D_{can})=Y_D(\mathcal Rf)$.

 A \emph{renormalization domain} $D$ is a canonical pre-renormalization domain such that 
\[f^N(D)\subset D\Leftrightarrow \mathcal Rf(Y_D)\subset Y_D,\quad \text{with } D=Y_D(\mathcal Rf)\; .\]
\end{defi}

In the parameter space of the Hénon family (see figure \ref{detailparahenonetshallow}), we see many blue strips (which go from the upper left side, toward the lower right side and cross the line $\{b=0\}$). These strips correspond to parameters for which the Hénon map is renormalizable and their renormalization displays a sink (this will be a consequence of Theorem \ref{Charthenonlike}, see also\cite{CLM, Ha11}). When $b=0$, the Hénon map $h_{a\, 0}$ is semi-conjugated to a quadratic map $Q_a$. For $a$ in this strip, $Q_a$ displays a sink. Thus these strips correspond to the hyperbolic continuation of a sink for a quadratic map. A quadratic map displays at most one sink, but these strips may overlap at $b\not = 0$. Actually, by adapting Newhouse proof \cite{Newhouse}, it is easy to show \cite{Ro83, BedeSi} the existence of a locally Baire generic set $N$ of parameters $(a,b)$ (with $b$ arbitrarily small) for which the Hénon map displays infinitely many sinks. These proofs use the fact that infinitely many aforementioned strips overlap at $(a,b)$. Van Strien \cite[Question 1.10]{VS} asked if this bifurcation analysis could explain any co-existence of infinitely many sinks in the Hénon family. In other words, if the co-existence of infinitely many sinks only occurs when there is an infinite overlapping of strips which all intersect $\{b=0\}$).
 Theorem \ref{thm B} gives a (sketchy) negative answer to his question.
The wings of the swallow are given by infinitely many strips shaped like parabolas. The above wing of the main swallow in Figure \ref{detailparahenonetshallow} are formed by a union of parabolic strips which all intersect the line $\{b=0\}$. On the other hand the bottom wings are formed by a union of parabolic strips which \emph{are all disjoint} to the line $\{b=0\}$. Then by implementing Newhouse argument to this setting, it is not hard to display of parameter which is in infinitely many parabolic strips, all of them being disjoint from $\{b=0\}$.
\subsection{Twin babies}\label{texsectionthmA}
When $Q_a$ is a quadratic map with an attracting cycle, it is easy to show that for $b$ small enough, $h_{a\, b}$ displays a unique attractor. Other examples of Hénon maps with a unique attractor exists, such as the Benedicks-Carleson's attractor \cite{BC2,berhen} and the odometer of an infinitely renormalizable Hénon-like map \cite{CLM, Ha11}. We saw above that some Hénon maps display infinitely many attractors. A natural question is:

\begin{ques}[Lyubich]\label{lyubich question} Does there exist a parameter in the Hénon family such that there exists a finite number of attractors which is not 1?\end{ques}
We give an answer to Lyubich's question. The following shows that there exist Hénon maps with exactly $2$ sinks. Moreover, a Hénon-like map can display two different renormalized Hénon maps (so called \emph{twin babies}) which attract Lebesgue a.e. point which does not escape to infinity:
\begin{theorem}[Second main result]\label{thm A}
For every $0\le d\le r-3\le \infty$, 
there exist $L>0$, $\delta>0$ and $\hat b>0$, so that for every $\delta^*>0$, for every $L$-wide, $\delta$-$C^{d,r}$-Hénon-like family $(f_{p})_{p\in I\times J}$, with $J\subset [-\hat b, \hat b]$, has the following property:

There exist a sub-domain $\cD\subset I\times J$ and diffeomorphisms
$P^+\colon \cD\to \R^2$ and $P^-\colon \cD\to \R^2$ such that for every $p\in \cD$, with $p^+=P^+(p)$ and $p^-=P^-(p)$, the map $f_{p}$ has two canonical prenormalization domains $D^+(p)$ and $D^{-}(p)$ associated to Hénon-like maps $\mathcal R^+ f_{p^+}$ and $\mathcal R^{-} f_{p^{-}}$ such that:
\begin{enumerate}[(1)] 
\item The families $(\mathcal R^+ f_{p^+})_{p^+}$ and $(\mathcal R^{-} f_{p^-})_{p^-}$ are $C^{d,r}$-$\delta^*$-Hénon-like with higher multiplicities.
\item The family $(\mathcal R^+ f_{p^+})_{p^+\in P^+(\cD)}$ is stretched.
\item The set $ P^{-}(\cD)$ is included in $[-0.1,0.1]\times [-\hat b, \hat b]$. In particular, for every $p\in \cD$,
the set $D^{-}(p)$ is always a renormalization domain which is included in the closure of the basin of a sink.
\item For every $p^+ \in P^{+}(\cD)$, if $a_{p^+}$
 is in $[-2+\delta^* , 1/4-\delta^*]$, then $D^+(p)$ is a renormalization domain and Lebesgue almost every point in $\R^2$ has its $f_p$-forward orbit which escapes to infinity or eventually lands into $D^+(p)$ or $D^{-}(p)$.
\end{enumerate}
\end{theorem}
\begin{rema} We set up $a_{p^-}$ nearby $0$ to make the dynamics of $\mathcal R^{-} f_{p^{-}}$ elementary: it has one attracting fixed point which attracts almost every point which does not escape to infinity. However, we could have chosen $a_{p^-}$ nearby any other number in $(-2,1/4)$ and even at any one-codimensional phenomena such as an infinitely renormalizable parameter or a parabolic attractor.
\end{rema}

%Here is a strong answer to Lyubich's question:
%\begin{coro}
%For every $N\ge 1$, there exists a Hénon map $h_{a\, b}$ which displays exactly $N$ sinks, and the union of the basins of theses sinks contains Lebesgue a.e. point which does not escape to infinity. Moreover $b$ can be assumed arbitrarily small. 
%\end{coro}
%\begin{proof} For $N=1$ the answer is well known. Let $N\ge 2$. Assume by induction that there exists a domain $\mathcal D'\subset [-b,b]\times \R$, a map $P^+$, a canonical prenormalization domains $D^+(p)$ associated to a Hénon-like map $\mathcal R^+ f_{p^+}$ 
%so that:
%\begin{itemize}
%\item for every $p\in \mathcal D'$ there are at least $(N-1)$ sinks $(S_{i\, p})_{1\le i<N}$,
%\item the family $(\mathcal R^+ f_{p^+})_{p^+\in P^+(\cD)}$ is wide,
%\item for every $p\in \cD$ such that the first coordinate of $P^+(p)$ is in $[-2+\eta, 1/4-\eta]$ and so $D^+(p)$ is a renormalization domain, then Lebesgue almost every point in $\R^2$ with bounded orbit 
%%which does not escape to infinity by $f_{p}$
% eventually lands into $D^+(p)$ or in the basin of some $S_{i\, p}$. 
%\end{itemize}
%The step $N=2$ is an immediate consequence of Theorem \ref{thm A}. The step $N\to N+1$ is obtained by applying Theorem \ref{thm A} to $(\mathcal R^+ f_{p^+})_{p^+\in P^+(\cD)}$.
%\end{proof}

\section{Concepts of pieces and affine-like representations}
\subsection{Yoccoz' puzzle pieces}\label{3.1}

Let $a<0$. We recall that $Q_a(x)=x^2+a$. The polynomial $Q_a$ displays two fixed points $\alpha= \alpha(a)<0<\beta=\beta(a)$.
We observe that $\alpha_0:= -\alpha$ belongs to $(0,\beta)$. 
Put $\R_\se:= [-\alpha_0, \alpha_0]$. 
\begin{defi}
A \emph{piece} is a pair $\sb=(\R_\sb,n_\sb)$ of a (non-empty) segment $\R_\sb\subset \R$ and an integer $n_\sb\ge 0$ so that 
the restriction 
$Q_a^{n_\sb}|\R_\sb$ is a diffeomorphism into $\R_\se$. A 
\emph{puzzle piece} $\sb$ is a piece such that 
$Q_a^{n_\sb}(\R_\sb)=\R_\se$. 
 \end{defi}
For instance $(\R_\se,0)$ is a puzzle piece. 
Let $\alpha_1$ be the preimage of $\alpha_0$ by $Q_a|\R^+$, and let $\alpha_2$ be the preimage of $\alpha_1$ by $Q_a|\R^+$. 
We notice that $\sw_+:= ([\alpha_0,\alpha_1],1)$ , $\sw_-:= ([-\alpha_1,-\alpha_0],1)$ and $\sw_{=}:= ([-\alpha_2,-\alpha_1],2)$ are puzzle pieces.

Puzzle pieces enable to describe the combinatory of the quadratic map $Q_a$ and select the parameter $a$ in function of this combinatory (as done in \cite{Y97}).

By \cite[Prop. 3.1. (2)]{Y97} there exist unique parameters $a_{2}<a_1$ so that $Q_{a_2}(0)=a_2= -\alpha_2(a_2)$ and $Q_{a_1}(0)=a_1= -\alpha_1(a_1)$. Moreover, for $a\in (a_2,a_1)$, the critical value $a=Q_{a}(0)$ is in the interior of $\R_{\sw_=}:=[-\alpha_2,-\alpha_1]$. 

\begin{rema}\label{remasupport} If $a\in (a_2,a_1)$, the complement of the basin of the infinity is $[-\beta,\beta]$. Also every point in $(-\beta,\beta)$ is sent by an iterate of $Q_a$ into $ [a ,Q_{a}(a)]$. The latter segment is in the interior of $\R_{\sw_=}\cup \R_{\sw_-}\cup \R_\se\cup \R_{\sw_+}$. 
\end{rema}

Let us assume that $a< a_1$. Let $-\tilde \alpha_2<0<\tilde \alpha_2$ be the two preimages of $-\alpha_1$ by $Q_a$. 
We notice that $\ss_+:= ([\tilde \alpha_2,\alpha_0],2)$ and $\ss_-:= ([-\alpha_0,-\tilde \alpha_2],2)$ are two puzzle pieces with segments included in $\R_\se$. 

When $a\in [a_2,a_1)$, we put $\boxdot= (\R_\boxdot, n_\boxdot)$, with $\R_{\boxdot}= \R_\se\setminus int(\R_{\ss_-}\cup \R_{\ss_+})$ and $n_\boxdot=3$. We notice that $\boxdot$ is not a piece since $Q_a|\R_\boxdot$ is not injective. 

\begin{rema}\label{imageRe} When $a\in [a_2,a_1)$, $Q_a$ sends the segment $\R_{\boxdot}$ into $ \R_{\sw_=}$ and the segments $\R_{\ss_-}$ and $\R_{\ss_+}$ into $\R_{w_-}$.\end{rema}

The $\star$-product is a binary operation on pieces which enables to construct new pieces. 

\begin{defi}
Let $\ss$ and $\ss'$ be pieces such that $\R_{\ss'}\subset Q_a^{n_\ss}(\R_\ss)$. Then the following is a piece:
\[\ss\star \ss':= (\R_\ss\cap Q_a^{-n_\ss}(\R_{\ss'}), n_\ss+n_{\ss'})\; .\]
If moreover $\ss'$ is a puzzle piece, then $\ss \star \ss'$ is a puzzle piece. 
\end{defi} 

\begin{Examples}
For instance $\sw_= \star \ss_+ $ is a well defined puzzle piece whereas $ \ss_+\star \sw_=$ is not well defined (it is not even a piece). 
\end{Examples}
\begin{fact} The $\star$-product is associative: if 
 $\ss,\ss',\ss''$ are pieces such that $(\ss\star \ss')\star \ss''$ and $\ss\star (\ss'\star \ss'')$ are well defined then:
 \[(\ss\star \ss')\star \ss''= \ss\star (\ss'\star \ss'')\; .\]
We will denote this product by $\ss\star \ss'\star \ss''$. 
\end{fact}

\subsection{Boxes }
A \emph{box} is a subset of $\R^2$ which is diffeomorphic to the filled square $[-1,1]^2$. 

Let $r\in [2,\infty]$. Let $f$ be a $\delta$-$C^r$-Hénon-like map (of multiplicity 1) for a parameter $(a,b)$ with $a\in (2,a_1)$ and $b$ small. We assume $b$ small and $a_1-a$ 
large compared to $\delta$ and $b$. 
\label{Puzzlepieceforhenon}

Let us define dynamically some boxes. Put $Q=Q_a$ and let $h:(x,y)\mapsto (Q(x),x)$. Let $K_0:=\{\alpha_0,-\alpha_0,\alpha_1,-\alpha_1, \tilde \alpha_2,-\tilde \alpha_2,\beta,-\beta\}$. The set $K_0$ is invariant and uniformly hyperbolic for $Q$. 
The product $\hat{\mathcal L}:=K_0\times [-3,3]$ is made by local stable manifolds of points in $\{(k,k): k\in K_0\}$, for the dynamics $h$. 

For every $k\in K_0$, we set $\hat{\mathcal L_k}:=\{k\}\times [-3,3]$.
As $f$ is $C^r$-close to $h$, this continuous family of local stable manifolds persists for $f$ as a union:
\[\mathcal L:= \cup_{k\in K_0} {\mathcal L}_k\quad \text{satisfying}\quad 
f({\mathcal L}_k)={\mathcal L}_{Q(k)}\quad \forall k\in K_0\; .\]
Furthermore, each $\mathcal L_k$ goes from the line $\{y=-3\}$ to $\{y=3\}$ and passes nearby $(k,0)$. 

Let $Y_\se$ be the set bounded by $\mathcal L_{-\alpha_0}$, $\mathcal L_{\alpha_0}$, and the two segments $\{y=\pm 3\}$. This set depends continuously on $f$. The set $Y_\se$ is a box.

Let $\partial^s Y_\se:= \mathcal L_{-\alpha_0}\cup \mathcal L_{\alpha_0}$ and let $\partial^u Y_\se:= Y_\se\cap \{|y|= 3\}$. 
 We notice that $\partial Y_\se= \partial^u Y_\se\cup \partial^s Y_\se$.

Likewise, for $\sd\in \{\ss_-, \ss+, \sw_-, \sw_+, \sw_=,\boxdot\}$ we define
$ \partial^s Y_\sd:= \cup_{k\in \partial \R_\sd} \mathcal L_{k}$ and the two segments $\partial^u Y_\sd$ of $\{(x,y): |y|= 3\}$ joining the endpoints of these curves. Let $Y_\sd$ be the box bounded by $ \partial^s Y_\sd\cup \partial^u Y_\sd$. Figure \ref{notation_zoologie_henon} depicts the boxes $Y_\sd$. 

We note that $f$ sends $Y_{\ss_\pm }$ into $Y_{\sw_-}$ and $Y_{\sw_= }$ into $Y_{\sw_+}$, 	and $Y_{\sw_\pm}$ into $Y_\se$. This implies:

\begin{fact}\label{boxinclusion}
Assume that $a_1-a$ is positive and large compare to $\delta $ and $b$. 
For every $\sd\in \{\ss_-, \ss_+, \sw_-, \sw_+, \sw_=\}$, the map $f^{n_\sd}$ sends $Y_\sd$ into $Y_\se$ and the two components of $\partial^sY_\sd$ are sent into different components of $\partial^s Y_\se$. 

If moreover $a>a_2$ and if $\delta$ and $b$ are small compared to $a-a_2$, the box $Y_\boxdot $ is sent by $f$ into $Y_{\sw_=}$ and by $f^3$ into $Y_\se$ ; both components of $\partial^s Y_\boxdot $ are sent by $f^3$ into a same component of $\partial^s Y_\se$. 
\end{fact}
\begin{coro}\label{remasupport2} If $a\in (a_1,a_2)$ and if $\delta$ and $b$ are small compared to  $\max(|a-a_2|,|a-a_1|)$, then every point not-attracted by $\infty$ is sent by an iterate of $f$ into $\mathcal L_\beta\cup Y_{\sw_-}\cup Y_{\sw_=}$. \end{coro}
\begin{proof}
By Remark \ref{remasupport}, every point not attracted by the infinity and not in $\cup_n f^{-n}(\mathcal L_\beta)$ is sent by an iterate of $f$ into $Y_\se\cup Y_{\sw_=}\cup Y_{\sw_-}\cup Y_{\sw_+}$. By fact \ref{boxinclusion}, one of its iterates is sent into $Y_\se$ and then into $ Y_{\sw_-}\cup Y_{\sw_=}$. 
\end{proof}

\begin{figure}
 \centering
 \includegraphics[width=15cm]{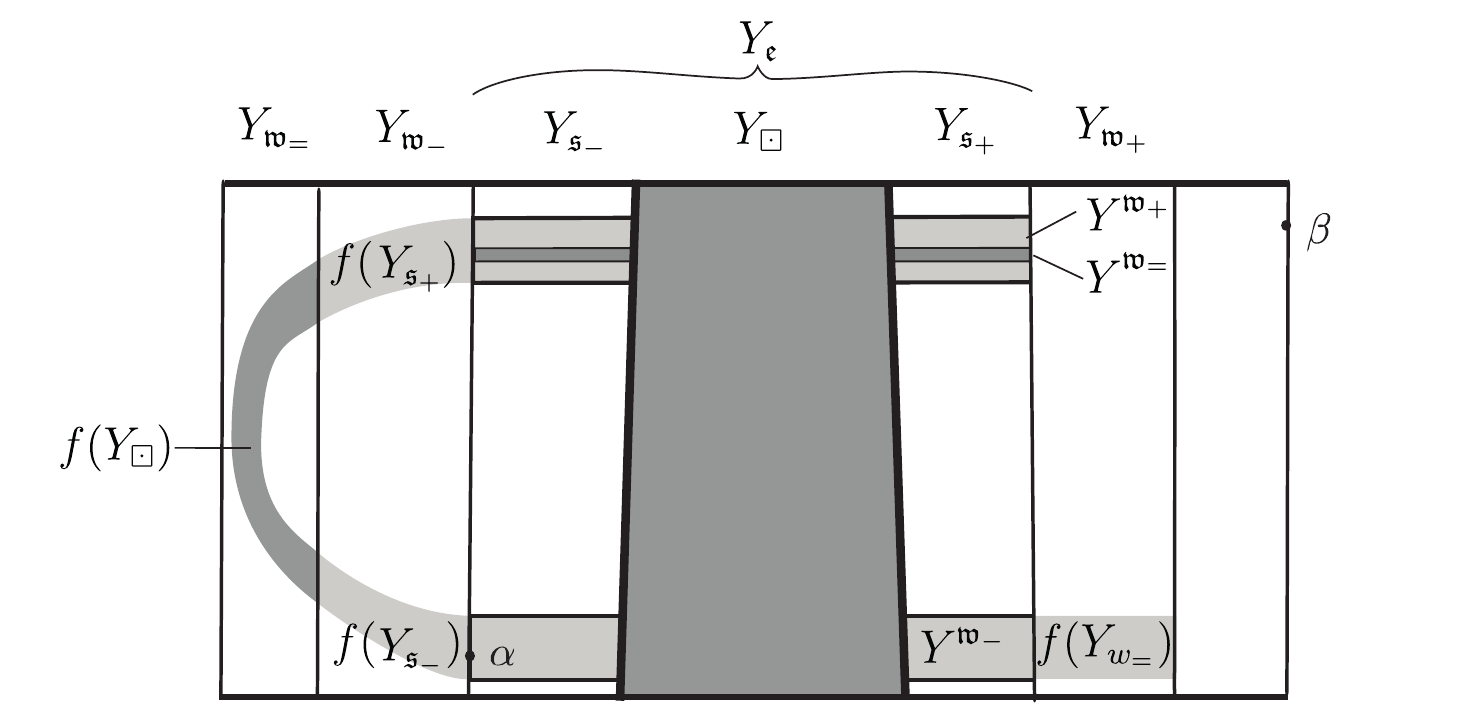}
\caption{ The boxes $Y_\sd$ and their images for $\sd\in \{\ss_-, s+, \sw_-, \sw_+, \sw_=, \boxdot\}$, for $a\in (a_2,a_1)$. 
}
\label{notation_zoologie_henon}\end{figure}
\subsection{Pieces for surface diffeomorphisms}\label{defCh}
Let us first state the general definitions of horizontal-vertical cones and of tame box. 

Given $c_h,c_v>0$ so that $c_h\cdot c_v< 1$, we define the horizontal cone $C_h$ and the vertical cone $C_v$ as:
\[ C_h:= \{(u_x,u_y): \|u_y\|\le c_h\| u_x\|\}\quad \text{and}\quad C_v := \{(u_x,u_y): \| u_x\| \le c_v \|u_y\|\}\; .\]

\begin{defi} 
A \emph{tame box} $Y$ is a box of the form 
$$Y= \{(x,y)\in \R^2 : \psi^{-}(x)\le y\le \psi^+(x)\qand \phi^{-}(y)\le x\le \phi^+(y)\}\; ,$$
where $\phi^{-},\phi^{+},\psi^{-}$ and $\psi^{+}$ are $C^r$-maps so that
$|\partial_t \phi^{-}|$ and $|\partial_t \phi^{+}|$ are at most $c_v$ and 
$|\partial_t \psi^{-}|$ and $|\partial_t \psi^{+}|$ are at most $c_h$ .

We put 
\[\partial^u Y= Y\cap \{(x,y): y\in \{\psi^{-}(x), \psi^+(x)\}\}\qand 
\partial^s Y= Y\cap \{(x,y): x\in \{\phi^{-}(y), \phi^+(y)\}\}\; .\]

\end{defi}
\begin{defi}
A \emph{piece} $\sd= (Y_\sd, n_\sd)$ is the data of a subset $Y_\sd\subset \R^2$ sent onto $Y^\sd:= f^{n_\sd}(Y_\sd)$ and such that:
\begin{itemize}
\item $Y_\sd $ and $Y^\sd$ are tame boxes. 
\item The map $f^{n_\sd}\colon Y_\sd \to Y^\sd$ satisfies the following cone condition:\end{itemize}
\begin{equation}\tag{$\mathcal C$}
Df^{n_\sd}(\R^2\setminus C_v)\subset C_h\quad \Leftrightarrow\quad (Df^{n_\sd})^{-1}(\R^2\setminus C_h)\subset C_v\; .
\end{equation}
\end{defi}

\begin{defi} A piece $\sd$ is \emph{hyperbolic} if the following condition holds true with $c:= \sqrt{c_vc_h}$. 
\begin{equation}
\tag{$H$} \forall z\in Y_\sd ,\; \forall (X_1,Y_1)=D_zf^{n_\sd} (X_0,Y_0),\quad\text{it holds} \quad |X_0|\le c |X_1|\qand |Y_1|\le c |Y_0|\; .\end{equation}
\end{defi}

Let us now consider these definitions in the context of a $C^r$-$\delta$-Hénon-like map $f$ with parameter $a<a_1$ and $b$ small. 

\begin{center}
\emph{ Let $\eta\in (0,\frac12\sqrt{a_1-a})$ and assume $\delta\le \eta/2$, and $|b|$ small compared to $\eta$.}
\end{center}
Then we observe that $Y_{s_-}$ and $Y_{s_+}$ are included in $\mathcal H_\eta:= 
Y_D\cap \{(x,y): |x|\ge \eta, |y|\le 3\}$. We fix $c_h$ and $c_v$ in function of $\eta$: 
\[c_h= 1/\eta\qand c_v= 1/(2\cdot c_h)=\eta/2\; .\]
\begin{defi}
A tame box $Y_\sd$ is a \emph{vertical strip} if $\partial^uY_\sd\subset\{y= \pm 3\}$.

A tame box $Y^\sd$ is a \emph{horizontal strip} if $\partial^s Y \subset \partial^s Y_\se$. 
\end{defi}

\begin{exam}\label{verticalstrip} For every $\sd\in \{\ss_-, \ss_+, \sw_-, \sw_+, \sw_=, \se\}$, the box $Y_\sd$ is a vertical strip.
\end{exam}

We are now ready to define the puzzle pieces for the Hénon-like map $f$ with $a<a_1$ and $b$ small as parameters.

\begin{defi}
A piece $\sd$ is a \emph{puzzle piece} if $Y_\sd$ is a vertical strip and $Y^\sd$ is a horizontal strip.
\end{defi}
The fact that $0< c_h\cdot c_v<1$ ensures that $C_h\cap C_v= \{0\}$. Hence, in particular $Df^{n_\sd}$ sends $C_h$ into $C_h$ and $Df^{- n_\sd}$ sends $C_v$ into $C_v$. This implies:

\begin{rema}\label{remacone} 
If $Y_\sd$ is a vertical strip, $Y_{\sd'}$ is a tame box and $n_\sd$ is such that
$f^{n_\sd}(\partial^s Y_\sd)\subset \partial^s Y_{\sd'}$ and $f^{n_\sd}|Y_\sd$ satisfies condition $(\mathcal C)$, then $(Y_\sd,n_\sd)$ is a piece. If moreover $\sd'=e$, then $\sd=(Y_\sd,n_\sd)$ is a puzzle piece. 
\end{rema}

 \begin{prop}\label{boxtame}
 For every symbol $\sd\in \{\ss_-, \ss_+, \sw_-, \sw_+, \sw_=, \se\}$, the associated pair $\sd:=(Y_\sd, n_\sd)$ is a puzzle piece.
 \end{prop}
\begin{proof}
In Example \ref{verticalstrip}, we already saw that $Y_\sd$ is a vertical strip. By Fact \ref{boxinclusion}, $f^{n_\sd}(\partial^s Y_\sd)\subset \partial^sY_\se$. By Remark \ref{remacone}, $\sd$ is a puzzle piece if the cone condition $(\mathcal C)$ is satisfied. Hence, it suffices to prove the cone condition. This is a consequence of the following Lemma, since $f^j (Y_\sd)$ is at least{ $2\eta$ }distant to $\{x=0\}$ for every $j\le n_\sd$:
\end{proof}
\begin{lemm}\label{coneta}
When $\delta\le \eta$ and $b$ is sufficiently small, for every $(x,y)\in \R^2$ such that $|x|\ge \eta$, the Hénon-like map $D_zf$ sends $\R^2\setminus C_v$ into $C_h$. \end{lemm}
\begin{proof}
Let $(u_x,u_y)$ be such that $\|u_y\|\le 1/c_v \|u_x\|= 2 c_h \|u_x\|$. Then for every $z\in Y_\sd$, it holds:
\[\partial_x f(x,y)= (2x +\partial_x \zeta(x,by), 1)\text{ and }\partial_y f(x,y) = (-b+b\partial_y \zeta(x,by), 0)\; .\]
 Hence with $(v_x,v_y)= Df(u_x,u_y)$, we have:
 \[v_x= (2x +\partial_x \zeta) u_x+ (-b+b\partial_y \zeta)u_y\text{ and }v_y= u_x\; .\]
 We recall that $|x|\ge \eta$, and so: 
 \[|v_x|\ge (2\eta-\delta)|u_x| -|b|(1+\delta)|u_y|\ge 
 (3\eta/2 -2|b|(1+\delta)/\eta)|u_x|=(3\eta/2 -2|b|(1+1/\eta))|v_y| \; .\]
As $b$ is small compared to $\eta$, it comes $|v_x|\ge \eta |v_y|= (1/c_h) |v_y|$.
\end{proof}

We extend here the definition of the $\star$-product from quadratic maps to Hénon-like maps:

\begin{prop}[Definition of the $\star$-product]\label{defstarhenon}
Let $\sd:= (Y_{\sd},n_{\sd})$ and $\sd':= (Y_{\sd'},n_{\sd'})$ be pieces so that $Y^\sd$ included in $Y_\se$ and such that $Y_{\sd'}$ is a vertical strip being between the two components of $\partial^s Y^\sd$. 
Then with $Y_{\sd\star \sd'}:= Y_{\sd}\cap f^{-n_{\sd}}(Y_{\sd'})$ and 
$n_{\sd\star \sd'}= n_{\sd}+n_{\sd'}$, the pair $\sd'= (Y_{\sd\star \sd'}, n_{\sd\star \sd'})$ is a piece. If moreover $\sd'$ is a puzzle piece, then $\sd\star \sd'$ is a puzzle piece.
\end{prop}
\begin{proof}
The cone property $(\mathcal C)$ for $ {\sd\star \sd'}$ is a direct consequence of the cone properties of $\sd$ and $\sd'$. We observe that $Y^{\sd}\cap Y_{\sd'}$ is a box bounded by two segments of $\partial^sY^{\sd}$ and two segments of $\partial^uY_{\sd'}$. Hence $Y_{\sd\star \sd'}:= f^{-n_{\sd}}(Y^{\sd}\cap Y_{\sd'})$ is a box, with $\partial ^u Y_{\sd\star \sd'}\subset \partial^u Y_\sd$, and $\partial ^s Y_{\sd\star \sd'}= (f^{n_\sd}|Y_\sd)^{-1}(\partial^s Y_{\sd'})$. As $C_v$ is sent by $D (f^{n_\sd}|Y_\sd)^{-1}$ into $C_v$, it follows that $Y_{\sd\star \sd'}$ is a vertical strip. 
By Remark \ref{remacone}, $Y^{\sd\star \sd'}$ is a tame box. Moreover, we note that if ${\sd'}$ is a puzzle piece then $\partial^s Y^{\sd\star \sd'}\subset \partial^s Y^{\sd'}$ is included in $\partial^s Y_\se$ and so $\sd\star \sd'$ is a puzzle piece. 
\end{proof}
\begin{defi} 
Given a $C^{d,r}$-family of surface dynamics $(f_p)_p$, a \emph{persistent piece} $\sd$ is the data for every $p$ of a piece $(Y_\sd, n_\sd)$ of $f_p$, so that $Y_\sd$ is bounded by $C^r$-curves $\partial^sY_\sd$ and $\partial^u Y_\sd$ which depend continuously of $p$.
\end{defi}
\subsection{Affine-like representation of pieces}
The affine-like representation of Palis-Yoccoz \cite{PY09} applies to (general) $C^r$-diffeomorphisms of $\R^2$. Let us fix the cone condition constants $0<c_h\cdot c_v<1$ for this section. 

\begin{prop}\label{ALsspara} Let $\sd=(Y_\sd,n_\sd)$ be a piece for a surface diffeomorphism $f$. Then, there are $C^r$-functions $A$ and $B$ satisfying for every $(x_0,y_0)\in Y_\sd$ and $(x_1,y_1)\in Y^\sd$:
\begin{equation}\tag{$\mathcal A$} f^{n_\sd}(x_0,y_0)= (x_1, y_1) \Leftrightarrow\left\{
\begin{array}{cc}
x_{0} =&A(x_1, y_0)\\
y_{1}=& B(x_1, y_0)\\
\end{array}\right. \end{equation}
\end{prop}
\begin{proof}
Let $ (x_0,y_0)\in Y_\sd$ and $(x_1,y_1)\in Y^\sd$. By the cone property, the image by $f^{n_\sd}$ of $\{y=y_0\}\cap Y_\sd$ is the graph of a function $x\mapsto B(x,y_0)$ intersected with $Y^\sd$. Likewise, the image by $f^{-n_\sd}$ of $\{x=x_1\}\cap Y^\sd$ is the transpose of the graph of a function $y\mapsto A(x_1,y)$ intersected with $Y_\sd$. 
We notice that the functions $(x_1,y_0)\mapsto A(x_1,y_0)$ and $(x_1,y_0)\mapsto B(x_1,y_0)$ are of class $C^r$. Observe also that by construction, the intersection point
$\{(x_1,B(x_1,y_0)\}= Graph\, B(\cdot ,y_0)\cap \{x=x_1\}$ is the image by $f^{n_\sd}$ of 
$\{(A(x_1,y_0),y_0)\}= \, ^t\! Graph\, A(x_1,\cdot)\cap \{y=y_0\}$ iff $f^{n_\sd}(x_0,y_0)= (x_1,y_1)$. 
\end{proof}
\begin{defi} The pair $(A,B)$ is the \emph{affine-like representation} of $\sd$. 
\end{defi}
An immediate consequence of the proof of Proposition \ref{ALsspara} is the following:
\begin{prop}Let $(f_p)_p$ is a $C^{d,r}$-family of surface mappings, with $d\le r$. Given any persistent piece $\sd$, the affine-like representations of $\sd$ form a $C^{d,r}$-family $(A_p,B_p)_p$.
\end{prop}

The affine-like representation was previously studied by Shilnikov \cite{Sh67} and its school (see for instance \cite{GST08}), and called \emph{the Shilnikov variable or cross map} in the $C^r$-smooth case. 

One of the interest of the affine-like representation is its times symmetry (times is reversed by swapping the variables $x$ and $y$). My main interest on this formalism are the very sharp distortion estimates we will obtained on the iterations of $f$, especially for the renormalizations problems. These bounds are obtained when we will look at the composition of hyperbolic pieces. 

A valuable fact already known from Shilnikov's work is that these variables give an effectively verifiable criterion of hyperbolicity:
\begin{lemm}[Lem. 3.2 \cite{PY01}]\label{lemYP}
The piece $\sd$ is \emph{hyperbolic} iff:
\[\frac{|\partial_x A|}{c} + \frac{|\partial_y A|}{c_v}\le 1\qand
\frac{|\partial_y B|}{c} + \frac{|\partial_x B|}{c_v}\le 1\quad \text{with } c=\sqrt{c_vc_h}\; .\]
\end{lemm}
Also the determinant is easily read in these variables:
\begin{fact}\label{determinant}
Under the hypothesis of Proposition \ref{ALsspara}, $\det\, D_{(x_0,y_0)}f^{n_\sd}= (\partial_y B/\partial_x A)(x_1,y_0)$. 
\end{fact}
\begin{proof}
We compute the differential $ D_{(x_0,y_0)}f^{n_\sd}$ and then we deduce its determinant: 
\[ dx_0= \partial_x A d x_1+ \partial_y A d y_0\Rightarrow 
dx_1 = (\partial_x A)^{-1} dx_0-(\partial_y A/\partial_x A)d y_0 
\]\[
 dy_1= \partial_x B d x_1+ \partial_y B d y_0 \Rightarrow 
 dy_1 = (\partial_x B/ \partial_x A)dx_0+(\partial_y B-\partial_x B \partial_y A /\partial_x A)dy_0
 \]\qedhere
\end{proof}

Let us now suppose that $f=f_p$ depends on a parameter $p$ in an open subset of $\R^k$, and that $(f_p)_p$ is a $C^{d,r}$-family of surface diffeomorphisms. 

The following expressions of the distortion bounds of pieces generalize those of given in \cite[\textsection3 and app. A]{PY09} in the case $d=2=r$, to any $d\le r< \infty$. Also a (new) special bound in the context of Hénon-like families is given (since these families are singular at $b\to 0$).

\begin{defi}[Distortion bound $\mathcal B^{d,r}$]
Let $\sd$ be a persistent piece for a family of maps $(f_p)_p$, and let $(A_p,B_p)_p$ be its affine-like representation. We define for $1\le d\le r<\infty$:
\[\begin{array}{rcl}
\mathcal B^{d,r}_0(\sd)&:=&\left\|(A_p, B_p)_p\right\|_{C^{d,r}},\\
\mathcal B^{d,r}_1(\sd)&:=&\left\|(D\log |\partial_x A_p|, D\log |\partial_y B_p|)_p\right\|_{C^{d',r-2}}\quad \text{with } d'= \min(d,r-2).
\end{array}\]
In the latter definition $D$ denotes the differential w.r.t. $x$ and $y$. 
The last bound regards the derivatives $\partial_p^k$:
\[\begin{array}{rcl}\mathcal B^{d}_0(\sd)&:=&\frac1{n_\sd}\max_{ 1\le k\le \min(d,r-1)}\left\| \partial_p^k\log |\partial_x A_p|,\partial_p^k\log |\partial_y B_p|\right\|_{C^0}\end{array}.\]
In the case where $(f_p)_p$ is a Hénon-like family $(f_{a\, b})_{(a, b)}$, we will see in section \ref{Affine-like Henon-like} that $ \partial_p \log |\partial_y B_p|$ is not bounded: an extra $\partial_p \log b$ appears. That is why, for every Hénon-like family of multiplicity $m\ge 1$, we will consider instead of $\mathcal B^{d}_0$ the following:
\[\mathcal B^{d}_{m}(\sd):=\frac1{n_\sd}\max_{ 1\le k\le \min(d,r-1)}\left\| \partial_p^k\log |\partial_x A_p|,\partial_p^k\log |\partial_y B_p/b^{m n_\sd}|\right\|_{C^0} ,\; \text{with } p=(a,b)\in \R^{k-1}\times \R.\]

We denote by $\mathcal B^{d,r}(\sd)$ the maximum of 
$\mathcal B^{d,r}_0(\sd)$ and 
$\mathcal B^{d,r}_1(\sd)$.
\end{defi}

The following improves several bound estimates from the Shilnikov school in the context of a hyperbolic set of surface maps (compare for instance with \cite[Lem. 7]{GST08}). Up to my knowledge, the bounds given by $(\mathcal B_1^{d,r})$ for $i+j+k=r-1$ are new when $r\ge 3$. This will be crucial to perform $C^{d,r}$-Hénon-like and Swallow-like renormalizations without losing derivatives. In the $C^2$-topology, the equivalent of the following proposition is \cite[Prop. 16 P.60]{PY09}. In the next subsection, we will give technical consequences of the bounds $\mathcal B_0^{d,r}$ and $\mathcal B_1^{d,r}$.
 
\begin{prop}\label{Birkhoff}
For every $(c_v,c_h)$, $1\le d\le r$, for every $K_0>0$, there exists $K_1\ge 1$ such that for every $C^r$-family of diffeomorphisms $(f_p)_{p}$ of $\R^2$, for every $m\ge 1$, for every sequence of hyperbolic pieces $(\sd_j)_{1\le j\le N}$, with $\sd:=\sd_1\star \sd_2 \cdots \star \sd_N$, it holds:
\begin{itemize}
\item if $\mathcal B^{d,r}_0(\sd_j)\le K_0$ for every $j$, then $\mathcal B^{d,r}_0(\sd)\le K_1$,
\item if $\mathcal B^{d,r}_0(\sd_j), \mathcal B^{d,r}_1(\sd_j)\le K_0$ for every $j$, then $\mathcal B^{d,r}_1(\sd)\le K_1$,
\item if $\mathcal B^{d,r}_0(\sd_j), \mathcal B^{d,r}_1(\sd_j), \mathcal B^{d}_m (\sd_j)\le K_0$ for every $j$, then $\mathcal B^{d}_m(\sd)\le K_1$.
\end{itemize}\end{prop}

The point of the above proposition is that $K_1$ depends only on $K_0,d,r, c_v$ and $c_h$. Thus $K_1$ does not depend on $N$, nor on $\|(f_p)_p\|_{C^{d,r}}$ nor on $\|(f^{-1}_p)_p\|_{C^{d,r}}$. 
Before proving this proposition let us simplify the set up. 
\begin{fact}
Up to the affine change of coordinates $(X',\frac{c}{c_v}Y')= (X,Y)$ with $c=\sqrt{c_h c_v}$, we can assume ${c}=c_v=c_h<1$.
\end{fact} 
\begin{proof} After the change of coordinates, we notice that condition $(H)$ is still satisfied. 
Also we recall that given $Df^{n_\sd}(X_0,Y_0)=(X_1,Y_1)$, the cone condition $(\mathcal C)$ is:
\[c_v|Y_0|< |X_0|\Rightarrow |Y_1|\le c_h|X_1| \qand 
c_h|X_1|\le |Y_1|\Rightarrow |X_0|\le c_v|Y_0|\; .\]
Hence 
 \[c_v \frac{c}{c_v} |Y'_0|< |X'_0|\Rightarrow |Y'_1|\le c_h\frac{c_v}{c} |X'_1| \qand 
c_h |X'_1|\le \frac{c}{c_v} |Y'_1|\Rightarrow |X'_0|\le c_v\frac{c}{c_v} |Y'_0|\; .\]
And we conclude thanks to the equality $c_vc_h=c^2$.
\end{proof}
Hence we assume ${c}=c_v=c_h<1$. This changes $K_0$ by a function depending only on $K_0$, $c_v$, $c_h$. For every $i$, let us denote by $(A_{i\, p},B_{i\, p})$ the affine-like representation of the piece $\sd_i$ of $f_p$. 

 An immediate consequence of Lemma \ref{lemYP} is:
\begin{equation}\tag{$H'$}
{|\partial_x A_{i\, p}|} + {|\partial_y A_{i\, p}|}\le c< 1\qand
{|\partial_y B_{i\, p}|}+ {|\partial_x B_{i\, p}|}\le c< 1\; .\end{equation}

For $(x_0,y_0)\in Y_\sd$, we define inductively $(x_i,y_i)$ by $(x_{i+1},y_{i+1})=f^{n_{\sd_{i+1}}}_p(x_i,y_i)$. By definition of the affine-like representation, it holds:
\begin{equation}\tag{$\star$} \left\{\begin{array}{c}
x_i= A_{i\, p}(x_{i+1}, y_i)\; ,\\
y_{i+1}= B_{i\, p}(x_{i+1},y_i)\; .\\
\end{array}\right. \end{equation}
The bound $\mathcal B_0^{d,r}$ on $\sd$ is given by the following:
\begin{lemm}\label{prebirkof1}
There is $K_0' >0$ depending only on $K_0:= \max_i \|(A_{i\, p}, B_{i\, p})_p\|_{C^{d,r}}$ and there are $C^{d,r}$-families of functions $(X_{i\, p},Y_{i\, p})_p$, 
 with $C^{d,r}$-norms bounded by $K_0'$ and which satisfy:
\[(x_i,y_i)=(X_{i\, p},Y_{i\, p})(x_N,y_0)\quad .\]
\end{lemm}
\begin{proof}
The system $(\star)$ invites us to consider the function:
\[\Psi_p:((x_N,y_0), (x_i,y_{i+1})_{0\le i\le N-1})\mapsto
\left((x_i,y_{i+1})- (A_i,B_i) (x_{i+1},y_{i})\right)_{0\le i\le N-1}\]

We endow $(\R^2)^N$ with the uniform norm:
\[\|(x_i,y_{i+1})_{0\le i\le N-1}\|= \max_{0\le i\le N-1} (|x_i|,|y_{i+1}|)\]
Then the family of functions $(\Psi_p)_p$ is of class $C^{d,r}$ with norm bounded in function of $K_0$. Also, by $(H')$, its derivative w.r.t. 
$(x_i,y_{i+1})_{0\le i\le N-1}$ is equal to the identity plus a contraction (by a factor $c$). Hence it is invertible, and the norm of its inverse is bounded independently of $N$. 
Then the implicit function theorem concludes the proof.
\end{proof}
To show the bounds on $\mathcal B_1^{d,r}$ and $\mathcal B_m^{d}$, we shall prove first:
\begin{lemm}\label{prebirkof2}Let $K_0:= \max_i \|(A_{ip},B_{ip})_p\|_{C^{d,r}}$. 
For every $\sigma\in (c,1)$, there is $C=C(K_0,c,\sigma)>0$ independent of $N$ such that for every $0\le i\le N$, with $D$ the differential w.r.t. $x_N$ and $y_0$, and with $\mu(i):=\min(i,N-i)$, it holds:
\[\| (D X_{i\, p}, DY_{i\, p})_{p}\|_{C^{\bar d,r-1}}\le C(c/\sigma)^{\mu(i)}\, \quad \text{with }\bar d:= \min (r-1, d).
\]
\end{lemm}
\begin{proof}
Put $X_{i\, p}':= X_{i\, p}/\sigma^{\mu(i)}$ and $Y_{i\, p}':= Y_{i\, p}/\sigma^{\mu(i)}$. We note that:
\[\left\{\begin{array}{c}
X_{i\, p}'= \sigma^{-\mu(i)}A_{i\, p}(\sigma^{\mu(i+1)}X'_{i+1\, p}, \sigma^{\mu(i)} Y'_{i\, p})=\sigma^{-\mu(i)}\cdot A_{i\, p}(X_{i+1\, p}, Y_{i\, p}) \\
Y'_{i+1\, p}= \sigma^{-\mu(i+1)}B_{i\, p}(\sigma^{\mu(i+1)}X'_{i+1\, p},\sigma^{\mu(i)}Y'_{i\, p})= \sigma^{-\mu(i+1)}\cdot B_{i\, p}(X_{i+1\, p},Y_{i\, p})\\
\end{array}\right. \]
With $\epsilon_i= \mu(i+1)-\mu(i)\in \{-1,0,1\}$, it holds:
\[
\left\{\begin{array}{c}
D X'_{i\, p}= 
\sigma^{\epsilon_i} \partial_x A_{i\, p}(X_{i+1\, p}, Y_{i\, p}) \cdot D X'_{i+1\, p} 
+ \partial_y A_{i\, p}(X_{i+1\, p}, Y_{i\, p}) \cdot D Y'_{i\, p} 
\\
D Y'_{i+1\, p}=
\partial_x B_{i\, p}(X_{i+1\, p}, Y_{i\, p}) \cdot D X'_{i+1\, p} 
+ \sigma^{-\epsilon_i}\partial_y B_{i\, p}(X_{i+1\, p}, Y_{i\, p}) \cdot D Y'_{i\, p}\\
\end{array}\right.
\]
The above system invites us to consider the following operator of $C^{r-1}$-functions from the domain of the implicit representation of $\sd$ into 
$(\mathcal L_2(\R))^N$, with $\mathcal L_2(\R^2)$ 
the space of linear maps of $\R^2$:
\[\Phi_p: \left(\begin{array}{c}L_{X\, i}\\ L_{Y\, i+1}\end{array}\right)_{0\le i\le N-1}\mapsto
\left(\begin{array}{c}
\sigma^{\epsilon_i} \partial_x A_{i\, p}(X_{i+1\, p}, Y_{i\, p}) \cdot L_{X\, i+1} 
+ \partial_y A_{i\, p}(X_{i+1\, p}, Y_{i\, p}) \cdot L_{Y\, i} 
\\
\partial_x B_{i\, p}(X_{i+1\, p}, Y_{i\, p}) \cdot L_{X\, i+1} 
+ \sigma^{-\epsilon_i}\partial_y B_{i\, p}(X_{i+1\, p}, Y_{i\, p}) \cdot L_{Y\, i}
\end{array}\right)_{0\le i\le N-1},
\]
with $L_{X\, N}:=D X_N'(x_N,y_0)=D X_N(x_N,y_0)$ and 
 $L_{Y\, 0}:=D Y_0'(x_N,y_0)=D Y_0(x_N,y_0)$. 
 
The family of functions $(\Phi_p)_p$ is of class $C^{\bar d, r-1}$ with bounded $C^{\bar d, r-1}$-norm by Lemma \ref{prebirkof1}. Also the differential of each $\Phi_p$ is $(c/\sigma)$-contracting. Thus by the implicit function theorem, the fixed point $(D X_{i\, p}',D Y_{i+1\, p}')_i$ is a function of $(x_n,y_0)$ of class $C^{r-1}$. Moreover as $(\phi_p)_p$ is of class $C^{\bar d, r-1}$, each family $(D X_{i\, p}',D Y_{i+1\, p}')_p$ is of class 
$C^{\bar d,r-1}$ with uniformly bounded $C^{\bar d,r-1}$ norm. 

Consequently, each $(D X_{i\, p},D Y_{i\, p})_p$ is a $C^{\bar d, r-1}$-$\sigma^{\mu(i)}$-bounded family of functions of $(x_N,y_0)$.
\end{proof}
\begin{proof}[Proof of Proposition \ref{Birkhoff}]
With the above notations, it holds:

\[\left\{\begin{array}{c}
x_0= A_{0\, p}(A_{1\, p}(\cdots A_{N-1\, p}(x_{N}, y_{N-1})\cdots, y_2), y_1)=A_p(x_{N},y_0)\; ,\\
y_{N}= B_{0\, p}(x_{1}, B_{1\, p}(x_2,\cdots B_{N-1\, p}(x_{N},y_0)\cdots))=B_p(x_N,y_0)\; .\\
\end{array}\right. \]

Consequently:
\[\left\{\begin{array}{cc}
\log \partial_{x_{N}} A_p =& 
\sum \log \partial_x A_{i\, p} (X_{i+1\, p},Y_{i\, p})\; , \\
\log \partial_{y_{0}} B_p= &\sum \log \partial_y B_{i\, p} (X_{i+1\, p},Y_{i\, p})\; .\\
\end{array}\right. \]
By time-symmetry of $\mathcal B_1^r$, we study only the bounds on $A$: 

\[D \log |\partial_{x_{N}} A_p| =
\sum D (\log |\partial_x A_{i\, p} |) (DX_{i+1\, p},DY_{i+1\, p})\; .\]

This sum is $C^{d',r-2}$-uniformly bounded since each if its terms is a product of a bounded term with an exponentially small one. This proves the bound $\mathcal B^{d,r}_1$. 

To obtain the bound $\mathcal B_m^{d}$ we observe that the following sum has its $i^{th}$-term dominated by $n_{\sd_i}$:
\[ \partial_{p}\log |\partial_{x_{N}} A_p| = 
\sum (D\log |\partial_x A_{i\, p} |)(\partial_p X_{i+1\, p},\partial_p Y_{i\, p})+
(\partial_p \log |\partial_x A_{i\, p} |)(X_{i+1\, p},Y_{i\, p})\; .\]
Although, we cannot use the time symmetry for the remaining bound on $(B_p)_p$ , the proof follows from the same observation:
\[ \begin{split} \partial_{p}\log |\partial_{y_0} B_p/b^{m\cdot n_\sd}| = 
\sum (D\log |\partial_y B_{i\, p} /b^{m n_{\sd_i}}|)(\partial_p X_{i+1\, p},\partial_p Y_{i\, p})+
(\partial_p \log |\partial_y B_{i\, p} /b^{ m n_{\sd_i}}|)(X_{i+1\, p},Y_{i\, p}).\end{split}\]
\end{proof}
\subsection{Useful consequences of the bounds $\mathcal B_0^{d,r}$ and $\mathcal B_1^{d,r}$
}
Let $(f_p)_{p}$ be a $C^{d,r}$-family of surface mappings with a piece $\sd$. 
Let $p\mapsto (x_0(p),y_0(p))$ be a $C^d$-function with image in $Y_\sd(p)$, and let $(x_1(p),y_1(p)) \in Y^\sd (p)$ be equal to $f^{n_\sd}(z_0(p))$. We put $z(p):= (x_1(p),y_0(p))$. 
Let $(A_p,B_p)_p$ be the affine-like representation of $\sd$.
We put also:
\[\left\{\begin{array}{c}
\breve A_p(x,y)= A_p(z(p)+(x, y))- A_p(z(p)+(0, y))\\
\breve B_p(x,y)= B_p(z(p)+(x, y))- B_p(z(p)+(x, 0))\end{array}\right.\]
\[\sigma_p:= \partial_x \breve A_p(0)= \partial_x A_p(z(p))\qand \lambda_p :=\partial_y\breve B_p(0)=\partial_y B_p(z(p)) \; .\]

\begin{lemm}\label{distorsionutil}
Let $d\le r-2$and $R>0$. If $\epsilon >0$ is such that $\epsilon R\cdot \mathcal B_1^{d,r}(\sd)(1+ \|p\mapsto z(p)\|_{C^d})^d$ is small, then the functions:
\[\Pi_{A\, p}:(x,y)\in [-R,R]^2\mapsto \frac{\breve A_p(\epsilon\cdot x, \epsilon \cdot y)}{\epsilon\cdot \sigma_p}\qand \Pi_{B\; p}: (x,y,p)\mapsto \frac{\breve B_p(\epsilon \cdot x, \epsilon \cdot y)}{\epsilon\cdot \lambda_p }\]
form families $(\Pi_{A\, p})_p$ and $(\Pi_{A\, p})_p$ which 
are $C^{d,r}$-close to respectively the first and the second coordinate projections. 
\end{lemm}
\begin{proof}
By time symmetry, we just need to prove that $(\Pi_{A\, p})_p$ is $C^{d,r}$-close to $((x,y)\mapsto x)_p$. First we observe that $\Pi_{A \, p}(0,y)=0$ for every $y$ and $p$. Hence we only need to show that $(\partial_x \Pi_{A\, p})_p$ is $C^{d, r-1}$-close to $((x,y)\mapsto 1)_p$. 
We note that $\partial_x \Pi_{A\, p}(0)= 1$ for every $p$. Thus:
$$\partial_x \Pi_{A\, p}(z) = \exp\left(\int_0^1 \partial_t (\log |\partial_x \Pi_{A\, p}|(tz)) dt\right)=
\exp\left(\int_0^1 D_{t \cdot z} (\log |\partial_x \Pi_{A\, p}|) (z) dt\right).$$
Put $z=(z_x,z_y)$. We conclude the proof by noting that $(D_{tz}(\log |\partial_x \Pi_{A\, p}|)(z))_p$ is equal to $ \epsilon (D_{z(p)+tz} \log \partial_x A_p)(z)$
 and so $C^{d,r-2}$-dominated by $\epsilon \cdot \mathcal B_1^{d,r}(\sd)(1+ \|p\mapsto z(p)\|_{C^d})^d\cdot R$ which is small.
 \end{proof}

A coarse consequence of $\mathcal B_1^{0,r}$ is the following:
\begin{lemm}\label{distorsionutil0}
Let $K:= \mathcal B^{0,r}_1(\sd)(\diam Y_\sd+\diam Y^\sd)(1+ \|p\mapsto z(p)\|_{C^d})^d$. 

Then for every $p$, the functions $z\mapsto \log |\frac{\partial_x\breve A_p}{\sigma_p}|(z)$ and $z\mapsto \log|\frac{ \partial_y\breve B_p}{\lambda_p}|(z)$ are $C^{r-1}$-dominated by $K$. 
Also, for every $(x,y)$, it holds $|{\partial_y\breve A_p}(x,y)| 
\le K' \sigma_p |x|$ and $|\partial_x \breve B_p(x,y)|\le K' \lambda_p |y|$, with $K'$ depending only on $K$. 
\end{lemm}
\begin{proof}
By time symmetry, we just need to prove the bound on $\breve A_p$. 
As $\partial_x \breve A_p$ is equal to the composition of $\partial_x A_p$ with a translation by $z(p)$, the $C^{r-1}$-norm of the derivatives of $z\mapsto \log|\partial_x \breve A_p|(z)$ are dominated by $\mathcal B_1^{0,r}(\sd)(1+ \|p\mapsto z(p)\|_{C^d})^d$. Thus:
 \begin{equation}
 \left\|\log\left|\frac{\partial_x \breve A_p }{\sigma_p}\right|\right\|_{C^{r-1}}\le K \qand 
 \left\|\partial_x \breve A_p\right\|_{C^{r-1}}\le K' |\sigma_p|\end{equation}
 \begin{equation}
 \partial_y \breve A_p =\int_{t=0}^1 x\cdot \partial_x \partial_y \breve A_p(tx, y) dt
\Rightarrow \|\partial_y \breve A_p\|_{C^{r-1}}
 \le K' |\sigma_p|\cdot |x|\; .\end{equation}
\end{proof}
{
 \begin{prop}\label{continuityinprod}
Let $d<r$. Let $\sd=\sd'\star \sd''$ be a piece equal to the $\star$-product of two hyperbolic pieces $\sd'$ and $\sd''$. Let $(A_p,B_p)$, $(A'_p,B'_p)$ and $(A''_p,B''_p)$ be the affine-like representations of $\sd,\sd'$ and $\sd''$ respectively. For some fixed real numbers $x_2$ and $y_0$, let $p\mapsto (x_0(p), x_1(p),y_1(p),y_2(p))$ be a $C^d$ function so that 
$(x_0(p),y_0)\in Y_{\sd}$ is sent by $f_p^{n_\sd}$ to $(x_1(p),y_1(p))\in Y^{\sd'}\cap Y_{\sd''}$ which is sent in turn by $f^{n_{\sd'}}_p$ to $(x_2,y_2(p))\in Y^{\sd}$.

Then there exists $K_1$ 
 depending only on $c_v\cdot c_h$, $\|(A''_{ p},B'_p)_p\|_{C^{d,r-1}}$ and $\mathcal B^{d,r}_1(\sd')$ , and there exists 
 $K_2$ depending only on 
$c_v\cdot c_h$, $\|(A''_{ p},B'_p)_p\|_{C^{d,r-1}}$ and $\mathcal B^{d,r}_1(\sd'')$ 
 such that:\begin{enumerate}
\item $\| (y\mapsto A_p(x_2,y)- A'_p(x_1(p),y))_p\|_{C^{d,r-1}}\le K_1 \|(\partial_x A'_p)_p\|_{C^{d,r-1}}$ ,
\item $\|(x\mapsto B_p(x,y_0)- B''_p(x,y_1(p)))_p\|_{C^{d,r-1}}\le K_2 \|(\partial_y B''_p)_p\|_{C^{d,r-1}}$.
\end{enumerate}
\end{prop}
\begin{proof}
First we recall that the (transverse) intersection point $(X_{1\, p},Y_{1\, p})$ of $\{(x, B'_p(x,y)): x\}$ with $\{(A''_p(x_2,y),y):y \}$ is a $C^{d,r}$-family of functions of $(x_2,y)$. Moreover $(y\mapsto X_{1\, p} (x_2,y))_p$ has its $C^{d,r-1}$-norm bounded by a function of the one of $(A''_{ p}, B'_{ p})_p$ and of $c_v\cdot c_h$. 
Also $A_p(x_2,y)= A'_p(X_1(x_2,y),y)$. Consequently the first statement follows from the following equality:
\[A_p(x_2,y)= A'_p(x_1(p),y) + \int_{x_1(p)}^{ X_{1\, p}(x_2,y)} \partial_x A'_p(s,y)ds\; .\]
The second statement is obtained by time-symmetry. 
\end{proof}} 
 
\subsection{Affine-like representation around uniformly hyperbolic set of Hénon-like maps}\label{Affine-like Henon-like}
Among $C^{d,r}$-Hénon-like families
 $(f_{p})_{I\times J}$, let us explain how the distortion bounds of persistent ``uniformly hyperbolic pieces" are independent of $\hat b:= \max_J |b|$ small and the multiplicity $m$ of the family.
We recall that we defined for $\eta>0$:
$$\mathcal H_\eta:= \{(x,y)\in Y_D: |x|\ge \eta \qand |y|\le 3\}\qand c_h=1/\eta,\; c_v=\eta/2\; .$$
Let $C\ge 1$ and $\kappa<1$. We will assume below the following hyperbolicity condition:
\begin{equation}\tag{$H_\eta$} 
 \| D_zf^j_{p} (u)\|\ge C\cdot \kappa^{-j}\cdot \|u\|\quad ,\quad \forall u\in C_h\; , \quad \forall z\in \bigcap_{k=0}^{j-1} f^{-k}_{p}(\mathcal H_\eta)\; , \quad \forall p\in I\times J\; .
 \end{equation}

\begin{prop}\label{unifhyp}
Let $\kappa>1$ and $C>0$. There exist $b_0=b_0(\eta)>0$ and $K_1=K_1(\kappa, \eta, C)$, such that the following property holds true for every $\delta\le \eta/2$ and for every $\delta$-$C^{d,r}$-Hénon-like family $(f_{p})_{I\times J}$ of multiplicity $m$ satisfying $(H_\eta)$ and $J\subset [-b_0,b_0]$.

If a piece $\sd$ satisfies $f^k_{p}(Y_\sd)\subset \mathcal H_\eta$ for every $k<n_\sd$, then
the distortions
$\mathcal B_0^{d,r}(\sd)$, $\mathcal B_1^{d,r}(\sd)$, $\mathcal B_m^d(\sd)$ are bounded by $K_1$.
\end{prop}
\begin{proof}We recall that by Lemma \ref{coneta}, if $b_0$ is sufficiently small in function of $\eta$, then it holds:
\begin{equation} \tag{$\mathcal C'$}
\forall p\in I\times J ,\quad \forall z\in \mathcal H_\eta,\quad Df_{p}(\R^2\setminus C_v)\subset C_h\; .
\end{equation}
Consequently the following is an immediate consequence of $(H_\eta)$ and $(\mathcal C')$: 
\begin{fact}\label{prehyperbolic}
Let $\sd$ be a piece of order $n_\sd\ge N_0:=\frac{\log(c/C)}{\log \kappa}$. If $f_{p}^k(Y_\sd)\subset \mathcal H_\eta$ for every $k<n_\sd$, then $\sd$ is hyperbolic. 
\end{fact}
Therefore, by Proposition \ref{Birkhoff}, Proposition \ref{unifhyp} is a direct consequence of the next Lemma.\end{proof}
 \begin{lemm}\label{distorsionordrefini}
For $b_0$ sufficiently small compared to $\eta$, there exists $K_0$ depending only on $\eta$ and $ N_0$ so that every piece $\sd$ of order $\le N_0$ and such that $f^k_{p}(Y_\sd)\subset \mathcal H_\eta$ for every $k< N_0$, the distortions
$\mathcal B_0^{d,r}(\sd)$, $\mathcal B_1^{d,r}(\sd)$, $\mathcal B_m^d(\sd)$ are bounded by $ K_0$. 
\end{lemm}
\begin{proof}[Proof of Lemma \ref{distorsionordrefini}]
We notice that the piece $\sd$ is the $\star$-product of at most $N_0$ pieces $\sd_i$ of order $1$ satisfying $Y_{\sd_i}\subset \mathcal H_\eta$. By the transversality given by the cone condition and 
since all the sums in the proof of Proposition \ref{Birkhoff} are made by a number of terms bounded by $N_0$, it suffices to show that each piece $\sd_i$ displays bounds $\mathcal B^{d,r}_0(\sd_i)$, $\mathcal B^{d,r}_1(\sd_i)$ and ${\mathcal B}^{d}_m(\sd_i)$ depending only on $\eta$.

 We observe that the affine-like representation of the piece $\sd_i$ is the solution of:
\[(A_{i\, p}(x_1,y_0)= x_0\qand B_{i\, p}(x_1,y_0)= y_1)\Leftrightarrow 
(x_1=x_0^2+a+\zeta_{p}(x_1,b^my_0)-b^my_0\qand y_1=x_0)\; ,\]
where $f_p(x,y)=(x^2+a-b^m y +\zeta_{p}(x,b^my), x)$ at the parameter $p=(a,b)$. Note in particular that $A_{i\, p}=B_{i\, p}$. Also $A_{i\, p}(x_1,y_0)=\hat A_{i\, p}(x_1,b^m y_0)$ where $\hat A_{i\, a}$ is the one of the two implicit solutions of:
 \[x_1=x_0^2+a+\zeta_{p}(x_0,y_0)-y\Leftrightarrow \pm \hat A_{i\, p}(x_1,y_0)= x_0 \; .\]
 By definition of $(\mathcal H_\eta)$, the bounds $\mathcal B_0^{d,r},\mathcal B_1^{d,r}$ and $\mathcal B_2^{d}$ of $(\hat A_{i\, p}, \hat A_{i\, p})$ depend only on $\eta$. As $(A_{i\, p}, B_{i\, p})$ is the composition of $(\hat A_{i\, p}, \hat A_{i\, p})$ with a non-expending map $(x,y)\mapsto (x,b^m\cdot y)$, the $C^r$-norm $\mathcal B_0^{d,r}(\sd_{i})$ depends only on $\eta$. To obtain the bounds $\mathcal B_1^{d,r}(\sd_{i})$ and ${\mathcal B}^{d}_m$, we compute:
\[\left\{\begin{array}{cc}
\log |\partial_{x } A_{i\, p}|(x,y) =& \log |\partial_x \hat A_{i\, p} (x,b^m y)|\\
\log |\partial_{y} B_{i\, p}/b^m|(x,y)= &\log |\partial_y \hat A_{i\, p} (x,
b^m y)|\\
\end{array}\right. \]
 Hence, for $b_0=b_0(\eta)$ sufficiently small, the bounds $\mathcal B_1^{d,r}$ and ${\mathcal B}^{d}_m$ on $\sd_i$ are given by the 
 $C^{\min(d,r-1),r-1}$-norm of $x\in [-2,2]\cap \hat A_{i\, p}^{-1}([-2,2]\setminus (-\eta,\eta))\mapsto (\log |\partial_x \hat A_{i\, p} (x,
0)|,\log |\partial_y \hat A_{i\, p} (x,
0)|)$ which is bounded by a function of $\eta$ independent of $\delta\le \eta/2$ and $p\in I\times J$. 
\end{proof}
The following will be useful:
\begin{prop}\label{bsmall} Under the assumptions of Proposition \ref{unifhyp}, with $(A_{p}, B_{p})_{I\times J}$ the affine-like representation of $\sd$, the families $(\partial_y A_{a\, b})_{I\times J}$ 
and $(\partial_y B_{a\, b})_{I\times J}$ are $O(\max_J |b|^m)$-$C^{\bar d, r-1}$-small whenever $n_\sd\ge 1$, with $\bar d= \min(d,r-1)$.
\end{prop}
\begin{proof}
With the notations of the previous proof, for $n_\sd=1$, this proposition is a consequence of the fact that $A_{i\, p}(x,y)=B_{i\, p}(x,y)=\hat A_{i\, a}(x,b^my)$. For $n_\sd> 1$, with $p=(a,b)$, the derivative $\partial_{y_0} \psi_p$ defined in Lemma \ref{prebirkof1} is dominated by $b^m$. Thus the functions $\partial_y X_{1\, p}$ and $\partial_y Y_{n_{\sd}\, p}$ in the statement of Lemma \ref{prebirkof1} are $C^{d,r-1}$-bounded by $\max_J |b|^m$. As $ A_{p}= X_{0\, p}$ and $ B_{p}=Y_{n_{\sd}\, p}$, the proposition follows. 
\end{proof}

\section{Renormalization Theorems \ref{Charthenonlike} and \ref{Charthenonlike2}}

\subsection{Hénon-like renormalization}
Let $0\le d\le r-3$. Let $m\ge 1$. 
Let $(f_{p})_{I\times J}$ be a two-parameters family of surface diffeomorphisms (not necessarily Hénon-like). Let $\sd$ be a piece with affine-like representation $(A_{p},B_{p})_{I\times J}$. Let us denote the parameter $p=:(a,b)$.

Let $p\mapsto c(p)$ be a $C^d$-function, and let us define
$z_0(p)= (A_p(c,c),c)(p)$, $z_1(p)= (c, B_p(c,c))(p)$ and:
$$V_p(y):= A_p(c(p),c(p)+y)\qand 
W^s_{p}:= f_{p}^{-n_\sd} (\{x=c(p)\}\cap Y^\sd)= \{(A_p(c(p), y),y):y\}\, ,
$$
$$H_p(x):= B_p(x+c(p), c(p)) \qand W^u_{p}:= f_{p}^{n_\sd} (\{y=c(p)\}\cap Y_\sd) = \{(x, B_{p}(x, c(p))):x\}
\; .$$

Assume the existence of a $C^d$-function $b\in J\mapsto a(b)\in I$ so that with $p(b):=(a(b),b)$ the following properties hold:
\begin{enumerate}[$(i)$]
\item Restricted to a neighborhood of $(x,y)\in \{z_0(p): p\}$ and $p\in \{p(b):b\}$ the family $(f_p)_p$ is 
 of the form $f_p(x,y)=(g_p(x,y),x)$.
\item For every $b$, the curve $W^u_{p(b)}$ is sent by $f_{p(b)}$ to a curve tangent to the curve $W^s_{p(b)}$ at the point $f_p(z_1(p))=f_{p}(c(p), H_p(0))= (V_p(0),c(p))$. 
Moreover, for every $p$ close to $\{p(b):b\}$, $x=0$ is a critical point of the distance $x\mapsto \|f_{p}(c(p)+x, H_p(x)) - (V_p(x),c(p)+x)\|$. In other words:

\begin{equation}\tag{$T$}\label{Defrp}
g_{p}(c(p)+ x, y+H_p(x))-V_p(x)=
 \mu_p+ q_p x^2-d_p \cdot y+r_p(x, d_p \cdot y) \; ,\end{equation}
with $ d_p:= \det \, D_{z_1(p)} f$, 
$p\mapsto (\mu_p, q_p)$ a $C^{d}$-function and $(r_p)_p$ a $C^{d,r}$-family of functions satisfying:
\[
0=r_p(0)= \partial_y r_{p}(0)= \partial_x r_{p}(0)=\partial_x^2 r_p (0),\quad \forall p \qand 
\mu_{p(b)}(0)=0\; ,\quad 
 \forall b\; .\]
 \item The tangency is quadratic and the unfolding of this tangency is non-degenerated: 
 \begin{equation}\tag{$ND$} \exists M\ge 1: \forall b\in J, \;\text{ at } p=p(b) \; , \quad q_p\not=0\qand 
\det \partial_p (\mu_p ,\sqrt[M]{ \det D_{z_0} f^{n_\sd +1}_p})\not = 0
 .\end{equation} 
 
\end{enumerate}
\begin{figure}[h!]
 \centering
 \includegraphics[width=7cm]{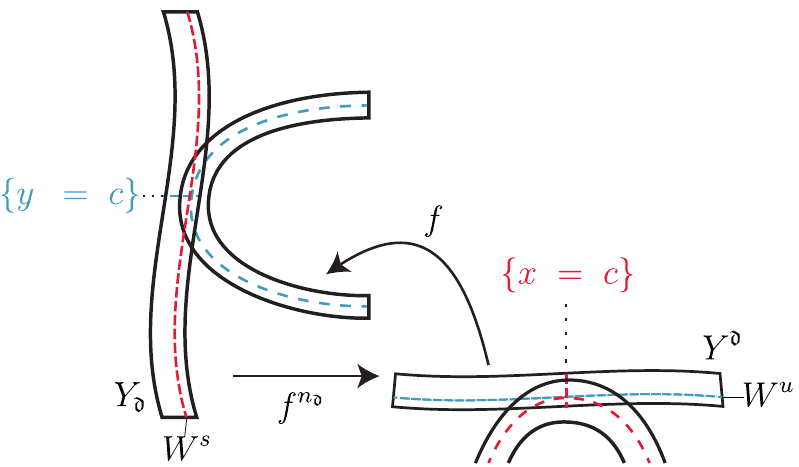}
 \caption{{ Preliminary setting for a Hénon-like renormalization}}
\label{notation_renormalization}\end{figure}
After the coordinate change $(x,y)\mapsto (c(p)+ x/ q_p,
c(p)+ y/q_p)$, we shall assume that for every $p$:
\[q_p=1\qand c(p)=0\; .\]
Here is the first renormalization Theorem, proved below after an example of applications of it.\label{sectionthmC}
\begin{theorem}\label{Charthenonlike}
Assume that $(i)-(ii)-(iii)$ hold true. Let $\sigma_p:= \partial_x A_p(0)$, $\lambda_p:=\partial_y B_p(0)$ and: 
\[\phi_p(X,Y)\mapsto (\sigma_p X, H_p(\sigma_p X) + \sigma_p \lambda _p Y)
\; . \]
Then $F_p:= \phi_p^{-1}\circ f^{n_\sd+1}_p\circ \phi_p$ is of the form:
$$F_{\bar p}(X,Y)=(X^2+Y+\bar a- \bar b^{ M} Y
,X) +(Z_{\bar p}(X,\bar b^{ M} Y), \Xi_{\bar p}(X,\bar b^{ M} Y))
\; ,$$
with $\bar a = \frac{\mu_p}{\sigma_p^2}$ , $\bar b:= 
\sqrt[M]{ \det \, D_{z_0} f^{n_\sd +1}_p}
$ and $\bar p:=(\bar a,\bar b)$. 

For every $\bar \delta>0$ and $R>0$, 
with $\Delta_p:= [-R,R]\times [- R/\bar b^M,R/\bar b^M]$, the family $(F_{\bar p})_{\bar p}$ of Hénon-like maps $F_{\bar p}$ restricted to $ \Delta_p\cap F_p^{-1}(\Delta_p)\cap \phi^{-1}_p(f_p^{-1}(Y_\sd))
$
is $\bar \delta$-$C^{d,r}$-like with multiplicity $M$, if $(\sigma_p)_p$ is sufficiently small in function of:
$$R, \quad \|\bar p\mapsto p\|_{C^d},\quad \|(r_p)_p\|_{C^{d,r}}\qand \mathcal B^{d,r}_1(\sd)\; . $$
\end{theorem}

\begin{rema} A key novelty in this proposition is that there is no loose of regularity in the Hénon-like family which is obtained (we start from a $C^{d,r}$-family and we obtain a $C^{d,r}$-family). 
We recall that by Proposition \ref{zetanull}, an ultimate coordinate change makes $\mathcal R\xi=0$ and so this proposition enables one to define infinitely renormalizable maps in the finitely and infinitely smooth settings. Let us also notice that this proposition is valid when $\bar b$ is not small, and even when $\bar b$ is large.

Interesting works on renormalizations of Hénon-like maps were done in \cite{CLM} and then in \cite{Ha11}. The techniques were designed for holomorphic Hénon-like maps with small determinant and for renormalizations of bounded periods. They actually used the renormalization chart equal to $\psi_p(X,Y)= (A_p(\sigma_p X, \sigma_p Y) , \sigma_p Y)$. Then $f^{n_\sd}\circ \psi_p(X,Y)= (\sigma_p X, B_p(\sigma_p X, \sigma_p Y))$ which is close to, but different to our chart. Then 
 their renormalizations leave invariant the space of Hénon-like maps of the form $f_p(x,y)=(x^2-by +\zeta(x,0)+ b y \zeta(x,y),x) $ which is larger than ours: they did not ask $\partial_y^j \zeta=O(b^j)$. I do not know how to prove the main theorems \ref{thm B} and \ref{thm A} for such general form. One of the problem is that for such a general form, there are no distortion bounds on $(B_p)_p$. 
\end{rema}
The following example is a study of the renormalized dynamics nearby the unfolding of a homoclinic tangency with two parameters. As pointed by the latter remark, the form obtained is more precise than the one of \cite{TLY, GST93, PT93}. Also, we are now able to give $C^r$-bound on the renormalization of $C^r$ dynamics, whereas these previous works gave no more than $C^{r-1}$-bound on this renormalization.
\begin{exam} \label{PTrevisited}
Let $0\le d\le r-3$ and consider a 2-parameters $C^{d,r}$-family $(f_p)_p$ of surface diffeomorphisms with a persistent saddle point $(S_p)_p$. Assume that 
$S_p$ displays a quadratic homoclinic tangency at the parameter $p=p_0$ and at the point $T_0$. Let $b:= \det D_S f_p$ and let $a$ be the distance 
between the local stable and unstable manifold nearby $T_0$. Assume that $p\mapsto (a,b)$ is a local diffeomorphism nearby $p=p_0$ (as it occurs typically). Then for any $\delta>0$, $R\ge 1$ and any $N\ge 1$ large, a renormalization of $(f_p)_p$ of period $N+1$ nearby the homoclinic tangency and the parameter $p=0$ is a $R$-wide, $\delta-C^{d,r}$-Hénon-like family of multiplicity $N+1$. 
\end{exam}
\begin{proof}Let us apply Theorem \ref{Charthenonlike}. Thanks to the local diffeomorphism $p\mapsto (a,b)$, we assume $p=(a,b)$ nearby $p_0$. For every $p$, we chose the coordinates so that $S_p=0$ , a local stable manifold of $S_p$ is $W^s_{loc}(S_p):=\{0\}\times (-1,1)$, a local unstable manifold of $S_p$ is $W^u_{loc}(S_p)=(-1,1)\times \{0\}$ and 
$T_0=(c_0,0)$ is sent by $f_{p_0}$ to $(0,c_0)$ with $c_0\in (-1,1)$. We can keep these properties and handle a deformation of the coordinate chart around $T_0$ so that at its neighborhood, the map $f_p(x,y)$ is of the form 
$(x,y)\mapsto (g_p(x,y),x)$. 
 Consequently hypothesis $(i)$ is satisfied. 

By the inclination lemma, when $N$ is large, 
the image by $f_p^{-N}$ of a small segment of the vertical line $\{x=c_0\}$ containing $T_0$ is a curve $\mathcal V_p$ which is $C^r$-close to $W^s_{loc}(S_p)$ and
the image by $f_p^{N}$ of a small segment of $\{y=c_0\}$ containing $f(T_0)$ is a curve $\mathcal H_p$ which is $C^r$-close to $ W^u_{loc}(S_p)$.
Let us prove that the families $(\mathcal V_p)_p$ and $(\mathcal H_p)_p$ are $C^{d,r}$-close to $(W^s_{loc}(S_p))_p$ and $(W^u_{loc}(S_p))_p$. Consider a tame box $Y_\sd\supset W^s_{loc}(S_p)$ sent by $f_p^N$ to a tame box $Y^\sd\supset W^u_{loc}(S_p)$. We observe that $\sd=(Y_\sd, N)$ is a hyperbolic puzzle piece, and it is equal to the $\star$-product of $N$-hyperbolic pieces of order 1 (for an adapted metric independent of $N$). Then by Proposition \ref{Birkhoff}, the affine-like representation $(A_p,B_p)_p$ of $\sd$ is $C^{d,r}$-bounded independently of $N$
and its derivatives are small. 
 As it satisfies $A_p(0,\cdot)=0$ and $B_p(\cdot,0)=0$, the family $(A_p, B_p)_p$ is $C^{d,r}$-small. Consequently $(\mathcal H_p)_p= (Graph B_p(\cdot ,c_0))_p$ is $C^{d,r}$-close to $(W^u_{loc}(S_p))_p$ and $(V_p)_p= (\, ^t Graph A_p(c_0, \cdot))_p$ is $C^{d,r}$-close to $(W^s_{loc}(S_p))_p$.

As the tangency is quadratic at $p=p_0$, for every $p$ close to $p_0$, there exists a unique point $z_0(p)$ in $\mathcal V_p$ close to $(0,c_0)$, whose image is at a minimal distance to $f_p(\mathcal V_p)$. As $d<r$, we observe that $p\mapsto z_0(p)$ is of class $C^d$ and $C^d$-bounded independently of $N$. Let $c(p)$ be the second coordinate of $z_0(p)$. Note also that $\sigma_p:= \partial_x A_p(c_p,0)$ is $C^d$-small, because all the derivatives of $(A_p)_p$ are small. 

Let $(H_p)_p$ and $(V_p)_p$ be the family of functions so that $\mathcal V_p:=\{(V_p(c(p)+y), y):y\in (-1,2)\}$ and $\mathcal H_p:=\{(x,H_p(c(p)+x):x\in (-1,2)\}$. Both are $C^{d,r}$-small, and so with $g_p$ the first coordinate of $f_p$, the family of maps $((x,y)\mapsto g_{p}(c(p)+ x, y+H_p(x))-V_p(x))_p$ 
satisfies condition $(T)$ with $p\mapsto \mu_p$ close to $p\mapsto a$,
a curve of parameters $\{p(b)=(a(b),b): b\}$ close to $\{(0,b): b\}$, 
$(\log |d_p|, \log |q_p|)_p$ $C^d$-bounded, and $(r_p)_p$ $C^{d,r}$-bounded when $N$ is large and $p$ close to $p_0$. 
Also the maps $(p\mapsto f_p^i (z_0(p)))_{0\le i\le N}$ are $C^d$-bounded and most of them are $C^d$- small by Lemma \ref{prebirkof2}. Consequently, $p\mapsto \sqrt[N+1]{|\det\, D_{z_0}f_{p}^{N+1}|}$ is $C^d$-close to $p\mapsto |b|$. Therefore condition $(ND)$ is also satisfied. Hence we can apply Theorem \ref{Charthenonlike} with $M=N+1$.
\end{proof}

\begin{proof}[Proof of Theorem \ref{Charthenonlike}]
By assumption, $\bar p\mapsto p$ is $C^d$-bounded, so it suffices to define $(Z_{ p}, \Xi_p)_{p}$ and show that this family is $C^{d,r}$-small when $(\sigma_p)_p$ is small.

To make the notations lighter, we omit to display the parameter $p$ in the indexes of the following computations. Let $(x_0,y_0)=f\circ \phi(X,Y)$ be sent by $f^{n_\sd}$ to 
 $(x_1, y_1)=\phi(X_1, Y_1)$. Observe that:
\begin{equation}\label{preR0} x_0=g(\sigma X, H(\sigma X) + \sigma \lambda Y) \; ,\; y_0=\sigma X\; ,\; x_1=\sigma X_1\qand y_1= H(\sigma X_1) + \sigma \lambda Y_1\; .\end{equation}
By (\ref{Defrp}):\begin{equation}\label{preprehenon}
g\circ \phi(X, Y)= V(\sigma X) + \mu+ \sigma^2 X^2- d\cdot 
\sigma \lambda Y+r(
\sigma X, d\cdot \sigma \lambda Y)\; .\end{equation} 

As $\phi(X,Y)$ is sent by $f^{n_\sd+1}$ to $\phi(X_1,Y_1)$, the point $(x_0,y_0)$ is sent by $f^{n_\sd}$ to $(x_1,y_1)$. Hence with $\breve A(x_1, y_0):= A(x_1,y_0)-V(y_0)$ and $\breve B(x_1, y_0):= B(x_1, y_0)-H(x_1)$ it holds:
\begin{equation}\label{R0}
x_0= \breve A(x_1, y_0) + V(y_0)\qand y_1= \breve B(x_1,y_0)+H(x_1)\; .
\end{equation}
By (\ref{preR0},\ref{preprehenon},\ref{R0}), it comes:
\begin{equation}\label{pourmainthm}
\left\{\begin{array}{c}
\breve A(\sigma X_1, \sigma X) = \breve A(x_1, y_0) = \mu+ \sigma^2 X^2- d \cdot \sigma \lambda Y+r(
\sigma X, d\cdot \sigma \lambda Y)\\
 \breve B(\sigma X_1,\sigma X)=\breve B(x_1,y_0)= \sigma \lambda Y_1
\end{array}\right. \end{equation}

We recall that by Fact \ref{determinant}, $\det D_{z_0}f^{n_\sd}= \lambda/\sigma$ and so 
$d \cdot \lambda/\sigma = \det D_{z_0}f^{n_\sd+1}= \bar b^M$. 
We put $\Pi_A(X_1,X):= \breve A(\sigma X_1, \sigma X)/ \sigma^2$ and $\Pi_B(X_1,X):= \breve B(\sigma X_1, \sigma X)/ (\sigma\lambda)$. Then the system is:
\[\left\{\begin{array}{c}
\Pi_A(X_1,X) = \bar a + X^2- \bar b^M Y+\sigma^{-2} r(
\sigma X, \sigma^2 \bar b^M Y)\\
 Y_1= \Pi_B (X_1,X)\end{array}\right. \]

We assume that $(X,Y)$ and $(X_1,Y_1)$ are in $\Delta$ and so $(X,X_1)$ are in $[-R,R]^2$. By Lemma \ref{distorsionutil}, the family of functions $(\Pi_A,\Pi_B)_p$ is $C^{d,r}$-close to the family constantly equal to the identity. 
Let $\rho_p(x,y):=\sigma^{-2}_p r_p(\sigma_p x,\sigma^2_p y)$. We notice that $(\partial_x ^i \partial_y ^j\partial_p^k \rho_p)_p$ is $C^0$-dominated by $\sigma$ for every $i+j+k\le r$, $k\le d$ and 
$i+2j\ge 3$.
 As $\rho_p(0)=\partial_y \rho_p (0)=\partial_x\rho_p (0)=\partial^2_x\rho_p (0)=0$ for every $p$, the family $(\rho_p|[-R,R]^2)_p$ is $C^{d,r}$-small when $\sigma$ is small.

By the implicit function theorem, $X_1$ is a bounded function of $(X,Y)$ and so there is a $C^{r}$-$\bar \delta$-small family of functions $(Z_{p})_{p}$ satisfying $Z_p(X,\bar b^M Y) = X_1-\Pi_A(X_1,X)+\rho_p(X,\bar b^M Y)$. Thus:
 \[X_1 =: X^2 + \bar a- \bar b^M Y + Z_p(X, \bar b^M Y)\; . \]
 By injecting the latter expression of $X_1$ into $Y_1= \Pi_B (X_1,X)$, it comes that $Y_1= X+\Xi_p(X,\bar b^M Y)$, for a family $(\Xi_p)_p$ which is $C^{d,r}$-small. 

\end{proof}

\subsection{Multi-Renormalization}
Let $(f_p)_p$ be a $N$-parameters $C^{d,r}$-family of surface diffeomorphisms, with $d\le r-3$. 
Let $(\sd_i)_{i\in \Z/N\Z}$ be pieces with affine-like representations $(A_{i\, p},B_{i\, p})_p$, and assume that each of their domains of definition contains $0$. Put $n_i:= n_{\sd_i}$. Let us define:
$$V_{i\, p}(y):= A_{i\, p}(0,y)\qand 
W^s_{i\, p}:= f_{p}^{-n_i} (\{x=0\}\cap Y^{\sd_i})= \{(A_{i\, p}(0, y),y):y\}\; .
$$
$$H_{i\, p}(x):= H_{i\, p}(x,0) \qand W^u_{i\, p}:= f_{p}^{n_i} (\{y=0\}\cap Y_{\sd_i}) = \{(x, B_{i\, p}(x, 0)):x\}
\; .$$ 
\begin{figure}[h!]
 \centering
 \includegraphics[width=12cm]{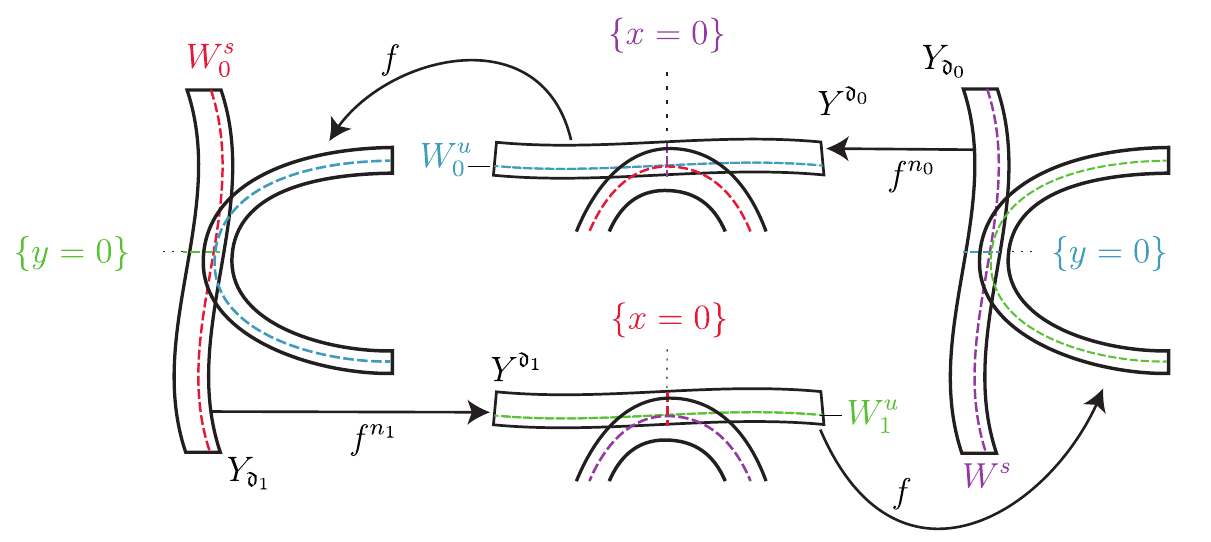}
 \caption{{ Preliminary setting for a multi-renormalization with $N=2$.}}
\label{notation_renormalization}\end{figure}
Note that $ z_i(p):=(0,H_{i\, p}(0))$ belongs to $W^u _{i\, p}\subset Y^{\sd_i}$.
Let us assume the following properties:
\begin{enumerate}[$(i)$]
\item There are neighborhoods $U_i$ of $ z_i(p)$ so that the restriction $f_p|\cup_i U_i$ is of the form $f_p:(x,y)\in \cup_i U_i\mapsto (g_p(x,y),x)$, for every $p$.
\item There exists a parameter $p=p_0$ so that the curve $W^u_{i\, p}$ is sent by $f_{p }$ to a curve tangent to $W^s_{i+1\, p}$ at the point $f_{p}(z_i(p))= (V_{i+1\, p}(0),0)$, 
and the curve $W'^u_p$ is sent by $f_p$ tangent to $W^s_p$ at the point $f_{p}(z'_i(p))= (V_{i+1\, p}(0),0)$. 
\item The tangencies are quadratic. In other words, it holds:
\begin{equation}\tag{$\mathbb T$}\label{Defrp12}
g_{p}(x, y+H_{i\, p}(x))- V_{i+1\, p}(x) = q_{i\, p} \cdot x^2+\mu_{i\, p}-d_{i\, p} \cdot y+r_{i\, p}(x, d_{i\, p}\cdot y) \; ,\end{equation}
 with $q_{i\, p}\not = 0$, 
$d_{i\, p}:= \det D_{z_i(p)} f_p$ and the family $(r_{i\, p})_p$ of class $C^{d,r}$ which satisfies for all $p$:
$$0= r_{i\, p}(0)= \partial_y r_{i\, p}(0)=
 \partial_x r_{i\, p}(0)=\partial_x^2 r_{i\, p}(0)\; .$$
\end{enumerate}

After a local coordinate change we shall assume that $q_{i\, p}=1$ for all $p$ and $i$. For every $i$, put:
\[\sigma_{i\, p}:= \partial_x A_{i\, p}(0)\; ,\quad \lambda_{i\, p}:= \partial_y B_{i\, p}(0)\; ,\quad \gamma_{i\, p}:= sign(\sigma_i)\cdot 
 \sigma_{i+1\, p}^{1/2}\cdot \sigma_{i+2\, p}^{1/2^2}\cdots \sigma_{i+j\, p}^{1/2^{j}}\cdots \; .\]
 We observe that $(\gamma_{i\, p})_{i\, p}$ are well defined, of class $C^d$, and satisfies $\gamma_{i\, p}^2= \gamma_{i+1\, p}\cdot \sigma_{i+1\, p}$. 
 \label{sectionthmD} 
\begin{theorem}[Third main result]\label{Charthenonlike2}
Assume that $(i - ii- iii)$ hold true.
 Put:
\[\phi_{i\, p}(X,Y)= (\gamma_{i\, p} X,
H_{i\, p}(\gamma_p X)+\gamma_{i+1\, p}\lambda_{i\, p} Y)\; .\]
The map $\mathcal R_i f_p= (\phi_{i+1\, p})^{-1}\circ f_p^{n_{{i+1}}+1}\circ \phi_{i\, p}$ is Hénon-like and of the form:
\[ \mathcal R_if_p(X,Y)= (X^2+\bar a_i-\bar b_i Y, X)+ 
(Z_{i\, p}, \Xi_{i\, p})(X,\bar b_iY)
 \; ,\]
 \[\text{with }\quad (\bar a_i,\bar b_i):= \left(\frac{\mu_p}{\gamma_{i\, p}^2}, \frac{\gamma_{i-1\, p}\cdot \sigma_{i\, p}}{\gamma_{i+1\, p}\cdot \sigma_{i+1\, p}} \cdot \det D_{f_p^{-n_i}(z_{i}(p))} f_p^{n_{i}+1}\right)\] 
 Moreover, for every $R>0$, with $\Delta_{i\, p}:= [-R,R]\times [-R/b_i,R/b_i]$, 
 the map $\mathcal R_i f_p$ restricted to $\Delta_{i\, p}\cap \phi^{-1}_{i\, p}(f_p^{-1}(Y_{\sd_{i+1}}))\cap \mathcal R_if_p^{-1}(\Delta_{i+1\, p})$, the families $(Z_p, \Xi_p)_p$ are $C^{d,r}$-small 
 if $(\gamma_{i\, p})_{i\, p}$ is $C^d$-small w.r.t.:
$$R,\quad \|(r_{i\, p})_{i\, p}\|_{C^{d,r}}\qand (\mathcal B^{d,r}_1(\sd_i))_i \; .$$

\end{theorem}
A related result to Theorem \ref{Charthenonlike2} in the context chain of flat heteroclinic tangencies is given in \cite{Tu15}, see also \cite[\textsection 6]{BT17}. In the conservative context and for $N=2$, \cite[Formula (2.3)]{DGLS13} shows that a composition of two Hénon-like maps appears as a renormalized dynamics nearby a chain of two heteroclinic tangencies. Given a $C^r$-dynamics, they gave $C^{r-1}$-bound on the renormalization, Theorem \ref{Charthenonlike2} and Proposition \ref{Birkhoff} improve these bounds to the $C^r$-topology.

Theorem \ref{Charthenonlike2} implies Theorem \ref{Charthenonlike}.
Before proving it, let us first state an immediate consequence of it with Proposition \ref{compo2Swallow}; this will be useful for the proof of Theorem \ref{thm B}. 
\begin{coro}\label{preThmB}
Under the hypotheses of Theorem \ref{Charthenonlike2} with $N=2$, for every $\delta>0$, if moreover the two following conditions hold true:
\begin{itemize}
\item the map $p\mapsto (a_0,a_1)$ is a diffeomorphism onto $[-R,R]^2$ with inverse $C^d$-close to a constant,
\item the families $(b_{0\, p})_p$ and $(b_{1\, p})_p$ are $C^d$-small,
\end{itemize}
 then $(\phi_{0\, p}^{-1}\circ f^{n_0 +n_1+2}\circ \phi_{0\, p})_{(a_0, a_1)\in [-R,R]^2}$ is conjugated to a $\delta$-$C^{d,r}$-swallow-like family.
\end{coro}
The Milnor swallows appear in many parameter space pictures of surface diffeomorphisms. To explain this, the reader might\footnote{Hint: proceed as in Example \ref{PTrevisited} with some tools of the forthcoming proof of Theorem \ref{thm B}.} show a positive solution to the following problem: 
\begin{prob} Let $f_0$ be a $C^{r}$-surface diffeomorphism which is area contracting and which displays a hyperbolic set with a robust homoclinic tangencies, with $r\ge 3$.
Show that for every $d\le r-3$, for an open and dense set of $2$-parameters $C^{d,r}$-families $(f_p)_p$ containing $f_0$, the conclusions of Theorem \ref{thm B} holds true for $(f_p)_p$. \end{prob}
\begin{proof}[Proof of Theorem \ref{Charthenonlike2}]
The proof is similar to the one of Theorem \ref{Charthenonlike}, \textsection \ref{sectionthmC}.
 We fix the parameter $p$ and remove its writing as index in the computations below. Let $(x_0,y_0)= f\circ \phi_i(X,Y)$ be sent by $f^{n_{{i+1}}}$ to $(x_1,y_1)= \phi_{i+1}(X_1,Y_1)$. Observe that: 
\begin{equation}\label{gxswal0}
x_0= g(\gamma_i X, H_i(\gamma_i X)+\gamma_{i-1} \lambda_i Y)\; ,\; y_0=\gamma_i X\; ,\; x_1 =\gamma_{i+1}X_1\; \qand y_1= H_{i+1}(x_1)+\gamma_i \lambda_{i+1} Y_1
\; .\end{equation}
 By (\ref{Defrp12}), we have:
\begin{equation}\label{gxswal}
x_0 = V_{i+1}(\gamma_i X) + \gamma_i^2\cdot X^2+\mu_i-d_i \cdot \gamma_{i-1}\lambda_i \cdot Y+r(\gamma_i X, d_i\cdot \gamma_{i-1}\lambda_i Y)\end{equation}
As $(x_0,y_0)$ is sent by $f^{n_{{i+1}}}$ to $(x_1,y_1)$, with $\breve A_{i+1}(x_1,y_0)=A_{i+1}(x_1,y_0)-V_{i+1}(y_0)$ and $\breve B_{i+1}(x_1,y_0)=B_{i+1}(x_1,y_0)-H_{i+1}(x_1)$, it holds:
\begin{equation}\label{NFSWAL}
x_0= \breve A_{i+1}(x_1, y_0) + V_{i+1}(y_0) \text{ and }y_1= \breve B_{i+1}(x_1,y_0)+H_{i+1}(x_1)\; .
\end{equation}
Then (\ref{gxswal0}-\ref{gxswal}-\ref{NFSWAL}) and $\det D_{f_p^{-n_i}(z_{i}(p))} f_p^{n_{i}+1}= d_i\lambda_i/\sigma_i$ imply:
\[\left\{
\begin{split}
\gamma_i^2\cdot X^2+\mu_i-\det D_{f_p^{-n_i}(z_{i}(p))} f_p^{n_{i}+1}\sigma_i \gamma_{i-1} \cdot Y+r(\gamma_i X, d_i\gamma_{i-1}\lambda_i Y) = \breve A_{i+1}(\gamma_{i+1} X_1,\gamma_i X) \\
 \gamma_i \lambda_{i+1} Y_1= \breve B_{i+1}(\gamma_{i+1} X_1,\gamma_i X) 
\end{split}\right.\]
We divide the upper equation by $\gamma_i^2 = \gamma_{i+1}\sigma_{i+1}$ and the lower equation by $\gamma_i \lambda_{i+1}$ to obtain:
\begin{equation}\label{finalswallow}\left\{\begin{array}{rl}(\sigma_{i+1}\gamma_{i+1})^{-1}\breve A_{i+1}(\gamma_{i+1} X_1 , \gamma_i X)&= a_i+ X^2-b_i\cdot Y+ (\gamma_i)^{-2} r(\gamma_i X,b_i \gamma_i^2 Y)\\
Y_1&= (\gamma_i \lambda_{i+1})^{-1}\breve B_{i+1}(\gamma_{i+1} X_1,\gamma_i X) 
\end{array}\right.
\end{equation}

By Lemma \ref{distorsionutil}, the maps $(\sigma_{i+1} \gamma_{i+1})^{-1}\breve A_{i+1}(\gamma_{i+1} X_1,\gamma_i X)$ and $ (\gamma_i \lambda_{i+1})^{-1} \breve B_{i+1}(\gamma_{i+1} X_1,\gamma_i X)$ form a family $C^{d,r}$-close to the identity. Thus, by the implicit function theorem, $(X_1,Y_1)=\mathcal R_i f(X,Y)$. Then we conclude as in the proof of Theorem \ref{Charthenonlike} to show the existence of $C^{d,r}$-small families $(Z_{i\, p}, \Xi_{i\, p})_{p}$ so that 
$\mathcal R_i f_p(X,Y)= (X^2+a_i-b_iY+Z_{i\, p}(X,b_i Y), X+ \Xi_{i\, p}(X,b_iY))$. 
\end{proof}
\begin{rema}
Note that instead of assuming $i\in \Z/N\Z$ in the statement of Theorem \ref{Charthenonlike2}, we can assume $i\in \N$, as far as $(\gamma_{i\, p})_{i\, p}$ is a well defined sequence of $C^d$-functions (and $C^d$-small by assumption of the Theorem). This may allow one to construct pathological wandering domains for some surface dynamics. 
\end{rema}
\begin{rema}
We notice that if we apply the theorem for an area conservative map with $N=2$, $|\gamma_0|\not= |\gamma_1|$, then the multi-renormalization is the composition of an area expanding Hénon-like map with an area contracting Hénon-like map.
\end{rema}

\section{Proof of Theorems \ref{thm B} and \ref{thm A}}
Both theorems will be proved for Hénon-like maps with parameter $a$ 
close to $a_2$. We recall that $a_2$ was defined in \textsection\ref{3.1} as satisfying $a_2=-\alpha_2(a_2)$, and so $Q_{a_2}^4(0)$ equal to the fixed point $\alpha$.
\subsection{Preliminary bounds for Hénon-like maps with parameter $a$ nearby $a_2$}\label{near a2} We recall:
\begin{prop}[Prop. 3.1 \cite{Y97}]\label{transpara1}
For every $a\le a_1$, the preimage $\alpha_2(a)$ of $-\alpha(a)$ by $Q_a^2$ satisfies $\partial_a \alpha_2(a)\in [-1/2, -1/3]\; .$
\end{prop}

\subsubsection{The one-dimensional picture}\label{271}
In \textsection \ref{3.1}, for $a< a_1$, we defined the pieces $\{\ss_+,\ss_-,\sw_+, \sw_-,\sw_=\}$ of $Q_a$, for every $a<a_1$. We recall that $a_2<a_1$. Let us define the following sequence of pieces:
\[\sc_0:=\sw_=,\quad \sc_1:= \sw_=\star \ss_+,\quad \sc_{m+1}:=\sc_m\star \ss_-= \sw_= \star \ss_+ \star \ss_-^{\star (m-1)} \quad \forall m\ge 1\; .\]
We notice that $\sc_m$ is a puzzle piece and that the left hand side endpoint of $\R_{\sc_m}$ is equal to $-\alpha_2$ for every $m\ge 0$. 

Let $\alpha_3$ be the preimage of $\alpha_2$ by $Q_a|\R^+$, and let 
$\sw_\equiv:= ([-\alpha_3, -\alpha_2],3)$. 
It is a puzzle piece. Also ${\bar \sc_m}:= \sw_\equiv \star \ss_-^{\star (m)}$ is a puzzle piece and the right hand side endpoint of $\R_{\bar \sc_m}$ is equal to $-\alpha_2$ for every $m\ge 0$. Furthermore:
 \[\{-\alpha_2\} = \bigcap_{m\ge 0} \R_{\sc_m}\cup \R_{\bar \sc_m} \qand \forall i\ge 0,\quad d(\alpha_2, \R_{\sc_m}\setminus \R_{\sc_{m+1}})=\leb\, \R_{\sc_{m+1}}
 >0 \]
In this one dimensional context, for every $m\ge 0$, the symbols $\bm (\sc_m-\sc_{m+1})$ and $\bp (\sc_m-\sc_{m+1})$ denote the pairs
 $(\R_{\bm (\sc_m-\sc_{m+1})},n_{\bm (\sc_m-\sc_{m+1})})$ and 
$(\R_{\bp (\sc_m-\sc_{m+1})},n_{\bp (\sc_m-\sc_{m+1})})$ defined by:
\[\left\{\begin{array}{c}
\R_{\bm (\sc_m-\sc_{m+1})}=(Q_a|\R^-)^{-1} cl (\R_{\sc_m}\setminus \R_{\sc_{m+1}})\qand
\R_{\bp (\sc_m-\sc_{m+1})}= (Q_a|\R^+)^{-1}cl (\R_{\sc_m}\setminus \R_{\sc_{m+1}})\; ,\\
n_{\bm (\sc_m-\sc_{m+1})}=1+n_{\sc_m}= 3+2m=n_{\bp (\sc_m-\sc_{m+1})}\; .\end{array}\right.
\]
For every $j\ge 0$, let $I_j$ be the interval of parameters 
$a<a_1$ which belong to $\R_{\sc_{j+2}}\cup \R_{\bar \sc_{j+2}}$. We observe that every pair in the following set is a piece:
$$\sA_j:= \{\ss_-,\ss_+\}\cup \bigcup_{m\le j} \{\bm(\sc_m-\sc_{m+1}), 
\bp(\sc_m-\sc_{m+1}) \}\; .$$
As $a\in \R_{\sc_{j+2}}\cup \R_{\bar \sc_{j+2}} $, 
for every $\sd\in\sA_j$, the distance from the segment $\R_\sd$ to $0$ is at least:
 $$
\eta_j:= \frac12 \min_{a\in I_j} \sqrt{\leb(\R_{\sc_{j+1}}\setminus \R_{\sc_{j+2}})}>0\; .$$
 We observe that:
\[\bigcup_{\sw\in \{\sw_-,\sw_+,\sw_=\}} \R_\sw \cup \bigcup_{\sd\in \sA_j} \R_\sd
\subset \R\setminus (-2\cdot \eta_j,2\cdot \eta_j) \; .\]
\begin{lemm} \label{hypQ}
There exists $\kappa<1$, so that for every $j\ge 1$ large, there exists $C>0$ satisfying for every $a\in I_j$ and $n\ge 0$: 
 \begin{equation}\label{hypQQ}
 \| D_zQ^n_a \|\ge C\cdot \kappa^{-n}, \quad \forall z\in \cap_{j=0}^n Q_a^{-j}(\R\setminus (-2\cdot \eta_j, 2\cdot \eta_j))\; .
 \end{equation}
Moreover, when $j$ is large, for any $\sq\in \{\bm (\sc_j-\sc_{j+1}), \bp (\sc_j-\sc_{j+1})\}$, $\log \eta_j\| D_zQ^{n_\sq}_a| \R_{\sq} \|$ is bounded. Also it holds $\log (\eta_j \min_{a\in I_J, x\in \R_\sq} |x|)=O(1)$ and $\log (\eta_j \max_{a\in I_J, x\in \R_\sq} |x|) =O(1)$ when $j\to \infty$. 
\end{lemm} 
\begin{proof}
Let us recall that by a famous result of Ma\~né, for every 
$j_0\ge 0$, there exist $\kappa_0<1$ and $C_0>0$ so that 
\[ \| D_zQ^n_{a_2} \|\ge C_0\cdot \kappa_0^{-n}, \quad x\in \cup_{j=0}^n Q_{a_2}^{-j}(\R\setminus (-\sqrt{\leb{\R_{\sc_{j_0}}}},\sqrt{\leb{\R_{\sc_{j_0}}}}))\]
Hence to prove the first statement of the Lemma, it suffices to show the existence of $\kappa\in(\kappa_0,1)$ so that for every 
$\sq\in \sA_j\setminus \sA_{j_0}$, for every $x\in \R_{\sq}$, it holds $
|DQ_a^{n_{\sq}}(x)| >\kappa^{-n_\sq}\; .$ 

We are going to show more than this. We shall prove that $\log |DQ_a^{n_{\sq}}(x)| \sim 
({n_\sq /2}) \log ({-2\alpha})\sim \log \max_{\R_\sq} |x|$ ; this will also implies the second statement of the lemma.

In order to do so, we apply Sternberg Linearization theorem at the fixed point $\alpha$ of $Q_a$. It implies that the length of $ \R_{c_{k+1}}$ is of the order of $(-2\alpha)^{-n_{\sc_k}}$ ; also $ |D Q_a^{n_{\sc_k}}|_{|\R_{c_k}}$ is of the order of $(-2\alpha)^{-n_{\sc_k}}$. Thus, when $j_0$ is large, for every $\sq\in \{\bm (\sc_k-\sc_{k+1}), \bp (\sc_k-\sc_{k+1})\}$ with ${j_0\le k\le j}$, the distance from 0 to any point of $\R_{\sq}$ is of the order of $(-2\alpha)^{-n_{\sc_k}/2}$, and so $ \|D Q_a^{n_{\sq}}|\R_{\sq}\|$ is of the order of $(-2\alpha)^{-n_{\sc_k}/2}$.
\end{proof}

\subsubsection{The bounds for pieces of Hénon-like maps}\label{henonpiece}
We continue with the interval $I_0$ defined in the latter section. We recall that 
if $\delta$ and $\hat b=\max_J |b|$ are sufficiently small, for every $\delta$-$C^{d,r}$-Hénon-like family $(f_{p})_{I_0\times J }$,
every symbol $\sd$ in the set $\{\ss_-,\ss_+,\sw_-,\sw_+,\sw_=\}$
defines a puzzle piece by Proposition \ref{boxtame} for any cone constants $c_h=1/\eta$, $c_v=\eta/2$ and $\eta\in (0,\frac12 \min_{I_0}\sqrt{a_1-a})$. 

The following is a direct consequence of Lemma \ref{distorsionordrefini}. 
\begin{fact}\label{K0} 
There exist $K_0>0$, $\delta_0>0$ and $b_0>0$ such that for every $m\ge 1$, for every $\delta_0$-$C^{d,r}$-Hénon-like family $(f_{p})_{I_0\times J }$ of multiplicity $m$ and with $J\subset [-b_0,b_0]$, the affine-like representation of any $\sd\in \{\ss_-,\ss_+,\sw_-,\sw_+,\sw_=\}$ satisfies that 
$\mathcal B^{d,r}_0(\sd)$, $\mathcal B^{d,r}_1(\sd)$ and $ {\mathcal B}^{d}_m(\sd)$ are at most $K_0$. 
\end{fact}
A similar bound exists on the $\star$-product of the above pieces:
\begin{prop}\label{hauteurbandes}
There exists $\hat K_0>0$ such that under the hypotheses of Fact \ref{K0}, for every $c\in \R_\boxdot$, it holds:
\begin{itemize}
\item[$(1)$] For every $\sd \in \{\ss_-,\ss_+,\sw_-,\sw_+,\sw_=\}$, the vertical line $\{x=c\}$ intersects $Y^\sd$ at a segment of length in $[\hat K_0^{-1}\cdot |b|^{m \cdot n_\sd}, \hat K_0\cdot |b|^{m \cdot n_\sd}]$. 
 The length of the bounded component of $\{x=c\}\cap (Y^{\sw_=}\cup Y^{\sw_-})$ is in $(\hat K_0^{-1},\hat K_0)$.\end{itemize}

\noindent And for all $j\ge 0$, $(\sd_i)_{i}\in \{\ss_- ,\ss_+\}^j$ and $\sw\in \{\sw_-,\sw_+,\sw_=\}$ and $\sd\in \{ \sd_1\star \cdots \star \sd_m, \sw \star \sd_1\star \cdots \star \sd_m\} $: 

\begin{itemize}
\item[$(2)$] The piece $\sd$ satisfies that $\mathcal B^{d,r}_0(\sd)$, $\mathcal B^{d,r}_1(\sd)$ and $ {\mathcal B}^{d}_m(\sd)$ are at most $\hat K_0$.
 \end{itemize} 
\end{prop}
\begin{proof}The first bound of $(1)$ is a direct consequence of the $\mathcal B^{d,r}_1$-bounds given by Fact \ref{K0}. 

The horizontal distance between the vertical strips $Y_{\sw_-}$ and $Y_{\sw_+}$ is close to $2\alpha_1>1$. Then we that $f_p$ maps a vertical segment to a horizontal segment. Thus the vertical distance between $Y^{\sw_-}$ and $Y^{\sw_+}\supset Y^{\sw_=}$ is in $(\hat K_0^{-1}, \hat K_0)$ for some uniform $\hat K_0$. 
This proves the second part of $(1)$.

Assertion $(2)$ is a direct consequence of Propositions \ref{Birkhoff} and Fact \ref{K0}. 
\end{proof}

\medskip 

Let $j\ge 1$ be large, and let $b_j, \delta_j>0$ be small in function of $j$. 
Let $f=f_{a\, b}$ be a $\delta_j$-$C^r$-Hénon-like map so that $a\in I_j$ and $|b|\le b_j$. We now fix $\eta$ (defining the cone constants) as equal to $\eta_j$ (which was defined for Lemma \ref{hypQ}).

\medskip 

By Propositions \ref{boxtame} and \ref{defstarhenon}, 
a puzzle piece $(Y_{\sc_k}, n_{\sc_k})$ is canonically associated to the symbol $\sc_k:= \sw_=\star \ss_+\star \ss_-^{\star (k-1)}$, for every $k\ge 1$. In particular $cl(Y_{\sc_m}\setminus Y_{\sc_{m+1}})$ is a vertical strip. Also for every $k< j$, and for $b_j$ and $\delta_j$ small enough, its preimage by $f$ intersects $\R\times [-3,3]$ at two vertical boxes $Y_{\bm (\sc_k-\sc_{k+1})}$ and $Y_{\bp (\sc_k-\sc_{k+1})}$. 
These boxes are close to respectively $\R_{\bm (\sc_k-\sc_{k+1})}\times [-3,3]$ and $\R_{\bp (\sc_k-\sc_{k+1})}\times [-3,3]$.

\begin{fact}\label{coneta2} 
Let $\sd$ be of the form $\sd=\bm(\sc_k-\sc_{k+1})$ or $\sd=\bp(\sc_k-\sc_{k+1}) $ for $k< j$.
Then $(Y_\sd,n_\sd)$ is a piece also denoted by $\sd$.
Furthermore, $Y_\sd$ is included in $Y_\boxdot$ and the set $\partial^s Y^\sd$ displays one component in $Y_{\ss_-}$ and one component in $Y_{\ss_+}$. 
\end{fact}
\begin{proof}
We already saw that $Y_\sd$ is a vertical strip. Also $\partial^s Y_\sd$ is sent by $f^{n_\sd}$ into 
the stable border of the vertical strip $cl(Y_\se\setminus Y_{\ss_+})$ if $k=1$ or 
$cl(Y_\se\setminus Y_{\ss_-})$ if $k>1$. By remark \ref{remacone}, $\sd$ is a piece if the cone condition is satisfied. 
To prove this, we note that $Y_\sd$ is disjoint from $[\eta_j,\eta_j]\times \R$, and the same holds for every $f^{i}(Y_\sd)$ for every $i< n_\sd$. Hence by Lemma \ref{coneta}, $f^{n_\sd}|Y_\sd$ satisfies the cone condition. 
\end{proof}

By Proposition \ref{defstarhenon}, for any $k'\ge k\ge 1$, for any $(\sd_i)_{1\le i\le k}\in \{\ss_- ,\ss_+\}^k$ and $(\sd_i)_{k+1\le i\le k'}\in (\sA_j\setminus \{\ss_- ,\ss_+\})^{k'-k}$ the piece
$\sd:= \sd_1\star \cdots \star \sd_k\star \sd_{k+1}\star \cdots \star \sd_{k'}$ is well defined.

\begin{prop}\label{propVp}
For every $j\ge 1$, there exist $\hat K_j$, $b_j$ and $\delta_j>0$ , so that for every $C^{d,r}$-$\delta$-Hénon-like family $(f_{p})_{I_j\times J}$ of multiplicity $m\ge 1$ and with $J\subset [-b_j,b_j]$, the following properties hold true. 
\begin{enumerate}[$(1)$]
\item for every $\sd\in \sA_j$, the bound $\mathcal B^r(\sd)$ is at most $\hat K_j$. Also when $j$ is large, uniformly among $z\in Y_{\bm(\sc_j-\sc_{j+1})} \cup Y_{\bp(\sc_j-\sc_{j+1})}$ and $p\in I_j \times J$, it holds that 
$ \log (\eta_j \| \partial_x f^{3+ 2j}_p(z)\|)$ is bounded and the distance from $z$ to $\{0\}\times \R$ is of the order of $\eta_j$,

Let $k'\ge k\ge 1$, $\sw\in \{\sw_-,\sw_+,\sw_=\}$, 
 $(\sd_i)_{1\le i\le k}\in \{\ss_- ,\ss_+\}^k$ and $(\sd_i)_{k+1\le i\le k'}\in (\sA_j\setminus \{\ss_- ,\ss_+\})^{k'-k}$, and put $\sd':= \sd_1\star \cdots \star \sd_k\star \sd_{k+1}\star \cdots \star \sd_{k'}$. 
\item Any piece $\sd\in \{\sd', \sw\star d'\}$ satisfies that 
 $\mathcal B^{d,r}_0(\sd)$, $\mathcal B^{d,r}_1(\sd)$ and $\bar {\mathcal B}^{d}_m(\sd)$ are $\le \hat K_j$.
\end{enumerate}
\end{prop}
\begin{proof}
As we saw in the previous proof, for every $\sd\in \sA_j$ and $k\le n_\sd$, $f^{k}(Y_\sd)$ is included in $\mathcal H_{\eta_j}=Y_D\setminus (-\eta_j,\eta_j)\times \R$.
As $\delta_j\le \eta_j/2$ and $b_j$ is small enough, Lemma \ref{hypQ} and Proposition \ref{Birkhoff} imply $(1)$ and $(2)$. 
\end{proof}

Let $(f_{a\, b\, m})_{I_0\times J}$ be a $\delta$-$C^{d,r}$-Hénon-like family with $\delta$ small and put $\hat b=\max_J |b|$ small. Let $-\alpha_2(f_{a\, b\, m})$ be the unique intersection point of the line $\{y=0\}$ with the left hand side component $\mathcal L_{-\alpha_2}$ of $\partial^s Y_{\sw_=}$. By hyperbolic continuation, the map $(a,b)\mapsto -\alpha_2(f_{a\, b\, m})$ is of class $C^d$. As for every $b$, the one-parameter family $(f_{a\, b\, m})_{a\in I}$ is $O(\hat b^m+\delta)$-$C^{d,r}$-close to $((x,y)\mapsto (x^2+a,x))_a$, by Proposition \ref{transpara1} it holds:
\begin{fact}\label{transpara3}
The map $a\in I_0\mapsto -\alpha_2(f_{a\, b\, m})$ is $O(\hat b^m+\delta)$-$C^d$-close to $a\in I_0\mapsto -\alpha_2(a)$ and its derivative is close to be in $[1/3,1/2]$.
\end{fact}

Here is a consequence of Propositions \ref{continuityinprod} and \ref{propVp}.
\begin{coro}\label{transpara2}
Let $d\le r-2$. When $\delta$ is small and $\hat b= \max_J |b|$ is small, for every $k$ large, 
with $(A_{a\, b}, B_{a\, b})_{I\times J}$ the affine-like representation of $\sc_k$, it holds
$(A_{a\, b})_{I\times J}$ is $O(|\partial_x A(0)|+\hat b^m)-C^{d,r}$-close to the family of constant mappings $((x,y)\mapsto -\alpha_2(f_{a\, b\, m}))_{I\times J}$.

For every $j\ge 1$, if furthermore $I\subset I_j$, $\delta\le \delta_j$, and $k$ is large in function of $j$, then for every $\sd\in \sA_j$, 
 the implicit representation $(A'_{a\, b},B'_{a\, b})_{I\times J}$ of $\sc_k\star \sd$ satisfies that $(A'_{a\, b})_{I\times J}$ is $O(\hat K_j\cdot|\partial_x A(0)|+\hat b^m)- C^{d,r}$-close to the family of constant mappings $((x,y)\mapsto -\alpha_2(f_{a\, b\, m}))_{I\times J}$.
\end{coro}
\begin{proof}
The second statement follows from the first by Propositions \ref{continuityinprod} and \ref{propVp}.(1). 

To show the first statement, we first recall that by Lemma \ref{distorsionutil0}, 
the families $(\partial_x A_{a\, b})_{I\times J}$ is $C^{d, r-1}$-$O(|\partial_x A(0)|)$-$C^{d,r-2}$-small, and by Proposition \ref{bsmall}, $(\partial_y A_{a\, b})_{I\times J}$ is $C^{d, r-1}$-$O(\hat b^m)$-$C^{d,r-1}$-small. Hence $(A_{a\, b})_{I\times J}$ is $O(|\partial_x A(0)|+|\hat b^m|)-C^{d,r-1}$-close to a family of constant mappings (since $d\le r-2$). We conclude the proof of the first statement by noting that $(-\alpha_2(f_{a\, b\, m}),0)$ belongs to the range of $A_{a\, b}$. 
\end{proof}
\subsection{Examples of Hénon-like renormalizations useful for the proof of Theorem \ref{thm A}}\label{examren}
Let $d\le r-3$ and $m\ge 1$. In \textsection \ref{henonpiece} and Proposition \ref{hauteurbandes}, we defined $I_0$, $b_0>0$, $\delta_0>0$ and $\hat K_0$. We recall that given a $\delta$-$C^{d,r}$-Hénon-like family $(f_{p})_{ I\times J}$, with $I\subset I_0$, $\delta<\delta_0$ and $\hat b:= \max_J |b|\le b_0$,
for every $k\ge 1$, the piece $\sc_k:= \sw_= \star \ss_+\star\ss_-^{\star (k-1)}$ is well defined. 

We will continue to denote $g_p(x,y):= x^2+ a- b^m y+\zeta_p(x,b^m y)$ the first coordinate of $f_p$, with $p=(a,b)$ and $m\ge 1$ the multiplicity of $(f_p)_p$. As $d\le r-3$, a coordinate change enables us to assume that $ \partial_x \zeta_p(0)=0$ for every $p$.

If furthermore, for $j\ge 0$, it holds $I\subset I_j$, $\delta<\delta_j$ and $\hat b:= \max_J |b|\le b_j$, with $I_j, \delta_j, \hat b_j$ defined in Proposition \ref{propVp}, then we define $\sc_{k\, j}$ by: 
\[ \sc_{k\, j}:= \sc_k\star 
\bm(\sc_j-\sc_{j+1})\star \bm(\sc_0-\sc_1) \quad \text{if } b^m<0\qand \sc_{k\, j}:= \sc_k\star \bp(\sc_j-\sc_{j+1})\star \bm(\sc_0-\sc_1)\quad \text{if } b^m>0.\]

In the proof of Theorem \ref{thm A} we will do two simultaneous renormalizations. Both are obtained thanks to Theorem \ref{Charthenonlike}, one of them will be done via the piece $\sc_k$ and the other thanks to the piece $\sc_{k\, j}$ for some $j\ge 1$. Let us explain these two renormalizations in the two following examples. 
 
\begin{exam}\label{examren1}
For all $0<\hat b<b_0$ and $\delta<\delta_0$ sufficiently small, for all $\bar \delta, R>0$ and $m\ge 1$, for every $k$ sufficiently large, the renormalization given by Theorem \ref{Charthenonlike} can be applied to any $\delta$-$C^{d,r}$-Hénon-like $(f_{p})_{I_0\times J}$ of multiplicity $m$ with $J\subset [-\hat b,\hat b]$ to the piece $\sc_k$. 

The renormalization is a $\bar \delta-C^{d,r}$-Hénon-like family $(\mathcal R^- f_{p_-})_{p_-}$ 
which is defined on 
$[-R,R]\times [-R/b_- ^M, R/b_-^M]$ whenever $|a_-|\le R$, with $M=m(n_{\sc_k}+1)$ and $p_-=(a_-,b_-)$ depending on $p=(a,b)$. Also $p\mapsto b_-/b$ is $C^d-O(\delta)$-close to $1$. 
\end{exam}
\begin{proof} 
By Proposition \ref{hauteurbandes}.(2), the affine-like representation $(A_p, B_p)_p$ of $\sc_k$ displays distortion bounds $\mathcal B^{d,r}_0$, $\mathcal B^{d,r}_1$ and $\mathcal B^{d}_m$ at most $\hat K_0$, for every $k$.

To apply Theorem \ref{Charthenonlike}, we shall verify conditions $(i)-(ii)-(iii)$. The first condition is obviously satisfied. To check the two other conditions, we consider the following map:
$$\psi_p(x, c) \mapsto g_p(c+x,B_p(c+x,c))-A_p(c,c+x)\; .$$
 We compute:
 $$\partial_x \psi_p(0, c) = 2c-b^m \partial_x B_p(c,c) + \partial_x \zeta_p(c,b^m B_p(c,c))
+\partial_y\zeta_p (c, b^m B_p(c,c))b^m \partial_x B_p(c,c) - \partial_y A_p(c,c)\; .
 $$
 We recall that Proposition \ref{bsmall} implies that $\partial_y A_p$ is dominated by $b^m$. As $\|\zeta_p\|_{C^{d,r}}\le \delta$, $\partial_c \partial_x \psi_p(0, c)$ is close to $c\mapsto 2c$. Furthermore, using that $ \partial_x \zeta_p(0)=0$ for every $p$, it comes:
 \[\partial_x \zeta_p(c,b^m B_p(c,c)) = 
c \int_{0}^1 \partial_x^2 \zeta_p(t\cdot c,0)dt 
+
 b^m B_p(c,c)\int_{0}^1 \partial_y \partial_x \zeta_p(c,t \cdot b^m B_p(c,c))dt\; .\]
Then the implicit function Theorem implies:
 \begin{fact}\label{maxx0} For every $p$ there exists a unique $c(p)$ so that $\partial_x \psi_p(0, c(p))=0$. Moreover, there exists $\bar K_0>0$ independent of $k$ large so that the function $p\mapsto c(p)$ is of class $C^d$ and of norm at most $\bar K_0 |\hat b|^m$.
\end{fact} 
Put $\sigma_{k\, p}:= \partial_x A_p(c,c)(p)$ and $\lambda_{k\, p}:= \partial_y B_p(c,c)(p)$. For $\delta$ and $\hat b$ sufficiently small, for any $k$ sufficiently large, 
Corollary \ref{transpara2} sates that $p\mapsto A_p(c(p),c(p))$ is $C^d$-close $p\mapsto -\alpha_2(f_p)$ whose derivative w.r.t. $a$ is close to be in $[1/3,1/2]$ by Fact \ref{transpara3}.
As $(A_p,B_p)_p$ is $C^{d,r}$-bounded, this implies:
\begin{fact}
\label{ND4henon}
The map $p\mapsto \mu_p:=g_{p} (c(p),B_{p}(c(p),c(p)))- A_p(c(p),c(p))$ is $O(\delta \hat b^m)$-$C^{d}$-close to $p\mapsto a-b^m B_{p}(0)-A_p(0)+\zeta_p(0)$. The derivative w.r.t. $a$ is close to be in $[1/2,2/3]$.\end{fact}
Thus the implicit function theorem implies the existence of a $C^d$-function $b\in J\mapsto a (b)\in I_0$ so that for every $b\in J$, at $p=(a (b),b)$ it holds $\mu_p=0$. 

A direct computation shows that condition $(ii)$ is satisfied with $p\mapsto q_p$ $O(\hat b^m+\delta)$-close to $1$ and $(r_p)_p$ $O(\hat b^m+\delta)$-$C^{d,r}$-small. Thus the tangency is quadratic as asked in $(iii)$. Also with $z_0(p):= (A_p(c,c), c)(p)$, the function $p\mapsto f^i_p(z_0(p))$ are uniformly $C^d$-bounded among $0\le i\le n_{\sc_k}$ by Lemma \ref{prebirkof2}, and $(z\mapsto b^{-m}\det D_zf_p)_{p=(a,b)}$ is $C^{d,r}$-$O(\delta)$-close to $(z\mapsto b)_{p=(a,b)}$. This implies that $p\mapsto b_-/b:= \frac1b \sqrt[M] {\det D_{ z_0(p)} f_p^{n_{\sc_k}+1}}$ is $C^d$-close to $1$. Together with Fact \ref{ND4henon}, this implies that Condition $(iii)$ is satisfied. Put $a_-:= \mu_p/\sigma^2_{k\, p}$.We notice that $p_-\mapsto p$ has $C^d$-norm close to $1$. As $(\sigma_{k\, p})_p$ is small when $k$ is large, we can apply Theorem \ref{Charthenonlike} which implies the existence of a renormalization chart $\phi^-_p$ and a renormalized dynamics $(\mathcal R^- f_{p_-})_{p_-}$ of multiplicity $M=m(n_{\sc_k}+1)$, with $p_-= (a_-,b_-)$.

Let us show that for any $R>0$, for every $k\ge 1$ large enough, the renormalized mapping $\mathcal R^- f_{p_-}$ is well defined on $[-R,R]\times [-R/b^M_-,R/b^M_-]$ for every $b\in J$ and $a_-\in [-R,R]$.
 To this end, it suffices to show that $B_p:=\phi_p(\Delta_p)$, with $\Delta_p:= [-R,R]\times [-R/b^M_-,R/b^M_-])$, is sent by $f_p$ into $Y_{\sc_k}$ for every $p$ such that $ a_-\in [-R,R]$ and $b\in J$.
 Indeed, the set $B_p$ is a tame box containing $z_1(p)$ and with width $2\sigma_{k\, p} R$ and height $2\sigma_{k\, p} \lambda_{k\, p} R/b_-^M= 2\sigma_{k\, p} ^2 R/ \det \, D_{z_1(p)} f_p$ by Fact \ref{determinant}. 
 Then $\phi_p(\Delta_p)$ is sent by $f_p$ into 
$z_0(p)+[-K \sigma^2_{k\, p} R , K \sigma^2_{k\, p} R] \times [-\sigma_{k\, p} R,\sigma_{k\, p} R]$ for a universal constant $K$. As the $y$-coordinate $c(p)$ of $z_0(p)$ is close to $0$, it holds $f_p\circ \phi_p(\Delta_p)\subset \R\times [-3,3]$. Also by the $\mathcal B_1^{0,r}$-bounds $\sc_k$, the distance from $z_0(p)$ to $\partial^s Y_{\sc_k}$ is of the order of $\sigma_{k\, p}$, and so $f_p\circ \phi_p(\Delta_p)\subset Y_{\sc_k}$. 
\end{proof}
Put $H_p: x\mapsto B_p(c(p)+x,c(p))$. By definition \ref{canrendo}, we have:
\begin{fact}\label{shape of} When $2+\bar \delta<a_-<1/2-\bar \delta$, the canonical renormalization domain $D_p^-$ is bounded by two segments of $Graph \, H_p\pm (0,\sigma_{k\, p} ^2 / (8\det \, D_{z_1(p)} f_p))$ and two segments of a curve close to $\{(x,y): (x-c(p))^2 +det D_{z_1(p)} f_p \cdot (y-H_p(x-c(p))) = \sigma_{k\, p}^2 \beta_p \}$ for a certain $\beta_p\in (1/2,2)$. 
\end{fact}

The following is shown similarly:
\begin{exam}\label{examren2}
For every $j\ge 0$, For $0<\hat b<b_j$ and $\delta<\delta_j$ sufficiently small, for all $\bar \delta, R>0$ and $m\ge 1$, for every $k$ sufficiently large, the renormalization given by Theorem \ref{Charthenonlike} can be applied to any $\delta$-$C^{d,r}$-Hénon-like $(f_{p})_{I_j\times J}$ of multiplicity $m$ with $J\subset [-\hat b,\hat b]$ to the piece $\sc_{k\, j}$. 

The renormalization is $\bar \delta-C^{d,r}$-Hénon-like family $(\mathcal R^+ f_{p_+})_{p_+}$ which is defined on $[-R,R]\times [-R/b^{M'}_+, R/b^{M'}_+]$ whenever $|a_+|\le R$, with $M':= m(n_{\sc_{k\, j}}+1)$ and $p_+=(a_+,b_+)$. Also $p\mapsto b_+/b$ is $C^d-O(\delta)$-close to $1$. 
\end{exam}
\begin{proof}
By Proposition \ref{propVp}.(2), the affine-like representation $(A'_p, B'_p)_p$ of $\sc_{k\, j}$ displays distortion bounds $\mathcal B^{d,r}_0$, $\mathcal B^{d,r}_1$ and $\mathcal B^{d}_m$ at most $\hat K_j$, for every $j$.
Beside this, the proof is exactly the same as before. Let us just recall the main facts for these settings. As before it holds:
\begin{fact}\label{maxx1} For every $p$, there exists a unique $c'(p)$ so that 
$\partial_x g_p(x,B'_p(x,c'(p))= \partial_x A'_p(c'(p), x)$ at $x=c'(p)$. Moreover, there exists $\bar K_j>0$ so that the function $p\mapsto c'(p)$ is of class $C^d$ and of norm at most $\bar K_j |\hat b|^m$.
\end{fact} 

Also by Corollary \ref{transpara2} and Fact \ref{transpara3}, the map $p\mapsto A'_p(c',c')(p)$ is $C^d$-close to $p\mapsto \alpha_2(f_p)$ 
and its derivative is close to be in $[1/3, 1/2]$. This gives:
\begin{fact}
\label{ND4henon2}
The map $p\mapsto \mu'_p:=g_{a\, b} (z'_1(p))- A_p'(c'(p),c'(p))$ is $O(K_j\delta \hat b^m)$-$C^{d}$-close to $p\mapsto a-b^m B'_{p}(0)+A_p'(0)$. The derivative w.r.t. $a$ is close to be in $[1/2,2/3]$.\end{fact}
Thus the implicit function theorem implies the existence of a $C^d$-function $b\in J\mapsto a'(b)\in I_1$ so that for every $b\in J$, at $p=(a'(b),b)$ it holds $\mu'_p=0$. Then Theorem \ref{Charthenonlike} implies the existence of a renormalization chart $\phi^+_p$ and a renormalized dynamics $(\mathcal R^+ f_{p_+})_{p_+}$ of multiplicity $M'=m(n_{\sc_{k\, j}}+1)$, with $p_+= (a_+,b_+)$. The same proof as in the previous example shows that for every $j\ge 1$ and then $k$ large enough, the renormalized mapping $\mathcal R^+ f_{p_+}$ is well defined on $[-R,R]\times [-R/b'^M_+,R/b'^M_+]$ for every $b\in J$ and $a_+\in [-R,R]$.
\end{proof}
\begin{rema} The canonical adaptation of Fact \ref{shape of} to this renormalization holds true. 
\end{rema}

\begin{figure}[h!]
 \centering
 \includegraphics[width=17cm]{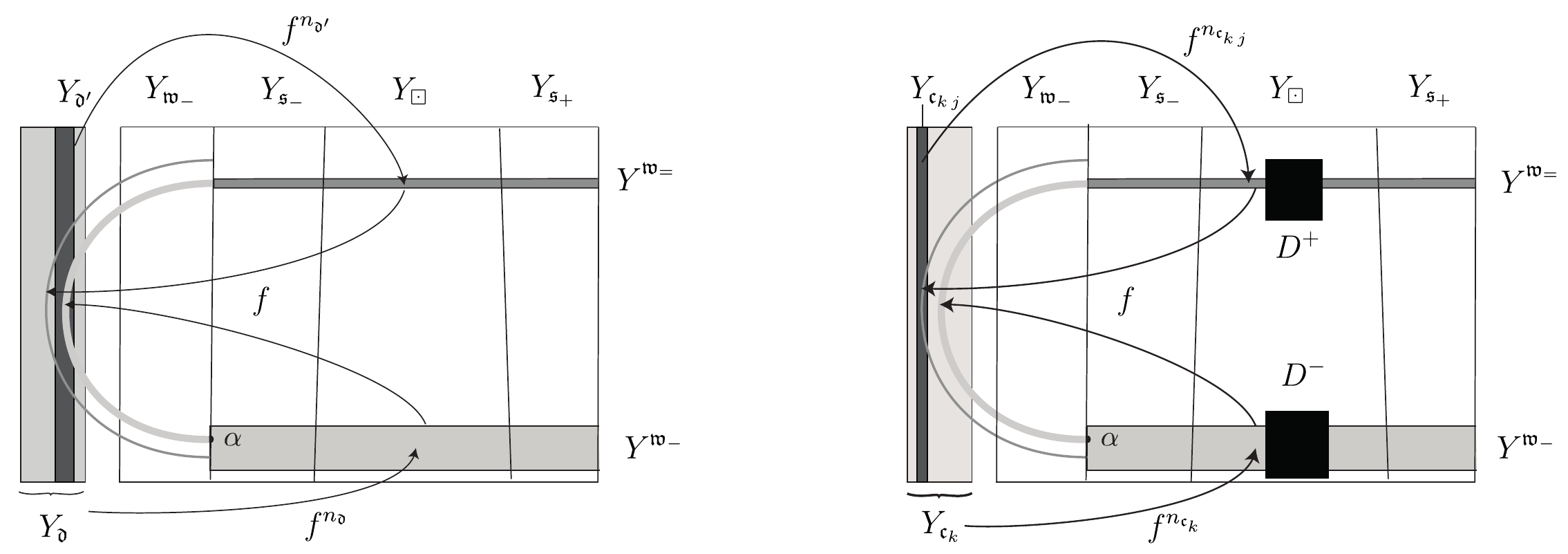}
 \caption{The figure on the left hand side depicts the setting for the multi-renormalization done to prove Theorem \ref{thm B}. The  figure on  the right hand side depicts the setting for the two renormalizations done for the proof of Theorem \ref{thm B}.}
 \end{figure}

\subsection{Proof of Theorem \ref{thm B} (stated in \textsection\ref{texsectionthmB})}
Let $0\le d\le r-3$ and let $(f_{p})_{ I\times J}$ be a $L$-wide, $\delta$-$C^{d,r}$ Hénon-like family of multiplicity $1$. Assume that $\delta<\delta_1$ and $\hat b:= \max_J |b|<b_1$, where $b_1$ and $\delta_1$ were defined in Proposition \ref{propVp}. 
We are going to prove Theorem \ref{thm B} by applying Theorem \ref{Charthenonlike2} \textsection \ref{sectionthmD} with $N=2$ and the two following pieces for a certain $k\ge 0$ large:
\[ \sd:= \sc_k= \sw_= \star \ss_+\star \ss_-^{\star (k-1)}
\qand \left\{ \begin{array}{cc} \sd':= \sc_k\star \ss_- \star \bm(\sc_0-\sc_1)& \text{ if } b>0\; ,\\
 \sd':= \sc_k\star \ss_+ \star \bm(\sc_0-\sc_1)& \text{ if } b<0\; . \end{array}\right.\]
Indeed, by Proposition \ref{defstarhenon} and Fact \ref{coneta2}, these pieces persist for every $p\in I_1\times J$ and every $k\ge 0$. 

Assumption $(i)$ of Theorem \ref{Charthenonlike2} is obviously satisfied. 
Let us show Assumptions $(ii)$ and $(iii)$. 

Let $(A_p, B_p)$ and $(A'_p, B'_p)$ be the affine-like representations of respectively $\sd$ and $\sd'$. By Proposition \ref{propVp}.(2), their distortion bounds $\mathcal B_0^{d,r}$, $\mathcal B_1^{d,r}$, $\mathcal B_1^{d}$ are at most $K_1$ for every $k\ge 0$. 
Let $g_p$ and $\zeta_p$ be such that $f_p(x,y)=(g_p(x,y),x)$ and $g_p(x,y)=x^2+a-b y +\zeta_p(x,b y)$. Put:
\[ \left\{\begin{array}{c} \Delta_p(x,y,c):= g_p(c +x,B_p(c+x,c)+y)-A'_p(c,c+x)\\
\Delta_p'(x,y,c):= g_p(c +x,B'_p(c+x,c)+y)-A_p(c,c+x)\end{array}\right.\]
With exactly the same proof as for Fact \ref{maxx0}, we show:
\begin{fact}\label{4.1B} If $\delta$ and $\hat b$ are sufficiently small, for every $p\in I_1\times J$, there exist unique real numbers $o_p $ and $o'_p $ such that $\partial_ x \Delta_p(0,o_p)=0=\partial_x \Delta'_p(0,o'_p) $. 
\end{fact}
We recall that $(o_p, B_p(o_p , o_p))\in Y^{\sd}$ is included in $Y^{\ss_-}\subset Y^{\sw_-}$ for every $k\ge 2$, whereas, $
(o'_p, B'_p(o'_p , o'_p))\in Y^{\sd'}$ is included in $Y^{\sw_=}\subset Y^{\sw_+}$. 
 By Proposition \ref{hauteurbandes}.(1), the vertical distance between $Y^{\sd}$ and $Y^{\sd'}$ is of the order of $1$. Consequently: 
\begin{fact}\label{fact4.5B} For every $p\in I_1\times J$, the value
$B'_p(o'_p , o'_p)- B(o_p , o_p)$ is of the order of $1$ and positive.
\end{fact} 

\begin{lemm}\label{coordchange}Up to a coordinate change we can assume that $o'_p=o_p=0$ for every $p$.\end{lemm}
\begin{proof} 
First by a coordinate change we can (temporary) assume that $\partial_x \zeta_p(0)=0$. Then with exactly the same proof as for Fact \ref{maxx0}, we obtain that $(o_p/b)_p$ and $(o_p'/b)_p$ are $C^d$-bounded. Thus $(f_p(o'_p, B'_p(o'_p , o'_p))-f_p(o_p, B_p(o_p , o_p)))/b= (O(b)+ B'_p(o'_p , o'_p)- B(o_p , o_p), O(1))$. Hence, there exist a universal constant $K>0$ and a $C^d$-family of functions $\pi_p$ with $C^{d,r}$-norm
 bounded by $K$ satisfying the following properties. 
 The differential $\partial_y \pi_p$ is at least $1/K$ and $\pi_p \circ f_p(o_p, B_p(o_p , o_p))= \pi_p \circ f_p(o'_p, B'_p(o'_p , o'_p))=0$.
We note that for $b$ small enough, $\psi_p:= (\pi_p\circ f_p, \pi_p)$ is a diffeomorphism with Jacobean discriminant lower bounded by a positive function of $K$. Therefore, the $\psi_p\circ f_p\circ \psi_p^{-1}$ have the requested property and form a $O(\delta)$-$C^r$-Hénon-like family.
\end{proof}
We note that Hypothesis $(iii)$ of Theorem \ref{Charthenonlike2} is satisfied with quadratic terms $q_p,q_p'$ and so functions $r_p$ and $r'_p$ satisfying:
 \[\begin{split} (q_p-1) dx^{\otimes 2}+ D^2 r_p (x,y)= 
 D^2(b B_p(x,0) +\zeta_p(x, b B_p(x,0)+by) -A'_p(0,x))\\
(q_p'-1) dx^{\otimes 2}+ D^2 r'_p (x,y)= 
 D^2(b B'_p(x,0)+\zeta_p(x, b B'_p(x,0)+by) -A_p(0,x))\end{split}\]
As $(A_p,B_p)_p$ and $(A'_p,B'_p)_p$ are $C^{d,r}$-bounded by $K_1$ and as $\partial_y A_p$ and $\partial_y A'_p$ are bounded by $\hat b$ by Proposition \ref{bsmall}, it comes 
that $(r_p)_p$ and $(r_p')_p$ are $O(\hat b+\delta)-C^{d,r}$-small.

Also we notice that $((q_p-1)/b,(q_p'-1)/b)_p$ are $O(K_j \delta)$-small. Thus, in the proof of Lemma \ref{coordchange}, by choosing a neat $\pi_p$ (one can vary the derivative of $\pi_p$ nearby $f_p(0,B_p(0))$ and $f_p(0,B'_p(0))$), we can indeed assume that $q_p=q'_p=1$ for every $p$.

Let us show Assumption $(ii)$ which states the existence of a parameter $p_0\in I_1\times J$ so that with $\mu_p:=\Delta_p(0)$ and $\mu'_p:=\Delta'_p(0)$ (in the new coordinates), it holds $\mu_{p_0}=0=\mu'_{p_0}$. In order to do so, we shall find a neat value of $k$ and we will assume $\hat b$ small in function of the wideness $L$. Note that:
\[ \mu_p:= a-b B_p(0)+\zeta_p(0, b B_p(0))-A_p'(0)\qand \mu_p':= a-b B'_p(0)+\zeta_p(0, b B'_p(0))-A_p(0)
\; .\]
By exactly the same proof as for Fact \ref{ND4henon}, we show that when $k$ is large and $\delta, \hat b$ are small, the derivative $\partial_a \mu_p$ is close to be in $[1/2, 2/3]$. Hence we have:
\begin{fact} For every $b\in J$, there exists a unique $a(b)\in I_0$ so that $\mu_{p(b)}=0$ with $p(b)=(a(b),b)$. Also $b\mapsto p(b)$ is $C^d$-bounded independently of $\hat b$ small and $k$ large.\end{fact} 
Also $(\zeta_p(0, b B'_p(0))-\zeta_p(0, b B_p(0)))_p$ is $O(\delta \hat b)-C^{d,r}$-small. Thus:
\begin{equation}
\label{difmu} \mu'_p-\mu_p= -b (B'_p(0)-B_p(0))-(A_p(0)-A'_p(0))+O(\delta \hat b)\; .\end{equation}
We recall that $(A_p,B_p)$ is the affine-like representation of $\sd= \sc_k$, and we have the freedom to chose $k$ large and $\hat b$ small independently. In view of Fact \ref{fact4.5B}, we shall define $k$ so that the distance $A_p(0)-A'_p(0)$ is of the order of $-\hat b/\sqrt L$ for every $b\in J$. To this end, we introduce the notations $\sigma_{k\, p}:= \partial_x A_p(0)$ and $ \hat \sigma_k:= \max_{b\in J} |\sigma_{k\; p(b) }| $, and we define $k$ as follows:
\begin{equation}\label{defk}
k:=k(\hat b):= \max \{i: \hat \sigma_i \le \hat b/ \sqrt L\}\; .\end{equation}
We notice that $k$ is large when $\hat b$ is small. Also remark that $ (\sigma_{k+1\, p}/\sigma_{k\, p})$ is in $[1/\|Df^2\|,\|Df^2\|]\subset [1/16, 16]$. Also by the bound $\mathcal B^{d}_1$ on $\sd= \sc_k$, 
$|\partial_b \log \sigma_{k\, p(b)}|\le K_1 n_{\sc_k}$ and since the length of $J$ is ${L}\hat b$, it comes:
$$\max_{b\in J} |\log |\sigma_{k\, p(b)}\sqrt L/\hat b||\le \log 16 + K_1 n_{\sc_k} |\hat b| L$$ 
By (\ref{defk}), the integer $n_{\sc_k}$ is dominated by $- \log (\hat b/\sqrt{ L})$,
thus we can assume $\hat b$ small in function of $L$ so that $ L |\hat b| \log |\hat b/\sqrt{ L}| $ is small. Then it comes:
\begin{fact}\label{fact4.1B} The value $\max_{b\in J} |\log |\sigma_{k\, p(b)}\sqrt L/\hat b||\le \log 17$.
\end{fact}

We recall that for $p\in \{p(b): b\in J\}$, the point $(A_p(0),0)$ is send by $f^{n_{\sc_k}}_p$ to $(0, B_p(0))$, whereas by definition of $\sd'$, the point $(A'_p(0),0)$ is send by $f^{n_{\sc_k}}_p$ into $Y_{s_\pm }$ where $\pm$ is the sign of $-b$. 
We infer that the distance from $Y_{\ss_\pm}$ to $\{0\}\times \R$ is of order of $1$, and we use the distortion bound $\mathcal B_1^{0,r}$ on $A_p$ to obtain that $(A_p(0)-A'_p(0))$ is of the order of $\sigma_{k\, p}$. Since $Q^{n_{\sc_k}}|\R_{\sc_k}$ is orientation preserving, the sign $\pm$ of $-b$ is the same as the one of $(A_p(0)-A'_p(0))$. Moreover by Fact \ref{fact4.1B}, $A_p(0)-A'_p(0)$ is of the order of $\hat b/\sqrt{L}$. 

As by Fact \ref{fact4.5B}, $-b (B_p(0)- B'_p(0))$ is of the order of $b$ and of the sign of $-b$, 
for $L$ sufficiently large (and then $\hat b$ small enough), there exists $b_0 \in [\hat b L^{-2/3} , \hat b L^{-1/3} ]$ so that with $p_0=(a_0(b_0),b_0)$, it holds $\mu_{p_0}=\mu'_{p_0}$. As $\mu_{p_0}=0$, Assumption $(iii)$ is proved. 

Hence we can apply Theorem \ref{Charthenonlike2}. Let $(\bar a_0,\bar b_0)$ and $(\bar a_1,\bar b_1)$ be the two renormalized parameters of the Hénon-like map $\mathcal R_0 f_p$ and $\mathcal R_1 f_p$. 

Observe that $\sigma_{k\, p}$ is of the order of $\sigma_{k\, p}':= \partial_x A'(0)$. Thus, both are also of the order of $\gamma_{0\, p}$ and $\gamma_{1\, p}$ (defined in the statement of Thm \ref{Charthenonlike2} \textsection \ref{sectionthmD}). Consequently, for any $R>0$, $\hat b$ small enough (or equivalently $k\ge 1$ large enough), by the same proof as for Example \ref{examren1}, the renormalized mappings $\mathcal R_0 f_{p}$ and $\mathcal R_1 f_{p}$ are well defined on respectively $[-R,R]\times [-R/\bar b_0,R/\bar b_0]$ and $[-R,R]\times [-R/\bar b_1,R/\bar b_1]$ whenever $(\bar a_0,\bar a_1)\in [-R,R]^2$. Then the following Lemma implies that we can apply Corollary \ref{preThmB} which in turn implies Theorem \ref{thm B}.
\begin{lemm}
For every $R>0$, when $\hat b$ is small, the map $p\mapsto (\bar a_0,\bar a_1)$ is a diffeomorphism from a small neighborhood of $p_0$ onto $[-R,R]^2$, and its inverse is $C^d$-close to be constant. Furthermore, the map $p\mapsto (\bar b_0, \bar b_1)$ is small.\end{lemm} 
\begin{proof}
We already saw that the same proof as for Fact \ref{ND4henon} shows that when $\delta$ and $\hat b$ are small (and so $k$ large), the derivatives $\partial_a \mu_p$ and $\partial_a \mu_p'$ are close to be in $[1/2, 2/3]$. Also $(A_p(0)-A_p'(0))_p$ is $C^{d}$-small by Proposition \ref{continuityinprod}. Thus by (\ref{difmu}), $\partial_b (\mu_p'-\mu_p)$ is equivalent to $B_p(0)-B'_p(0)$ which is of the order of $1$. 
This implies the first statement of the lemma since $\mu_p = \gamma_{0\, p}^2\cdot \bar a_0 $ and $\mu_p' = \gamma^2_{1\, p}\cdot \bar a_1 $, and $ (\gamma_{0\, p},\gamma_{1\, p})_p$ is small when $k$ is large. 

By symmetry, we only show that $\bar b_0$ is small. We recall that 
$\bar b_0= (\sigma_{0\, p}/\sigma_{1\, p}) \det D_{f^{-n_\sd}_p(z_0(p))} f_p ^{n_{\sd}+1} 
$, with $z_0= (0,B_p(0))$. Again by Proposition \ref{continuityinprod}, $p\mapsto \sigma_{0\, p}/\sigma_{1\, p}$ is $C^d$-bounded. On the other hand, we recall that $\det Df_p$ is $O(b(1+\delta)|)$-$C^d$-close to be constant. Also by the same argument as for example \ref{examren1}, $(f^j_p(z_0(p)))_{1\le j\le n_{\sd'}+1}$ is $C^d$-bounded. Thus $\det D_{f^{-n_\sd}_p(z_0(p))} f_p ^{n_{\sd}+1} $ is $C^d$-small. 
\end{proof}
\subsection{Proof of Theorem \ref{thm A} (stated in \textsection\ref{texsectionthmA})}
Let $j\ge 1$ be large. Let $(f_{p})_{ I\times J}$ be a $L$-wide, $\delta$-$C^{d,r}$ Hénon-like family of multiplicity $m$, with $\delta$ and $\hat b:= \max_J |b|$ sufficiently small so that Examples \ref{examren1} and \ref{examren2} apply for every $k$. 
We recall that in these examples, the renormalization is applied with the pieces $\sc_k:= \sw_= \star \ss_+\star\ss_-^{\star (k-1)}$ and $ \sc_{k\, j}:= \sc_k\star \sd_j$, where $\sd_j:= 
\bm(\sc_j-\sc_{j+1})\star \bm(\sc_0-\sc_1)$ {if } $b^m<0$ and 
$\sd_j:= \bp(\sc_j-\sc_{j+1})\star \bm(\sc_0-\sc_1)$ {if } $b^m>0$. 

We recall that $(A_p,B_p)$ and $(A'_p, B_p')$ denote the affine-like representations of $\sc_k$ and $\sc_{k\, j}$. The proofs of Examples \ref{examren1} and \ref{examren2} define $C^d$ families of points $(c_p)_p $ and $(c'_p)_p$, horizontal contractions $\sigma_{k\, p}:= A_p(c_p,c_p)$ and $\sigma_{k\, j\, p}:= A'_p(c'_p,c'_p)$ and parameter curves denoted by $\{p(b)=(a (b),b): b\in J\}$ and $\{p'(b)=(a' (b),b): b\in J\}$ respectively. 

We want to find a neat value of $k$ so that the two later parameter curves intersect each other at a unique point. The renormalizations associated to $\sc_k$ and $\sc_{k\, j}$ will be those claimed in statement of Theorem \ref{thm A}, for $p$ nearby $p_0$. 

We remark that for $\hat b,\delta $ sufficiently small and $k>j+2$,
the parameter $a(b)\in I_{k-2}\Subset I_j$ (these intervals are defined in definition \textsection \ref{271}), and so the piece $\sd_j$ persists for every $p\in \{p(b)=(a (b),b): b\in J\}$. 

 We are going to define $k$ in function of $\hat b$ as we did in the proof of Theorem \ref{thm B}. To this end we put
 $\hat \sigma_k:= \max_{b\in J} |\sigma_{k\; p(b) }|$ and
 $k:=k(\hat b):= \max \{i: \hat \sigma_i \le \hat b^m /\eta_j\}\; ,$ with $\eta_j$ defined in \textsection \ref{271}. 
By proceeding as in the proof of \ref{fact4.1B} (and using this time the bound $\mathcal B^{d,r}_m$ on $\sc_k$), we can assume $\hat b$ small in function of $L$ so that:
\begin{fact}\label{fact4.1} It holds $\max_{b\in J} |\log |\sigma_{k\, p(b)}\eta_j /\hat b^m ||\le \log 17$ and so $k$ is large when $\hat b$ is small.
\end{fact}
We recall that by Examples \ref{examren1} and \ref{examren2}, the renormalizations $\mathcal R^- f_{p^-}$ and $\mathcal R^+ f_{p^+}$ obtained from the pieces $\sc_k$ and $\sc_{k\, j}$ form $\bar\delta$-$C^{d,r}$-families with $\bar b$ small when $k$ is large. 
This proves the first conclusion of Theorem \ref{thm A}.
\subsubsection{Proof conclusions (2-3) of Theorem \ref{thm A}}
Let us prove:
\begin{prop}
For every $m\ge 1$, for every $j$ large, for every $\delta$ small enough in function of $j$ and every $\hat b$ small enough in function of $j$ and $\delta$, the following property holds true for any $(f_{p})_{ I\times J}$ $L$-wide, $\delta$-$C^{d,r}$ Hénon-like family with multiplicity $m$, $\hat b= \max_J |b|$ and $L\ge \eta_j^{-2}$.

There exists a unique $b_0\in [\hat b \eta_j ^{1/(2m)}, \hat b \eta_j ^{3/(2m)}]$ so that $a(b_0)=a'(b_0)$. 
\end{prop}
\begin{proof}
We go back to the notations used in Examples \ref{examren1} and \ref{examren2}. By the same proof as for Lemma \ref{coordchange} (where $b$ is replaced by $b^m$), we can assume that $c_p=0$ and $c'_p=0$ for every $p$. Then for any $b\in J$, the parameters $p(b):= (a(b), b)$ and $p'(b):=(a'(b),b)$ are implicitly defined by:
\[a=a(b)\Leftrightarrow \mu_{a\, b}=0\qand 
a=a'(b)\Leftrightarrow \mu'_{a\, b}=0\; .\]
\[\text{with}\quad 
\mu_p:= a-b^m B_p(0)+\zeta_p(0, b^m B_p(0))-A_p(0)\qand
\mu'_p:= a-b^m B'_p(0)+\zeta_p(0, b^m B'_p(0))-A'_p(0)\;. \]
Note that $(0, B_p(0))\in Y^{\sc_k}\subset Y^{\sw_-}$ for every $k\ge 1$, whereas $(0, B_p'(0))\in Y^{\sc_{k\, j}'}\subset Y^{\sw_=}$. By Proposition \ref{hauteurbandes}.(1), the vertical distance between $Y^{\sc_k}$ and $Y^{\sc_{k\, j}}$ is of the order of $1$. Consequently: 
\begin{fact}\label{fact4.5} For every $p\in I_1\times J$, the value
$B'_p(0)- B(0)$ is of the order of $1$ and positive.
\end{fact} 

Also by Proposition \ref{propVp}.(2), $\|(A_p,B_p,A_p',B'_p)\|_{C^{d,r}}$ is bounded by $\hat K_j$ (which is independent of $k$ and $\hat b$). Moreover, $(\zeta_p)_p$ is $\delta$-$C^{d,r}$-small. Consequently for the $C^{d}$-topology:
\begin{equation}
\label{difmuA} (p\mapsto \mu_p-\mu'_p)= (p\mapsto -b^m (B_p(0)-B'_p(0))-(A_p(0)-A'_p(0)))+O(\delta \hat K_j \hat b^m)\; .\end{equation}
By Proposition \ref{continuityinprod}, the function $p\mapsto (A_p(0)-A'_p(0))$ has a $C^d$-norm bounded by $\hat \sigma_k \hat K_j\asymp \hat K_j \hat b^m/\eta_j $ which is small compared to $\hat b^{m-1}$. Using 
 the bound $\|p\mapsto B_p(0)-B'_p(0)\|_{C^d}\le \hat K_j$, it comes:
\begin{fact}\label{fact4.4} The function $b\mapsto \partial_b(\mu_{p(b)}-\mu'_{p(b)})+m b^{m-1} (B_{p(b)}(0)-B_{p(b)}'(0))
$ is $O(\hat K_j \hat b^m/\eta_j)$-$C^{d-1}$-small. 
\end{fact}
Thus by Facts \ref{fact4.4} and \ref{fact4.5}, the derivative of $b\mapsto \mu_{p(b)}-\mu'_{p(b)}$ does not vanish, and so the curves $\{(a(b),b): b\in J\}$ and $\{(a'(b),b): b\in J\}$ are transverse, with at most one intersection point. 

It remains only to show the existence of an intersection point. We recall that the point $(A_p(0),0)$ is sent by $f^{n_{\sc_k}}$ to $(0, B_p(0))$, whereas by definition of $\sc_{k\, j}$, the point $(A'_p(0),0)$ is sent by $f^{n_{\sc_k}}$ into $Y_{\sd_j}$. By Proposition \ref{propVp}.(1), 
 the distance from any point in $Y_{\sd_j}$ to $\{0\}\times \R$ is of the order of $\eta_j$. Thus by the bound $\mathcal B_1^{d,r}$ of $\sc_k$, the value $A_p(0)-A_p'(0)$ is of the order of $\eta_j\sigma_{k\, p}$. Also the sign of 
$A_p(0)-A_p'(0)$ is the same as the one of $-b^m $ since $Q^{n_{\sc_k}}|\R_{\sc_k}$ is orientation preserving.

We infer Fact \ref{fact4.5} which states that $b^m (B_p(0)- B'_p(0))$ is of the order of $b^m$ and of the sign of $-b^m$, to conclude the existence of $b_0$ so that $b^m_0\asymp \hat b^m\eta_j \in [\hat b^m \eta_j^{1/2} , \hat b \eta_j^{3/2}]$ so that $\mu_{p(b_0)}=\mu'_{p(b_0)}$. 
\end{proof}

Theorem \ref{Charthenonlike} gives an explicit definition of the renormalized parameters $(a_+, b_+)=P^+(p)= (\frac{\mu'_p}{\sigma_{k\, j\, p}}, \sqrt[M' ]{\det D_{(A'_p (0),0)} f^{n_{\sc_{k\, j}+1}}}) $ and $(a_-, b_-)=P^-(p)=(\frac{\mu_p}{\sigma_{k\, p}}, \sqrt[M ]{\det D_{(A_p (0),0)} f^{n_{\sc_{k}+1}})}$, with $M':= m+m n_{\sc_{k\, j}}$ and $M:= m+mn_{\sc_{k}}$.
\begin{center}
Let $\mathcal D$ be a small neighborhood of $\{p=(a,b_0): |a_+(p)| \le 3 \}$. 
\end{center}
By definition of $\mathcal D$, conclusion (2) of Theorem \ref{thm A} is satisfied. 
We recall that $\sigma_{k\, j\, p}\asymp \eta_j \sigma_{k\; p}$ is small compared to $\sigma_{k\, p}$ by Proposition \ref{propVp}.(1). Thus for every $p=(a,b_0) \in \mathcal D$, the value $|a_-(p)|$ is bounded by $|a_-(p_0)| + 3|\partial_{a_+} a_-|\le 4 \max \sigma_{k\, j\, p}^2/\sigma_{k\, p}^2$ which is small. 
 Thus Conclusion (3) of Theorem \ref{thm A} is satisfied. 

\subsubsection{Proof of Conclusion (4) of Theorem \ref{thm A}} We continue in the settings of the latter subsection. Moreover, we fix a parameter $p\in \mathcal D$ such that $a^+(p)\in [-2+2\delta^* ,1/4-2\delta^*]$. In particular, the Hénon-like map $\mathcal R^+ f_p$ sends $Y_D$ into itself. We are going to show that the maximal invariant of the complement of the renormalization domain of $\mathcal R f^+_p$ and $\mathcal R^- f_p$ is uniformly hyperbolic. Then a classical argument implies conclusion (4) of Theorem \ref{thm B}.
 
As we will not vary the parameter $p$ anymore, we shall not display it in the indexes. We recall that $(A,B)$ is the affine-like representation of $\sc_k$ and $(A',B')$ is the affine-like representation of $\sc_{k\, j}$.

By Proposition \ref{propVp}.(1), it holds $\sigma_k \eta_j\asymp \sigma_{k\, j}$. This enables us to summary as follows the setting of the constants in the previous subsection:
\begin{equation}\label{order} \delta^*, \delta\gg \eta_j\gg \sigma_k \gg b^m\asymp \sigma_k \eta_j\asymp \sigma_{k\, j}\; .\end{equation}

Let $\mathcal K$ be the set of points which are not attracted by $\infty$. 
We recall that $|a-a_2|$ is of the order of $\sigma_k$ which is large compared to $b^m$. Thus by Corollary \ref{remasupport2}, Lebesgue a.e. point in $\mathcal K$ eventually lands in $Y^{\sw_-} \cup Y^{\sw_=}$. 
Let $\mathcal K'$ be $\mathcal K\cap ((Y^{\sw_-}\setminus D^-) \cup (Y^{\sw_=}\setminus D^+))$ where $D^-$ and $D^+$ denote the renormalization domain of $\mathcal R^-f_p$ and $\mathcal R^+ f_p$. 

\begin{fact} \label{defbeta}
There exists $e>0$, so that for every $(x,y)\in \mathcal K'\cap Y^{\sw_-}$, it holds $|x|>\sigma_k(1+e)/2$ and 
for every $(x,y)\in \mathcal K'\cap Y^{\sw_=}$, it holds $x>\sigma_{k\, j}(1+e)/2$.
\end{fact}
\begin{proof} 
We recall that by Fact \ref{shape of}, the renormalization domain of $D^-$ is vertically bounded by two translations of $Graph B(\cdot , 0)$ by $\sigma_k^2/ 4\det D f\gg b^m$. By Proposition \ref{hauteurbandes}, $Y^{\sw_-}$ is included in $Graph B(\cdot , 0) +\{0\}\times [-\hat K_0\cdot |b|^m, \hat K_0\cdot |b|^{m }]$. 
The renormalization chart sends $Y^{\sw_-}\setminus D^-$ at the complement of the stable manifold of the fixed point $\beta$ for the normalized dynamics. As $a^-\le 1/4-\delta^*$, the $x$-coordinate of the fixed $\beta$ point is $\ge 1/2+e$ for a certain $e>0$ depending only on $\delta^*$. This proves that for every $(x,y)\in Y^{\sw_-}\setminus D^-$ it holds
$|x|\ge \sigma_k (1+e)/2$. 
We deduce the second assertion similarly by using that $Y^{\sw_=}$ is included in $Graph B'(\cdot , 0) +\{0\}\times [-\hat K_0\cdot |b|^{2m}, \hat K_0\cdot |b|^{2m} ]$ and $|b|^{2m}$ is small compared to $\sigma_{k\, j}^2/\det Df$. 
\end{proof}
We remark that $\mathcal K'$ is covered by $\mathcal K_1\cup\mathcal K_2\cup\mathcal K_3$ with:
\[\mathcal K_1:= \cup_{\sa\in \sA_j} Y_{\sa}\; , 
\quad \mathcal K_2= f^{-1}(Y_{\sc_j} \setminus Y_{\sc_k})
\bigcup (Y^{w_=} \cap f^{-1}(Y_{\sc_{j}} \setminus Y_{\sc_{k\, j}}))
 \; ,
\]
\[\mathcal K_3= 
(Y^{w_-} \cap f^{-1}(Y_{\sc_{k}})\setminus D^-)
\sqcup 
(Y^{w_=} \cap f^{-1}(Y_{\sc_{k\, j}})\setminus D^+)\; .
\]
This covering induces a dynamics $T=f^\tau$ for the following return times $\tau$:
\[
\tau \colon \mathcal K'\mapsto \left\{\begin{array}{cl}
 \tau(z)=n_\sa & \text{ if } z\in \mathcal K_1 \cap Y_{\sa} , \; \sa\in \sA_j\; ,\\
 \tau(z)= 1+n_{\sc_q}& \text{ if }
z\in \mathcal K_2 \text{ and $q\le k$ is maximal s.t. }
z\in f^{-1}(Y_{\sc_q})\; ,\\
 \tau(z)= 1+n_{\sc} & \text{ if }z\in \mathcal K_3\cap f^{-1}(Y_{\sc}),\; \sc\in \{\sc_k,\sc_{k\, j}\}\text{ with $n_\sc$ maximal.}\end{array}\right. 
\] 
We consider the following cone field:
\[C'_h:= z\in Y^{\sw_-}\sqcup Y^{\sw_=}\mapsto 
\left\{\begin{array}{cc}C_{\sc_k}(z):= \{u\cdot (1,\partial_x B(x,0)) + (0,v): |v|< \eta_j^2
|u|\}&\text{ if } z\in Y^{\sw_-}\\
C_{\sc_{k\, j}}(z):=\{u\cdot (1,\partial_x B'(x,0)) + (0,v): |v|< \eta_j^2 |u|\}&\text{ if } z\in Y^{\sw_=}\end{array}\right. \]

We are going to show that every vector in $C'_h$ is uniformly expanded by $T$. This will imply that $\mathcal K= \cap_{n\in \Z} T^n(\mathcal K')$ is a uniformly horseshoe, and so this will achieve the proof of Theorem \ref{thm B} (the stable set of a uniformly hyperbolic horseshoe of a $C^2$-map has Lebesgue measure equal to zero). 

The hyperbolicity of $\mathcal K$ is the consequence of the three following propositions shown below:

\begin{prop}\label{cone1HLrenor} There exists $\beta>0$ such that for every $z\in \mathcal K_3$, for every vector $w_0\in C_h'(z)$:
\begin{itemize}
\item the vector $w_1:=D_z T (w_0)$ belongs to $C_h'(T(z))$ and satisfies $\|w_1\|\ge (1+\beta)\|w_0\|$.
\item if $T(z)\notin \mathcal K_3$, then $\|w_1\|\ge \|w_0\|/\eta_j$. 
\end{itemize}
\end{prop}
\begin{prop}\label{cone1HLrenor2} There exists $K$ independent of $j$ such that for every $z\in \mathcal K_2$, for every vector $w_0\in C_h'(z)$, the vector $w_1:=D_z T (w_0)$ belongs to $C_{h}'(T(z))$ and satisfies $\|w_1\|\ge \|w_0\|/(K\eta_j)$.

\end{prop}

\begin{prop}\label{cone1HLrenor3} 
There exist constants $C>0$ and $\kappa<1$ 
independent of $j$ so that for every $z\in \mathcal K_1\cap T^{-1}(\mathcal K_1)\cdots \cap T^{-q}(\mathcal K_1)$, for every vector $w_0\in C_h'(z)$, the vector $w_1:=D_z T^q (w_0)$ belongs to $C_{h}'(T(z))$ and satisfies $\|w_1\|\ge C \kappa^{-N}\|w_0\|$, with $N:=\tau(z)+\cdots +\tau(T^{q-1}(z)) $.
\end{prop}
Indeed, the maximal invariant of $\mathcal K_1$ is uniformly hyperbolic, and the first return times of $T$ in $\mathcal K_2\cup \mathcal K_3$ is uniformly hyperbolic.

\begin{proof}[Proof of Proposition \ref{cone1HLrenor3}]
We recall that we defined in \textsection \ref{defCh}, the cone $C_h:= \{(u,v)\in \R^2: |v|\le c_h |u|\}$, and in \textsection \ref{henonpiece} we fixed $c_h:= 1/\eta_j$.

\begin{fact}\label{diffdeB} For every $z\in Y^{\sw_-}\cup Y^{\sw_+}$, the cone $C_{h}'(z)$ is included in $C_h$. 
Also, for every $z\in Y_{\sw_-}\cup Y_{\sw_+}$, it holds 
$D_zf(C_h)\subset C_{h}'$. 
\end{fact}
\begin{proof} 
Let $(A'',B'')$ be the affine-like representation of $\sw_-$ (resp. $\sw_+$). By Proposition \ref{bsmall}, $\partial_y B''$ is $O(b^m)$-small. Thus by Proposition \ref{continuityinprod}.(2), $B$ is $O(b^m)$-$C^d$-close 
(resp. $B'$ is $O(K_jb^m)$-$C^d$-close)
to $B''$. As $B''$ is bounded, $\partial_x B$ and $\partial_x B'$ are bounded independently of $j$. Thus $C_{h}'$ is included in $C_h$ because the angle $c_h=1/\eta_j$ is large. Also for $z\in Y_{\sw_-}$ (resp. $\in Y_{\sw_+}$), $D_zf(C_h)$ is $O(b^m/\eta_j)$-close to $(1,\partial_x B'')$, and so to $(1,\partial_x B)$ (resp. $(1,\partial_x B')$). As $b^m$ is small compare to the angle $\eta_j^2$ of $C_h'$, it comes $D_zf(C_h)\subset C_{h}'$.
\end{proof}

We recall that for every $\sd\in \sA_j$, for every $z\in Y_\sd$ is send by $f^{n_{\sd}-1}(z)$ into $Y_{\sw_+} \sqcup Y_{\sw_-}$. Also for every $q<
n_{\sd}$, the point $f^q(z)$ is $\eta_j$-distant to $\{0\}\times \R$ and so $D_z f^{n_{\sd}-1}(C_h)$ is included in $C_h$ by Lemma \ref{coneta}. Thus the latter fact gives:
\begin{fact}\label{pourcone} For every $\sd \in \sA_j$, for every $z\in Y_\sd$, $D_zT= D_zf^{n_\sd}$ sends $C'_h$ into $C_{h}'$.
\end{fact}
Then the Proposition is an immediate consequence of inequality (\ref{hypQQ}) of Lemma \ref{hypQ}, whenever $\hat b$ sufficiently small compared to $\eta_j$. 
\end{proof}

\begin{proof}[Proof of Proposition \ref{cone1HLrenor}] We use the notations of the proof of Theorem \ref{Charthenonlike}, where $H$, $V$, $\sigma$, $r$ were defined. In particular, we continue to assume that $c=0$, $q=1$. 
Put $z=(x,y+H(x))$. By Fact \ref{defbeta}, $|x|\ge |\sigma|\cdot (1+ e)/2$,  for some $e>0$.
We recall that $w_0\in C_h'(z)$ and $w_1:= D_zf^{n_\sc+1} (w_0)$. Put 
 $w_0=: (1,D H(x)) u_0+(0,v_0)$ and 
 $w_1=: (1,D H(x))u_1 +(0,v_1)$. Note that $|v_0|\le \eta_j^2| u_0|$. 
By Equation (\ref{pourmainthm}) of the proof Theorem \ref{Charthenonlike}, it holds:
\[\partial_x \breve A\cdot u_1 +\partial_y \breve A\cdot u_0= 2 x \cdot u_0 - d \cdot v_0 +Dr (u_0, d \cdot v_0)\qand \partial_x \breve B \cdot u_1 +\partial_y \breve B \cdot u_0 = v_1\; .\]
We recall that $d\sim b^m$ and $\sigma(1+e)/2\le |x| =O(\eta_j)$. 
We infer that $\sigma\gg \eta_j^2 b^m$ (see (\ref{order})), and so $|(\partial_y r-1)\cdot d \cdot v_0|\le 2| b^m \eta_j^2 u_0| \ll |x|\cdot |u_0|$. Also $|\partial_x r \cdot u_0| = O(|x|^2+|b^m y|)) \cdot |u_0|$, and $y$ is at most of the order of $b^m$ and $|x|$ is large compare to $x^2$ and $b^m y$. Thus $|\partial_x r \cdot u_0|\ll |x|\cdot |u_0|$ and:
\begin{equation}\label{u0u1} |\partial_x \breve A\cdot u_1 +\partial_y \breve A\cdot u_0|=
(2-o(1))\cdot |x|\cdot |u_0|.\end{equation}

By Lemma \ref{distorsionutil0}, there exists $K'$ depending only on the $\mathcal B^{0,r}_1$-bound of the piece $\sc$ such that:
 \begin{equation}\label{estimenoncroise}
 |\partial_x \breve A | \ge |\sigma|(1 - K' |x|)
\quad , \quad 
 |\partial_y \breve A | \le |\sigma| K' |x| 
\qand |D\breve B| \le K' |\lambda|\; .
 \end{equation}

 Thus (\ref{u0u1}) , (\ref{estimenoncroise}) and $|x|\in [\sigma(1+e)/2, \eta_j]$ give:
 \begin{equation} \label{pourlatresfin}
 |u_1|\ge (2-o(1))\frac{|x|}{|\sigma|} |u_0|
 \qand 
|v_1|\le K|\lambda| (|u_1| +|u_0|)\le 2 K |\lambda| |u_1| \; .\end{equation}
This implies the first statement of the Proposition since $2K |\lambda| $ is smaller than $\eta_j^2$ and $ |u_1|\ge |u_0|(1+e-o(1))$ (we take $\beta=e/2$).
Furthermore, if $f^{n_\sc+1}(z)$ does not belong to $f^{-1}(Y_\sc)$, then $x$ is of the order of $\sqrt{|\sigma|}$. Then $|u_1|/|u_0| $ is of the order of $1/\sqrt{|\sigma|}$ which is large compared to $1/\eta_j$. \end{proof}

\begin{proof}[Proof of Proposition \ref{cone1HLrenor2}]
The proof is similar to the one of the latter proposition. Again we will not display the parameter $p$ in the indexes. This time we use the affine-like representation $(A_\sq, B_\sq)$ of the piece $\sc_q$.
 Let $\sc= \sc_k $ if $z\in Y^{\sw_-}$ and $\sc= \sc_{k\, j}$ if $z\in Y^{\sw_+}$. We recall that $(A,B)$ denotes the affine-like representation of $\sc$; it defines $H$, $V$ and $\sigma$. 
 
Similarly to Example \ref{examren1}, there exist $\mu_\sq$, $q_\sq$ close to $1$, $d_\sq$ of the order of $b^m$, $c=O(b^m)$ and a function $r_\sq$ with $r_\sq(0)=Dr_\sq= \partial_x^2r_\sq(0)=0$ so that:
 \[ g(c+x,H (c+ x)+y)-V_\sq(c+x)= \mu_\sq + q x^2-d_\sq y+r_\sq (x,d_\sq y)\; .\]
 
 We write $z=(c+x,y+H(c+x))$. It belongs to $ Y^\sw \cap f^{-1}(Y_{\sc_q}\setminus Y_{\sc'_{q}})$ with $\sc'_{q}=\sc_{q+1}$ and $\sw\in \{\sw_-,\sw_=\}$ if $q<k$ and with $\sc'_{q}=\sc_{k\, j}$ and $\sw=\sw_=$ otherwise.

 \begin{lemm}\label{dist to 0} There exists $K>0$ independent of $j$ and $\hat b$ so that 
$ |x|\ge \sqrt{|\sigma_\sq|}/K$ if $q<k$ and $ |x|\ge \sqrt{\eta_j |\sigma_k|}/K$ if $q=k$.
\end{lemm}
\begin{proof}
 
If $q<k$, by the bound $\mathcal B_1^{0,r}$, the distance between the curve $\mathcal V$ and $Y_{\sc_q}\setminus Y_{\sc_{q+1}}$ is at least of the order of $\sigma_q$. If $q=k$, the distance between the curve $\mathcal V$ and $Y_{\sc_k}\setminus Y_{\sc_{k\, j}}$ is at least of the order of $\sigma= \sigma_k\eta_j $.

As the renormalization parameters $a_-, a_+$ are in $[-3,3]$, 
 the point $(0, H(0))$ is sent by $f$ into $Y_{\sc}$ and even $O(\sigma^2)$-close to be tangent to $\mathcal V$. Hence $g(0, H(0))$ is $o(\sigma)$-close to $V(0)$. As $c=O(b^m)$ and $\partial_y V=O(b^m)$, it comes that $g(c, H(c))$ is $O(b^{2m})\le o(\sigma)$-close to $V(c)$.

Let $\bar z:= f(z)=(g(x,H (x)+y),x)$. We have $\bar z\in Y_{\sc_q}\setminus Y_{\sc_{q}'}$, and so $g(c+x,H_\sq (c+x)+y)-g(c,H_\sq (c))$ is at least of the order of $\sigma_q\ge \sigma$ if $q<k$, and at least of the order of $ \sigma\asymp \sigma_k \eta_j$ if $q=k$. In other words, the value $q x^2-d_\sq y+r_\sq (x, d_\sq y)$ is at least of the order of $\sigma_q$ if $q<k$ and $\sigma_k \eta_j$ if $q=k$.
If $q<k$, then $|y|\le 3$ and $d_q |y|$ is dominated by $b^m=o(\sigma)$. If $q=k$ then $|y|=O(b^m)$ and $d_q |y|$ is dominated by $b^{2m}=o(\sigma)$. The three latter sentences imply the Lemma.
 \end{proof}
 We recall that $w_0\in C_h'(z)$ and $w_1:= D_zf^{n_{\sc_q}+1} (w_0)$. Put $w_0=: u_0\cdot (1,D H(x))+(0,v_0)$ and $w_1=: u_1\cdot (1,D H_\sq(x))+(0,v_1)$. By definition of $C_h'(z)$, $|v_0|\le \eta_j^2 |u_0|$. In the proof of Fact \ref{diffdeB}, we saw that the function $B_\sq$ and $B$ are $O(b^m)$-close to $B''$. Hence, to prove the cone condition it suffices to show that $|v_1|\le \eta_j^2 |u_1|/2$.

Put $\breve A_\sq= A_\sq -V_\sq$ and $\breve B_\sq = B_\sq-H_\sq$. Then it holds:
\[\partial_x \breve A_\sq\cdot u_1 +\partial_y \breve A_\sq\cdot u_0= 2 x \cdot u_0 - d_\sq \cdot v_0 +Dr_\sq (u_0, d_\sq \cdot v_0)\qand \partial_x \breve B_\sq \cdot u_1 +\partial_y \breve B_\sq \cdot u_0 = v_1\; .\]
Similarly to the proof of (\ref{pourlatresfin}), we have:
 \begin{equation} \label{pourlatresfin2}
 |u_1|\ge (2-o(1))\frac{|x|}{|\sigma_\sq|} |u_0|
 \qand 
|v_1|= O(|\lambda_\sq| |u_1|)= o(\eta_j^2 |u_1|)\; .\end{equation}
The second inequality implies the cone condition. The first inequality and Lemma \ref{dist to 0} imply that the expansion of the vectors in the cone by $DT$ is at least $\min(|\sigma_\sq|^{-1/2}, |\sigma_{\sc_k}|^{-1/2}\eta_j)/K$ which is large because $q\ge j$ is large and $k$ is large compared to $j$. 
\end{proof}

\medskip 
\subsubsection*{
I thanks A. de Carvalho and D. Turaev for important conversations. }
\bibliographystyle{alpha}

%\nocite{*}
\bibliography{references}

\newcommand{\etalchar}[1]{$^{#1}$}
\def\cprime{$'$} \def\cprime{$'$} \def\cprime{$'$}
\begin{thebibliography}{DCLM05}

\bibitem[BC91]{BC2}
M.~Benedicks and L.~Carleson.
\newblock The dynamics of the {H}\'enon map.
\newblock {\em Ann. Math.}, 133:73--169, 1991.

\bibitem[BDS16]{BedeSi}
P.~Berger and J.~De~Simoi.
\newblock On the {H}ausdorff dimension of {N}ewhouse phenomena.
\newblock {\em Ann. Henri Poincar\'e}, 17(1):227--249, 2016.

\bibitem[Ber11]{berhen}
P.~Berger.
\newblock Abundance of one dimensional non uniformly hyperbolic attractors for
  surface endomorphisms.
\newblock {\em arXiv:0903.1473v2}, 2011.

\bibitem[BT17]{BT17}
Pierre Berger and Dimitry Turaev.
\newblock On herman's positive entropy conjecture.
\newblock {\em arXiv preprint arXiv:1704.02473}, 2017.

\bibitem[BY14]{BY}
P.~Berger and J.-C. Yoccoz.
\newblock Strong regurlarity.
\newblock pages 1--13, 2014.

\bibitem[DCLM05]{CLM}
A.~De~Carvalho, M.~Lyubich, and M.~Martens.
\newblock Renormalization in the {H}\'enon family. {I}. {U}niversality but
  non-rigidity.
\newblock {\em J. Stat. Phys.}, 121(5-6):611--669, 2005.

\bibitem[DGS{\etalchar{+}}13]{DGLS13}
A.~Delshams, S.~V. Gonchenko, Gonchenko~V. S., J.~T. L\'azaro, and O.~Sten'kin.
\newblock Abundance of attracting, repelling and elliptic periodic orbits in
  two-dimensional reversible maps.
\newblock {\em Nonlinearity}, 26(1):1, 2013.

\bibitem[EM82]{EHM82}
H.~{El Hamouly} and C.~{Mira}.
\newblock {Singularites dues au feuilletage du plan des bifurcations d'un
  diff\'eomorphisme bi-dimensionnel.}
\newblock {\em {C. R. Acad. Sci., Paris, S\'er. I}}, 294:387--390, 1982.

\bibitem[GST93]{GST93}
S.~V. Gonchenko, L.~P. Shil{\cprime}nikov, and D.~V. Turaev.
\newblock On the existence of {N}ewhouse domain in a neighborhood of systems
  with a structurally unstable {P}oincaré homoclinic curve (the
  higher-dimansional case).
\newblock {\em Russian Acad. Sci. Dokl. Math.}, 47(2):268--273, 1993.

\bibitem[GST08]{GST08}
S.~V. Gonchenko, L.~P. Shilnikov, and D.~V. Turaev.
\newblock On dynamical properties of multidimensional diffeomorphisms from
  {N}ewhouse regions. {I}.
\newblock {\em Nonlinearity}, 21(5):923--972, 2008.

\bibitem[Haz11]{Ha11}
P.~E. Hazard.
\newblock H\'enon-like maps with arbitrary stationary combinatorics.
\newblock {\em Ergodic Theory Dynam. Systems}, 31(5):1391--1443, 2011.

\bibitem[Mil92]{Mil92}
J.~Milnor.
\newblock Remarks on iterated cubic maps.
\newblock {\em Experiment. Math.}, 1(1):5--24, 1992.

\bibitem[New74]{Newhouse}
S.~E. Newhouse.
\newblock Diffeomorphisms with infinitely many sinks.
\newblock {\em Topology}, 12:9--18, 1974.

\bibitem[PT93]{PT93}
J.~Palis and F.~Takens.
\newblock {\em Hyperbolicity and sensitive chaotic dynamics at homoclinic
  bifurcations}, volume~35 of {\em Cambridge Studies in Advanced Mathematics}.
\newblock Cambridge University Press, 1993.

\bibitem[PY01]{PY01}
J.~Palis and J.-C. Yoccoz.
\newblock Implicit formalism for affine-like maps and parabolic composition.
\newblock In {\em Global analysis of dynamical systems}, pages 67--87. Inst.
  Phys., Bristol, 2001.

\bibitem[PY09]{PY09}
J.~Palis and J.-C. Yoccoz.
\newblock Non-uniformly hyperbolic horseshoes arising from bifurcations of
  {P}oincar\'e heteroclinic cycles.
\newblock {\em Publ. Math. Inst. Hautes \'Etudes Sci.}, (110):1--217, 2009.

\bibitem[Rob83]{Ro83}
C.~Robinson.
\newblock Bifurcation to infinitely many sinks.
\newblock {\em Comm. Math. Phys.}, 90(3):433--459, 1983.

\bibitem[TLY86]{TLY}
L.~Tedeschini-Lalli and J.~Yorke.
\newblock How often do simple dynamical processes have infinitely many
  coexisting sinks?
\newblock {\em Comm. Math. Phys.}, 106(4):635--657, 1986.

\bibitem[Tur15]{Tu15}
D.~Turaev.
\newblock Maps close to identity and universal maps in the {N}ewhouse domain.
\newblock {\em Comm. Math. Phys.}, 335(3):1235--1277, 2015.

\bibitem[vS10]{VS}
S.~van Strien.
\newblock One-dimensional dynamics in the new millennium.
\newblock {\em Discrete Contin. Dyn. Syst.}, 27(2):557--588, 2010.

\bibitem[vSn67]{Sh67}
L.~P. \v~Sil\cprime~nikov.
\newblock On a problem of {P}oincar\'e-{B}irkhoff.
\newblock {\em Mat. Sb. (N.S.)}, 74 (116):378--397, 1967.

\bibitem[Yoc97]{Y97}
J.-C. Yoccoz.
\newblock Jackobson' theorem.
\newblock {\em Manuscript}, 1997.

\end{thebibliography}

\end{otherlanguage}

\end{document}